\documentclass[12pt]{article}

\usepackage{mathrsfs}
\usepackage{amsmath}
\usepackage{amsfonts}
\usepackage{amssymb}
\usepackage[all,cmtip]{xy}

\usepackage{ulem}
\usepackage{eufrak}

\usepackage{amsmath, amsthm, amsfonts, amssymb}
\def\<{\langle}
\def\>{\rangle}

\def\ra{\rightarrow}

\def\a{\alpha}
\def\wh{\widehat}
\def\wt{\widetilde}

\def \sm{\setminus}
\def\-{\overline}

\def\o{\omega}
\def\ov{\overline}

\def\endpf{\hbox{\vrule height1.5ex width.5em}}

\def\a{\alpha}
\def\endpf{\hbox{\vrule height1.5ex width.5em}}

\def\BB{{\mathbb B}}

\def\-{\overline}

\def\o{\omega}

\def\sm{\setminus}

\def\wt{\widetilde}
\def\ra{\rightarrow}

\def\endpf{\hbox{\vrule height1.5ex width.5em}}

\def\a{\alpha}

\def\Ol{\overline}

\def\a{\alpha}

\def\endpf{\hbox{\vrule height1.5ex width.5em}}
\def\a{\alpha}

\def\BB{{\mathbb B}}

\def\ra{\rightarrow}

\def\wt{\widetilde}
\def\-{\overline}

\newtheorem{theorem}{Theorem}[section]
\newtheorem{lemma}[theorem]{Lemma}

\newtheorem{proposition}[theorem]{Proposition}

\newtheorem{definition}[theorem]{Definition}

\newtheorem{remark}[theorem]{Remark}
\newtheorem{claim}[theorem]{Claim}
\date{}

\begin{document}

\title{\bf  Volume-preserving  mappings between Hermitian symmetric spaces of compact
type}

\author{Hanlong Fang, Xiaojun Huang \footnote{Supported in part by  NSF grant DMS-1665412} \ and Ming Xiao\footnote{Supported in part by NSF grant DMS-1800549 }}

\vspace{3cm} \maketitle

\tableofcontents

\bigskip\bigskip
\section{Introduction}
Let $M$ be an irreducible $n-$dimensional Hermitian symmetric space
of compact type, equipped with a canonical K\"ahler-Einstein metric
$\omega$. Write $\omega ^n$  for the associated volume form (up to a
positive constant depending only on $n$). The purpose of this paper
is to prove  the following rigidity theorem:

\begin{theorem}\label {mainthm} Let  $(M, \omega)$ be an irreducible $n-$dimensional Hermitian symmetric space
of compact type as above. 
Let $F=(F_{1},...,F_{m})$ be a holomorphic mapping from a connected
open subset $U \subset M$ into  the $m$-Cartesian product $M \times
...\times M$ of $M$. Assume that each $F_{j}$ is generically
non-degenerate in the sense that $F^*_j(\omega^n)\not\equiv 0$ over
$U$.
Assume that $F$ satisfies the following volume-preserving (or
measure-preserving) equation:
\begin{equation}\label{eq1}
 \omega^n=\sum_{i=1}^m \lambda_i F_{i}^*(\omega^n),
\end{equation}
for certain constants $\lambda_j > 0.$ Then for each $j$ with $1\le
j\le m$, $F_{j}$ extends to a holomorphic  isometry of $(M,\omega)$.
In particular, the conformal factors satisfy the identity:
$\sum_{j=1}^{m}\lambda_j=1.$
\end{theorem}

Rigidity properties are among the  fundamental phenomena in Complex
Analysis and Geometry of several variables, that study the global
extension and  uniqueness for various holomorphic objects up to
certain group actions.  The rigidity problem that we consider in
this paper was initiated by a celebrated paper of Calabi \cite{Ca}.
In \cite{Ca},  Calabi studied the global holomorphic extension and
uniqueness (up to the action of the holomorphic isometric group of
the target space) for a local holomorphic isometric embedding from a
K\"ahler manifold into a complex space form.
He established the global extension and the Bonnet type rigidity
theorem for a local holomorphic isometric embedding from a complex
manifold with a real analytic K\"{a}hler metric into a standard
complex space form. The phenomenon discovered by  Calabi \cite{Ca}
has been further explored in the past several decades due to its
extensive connection with  problems in Analysis and Geometry. (See
\cite{U} \cite{DL} \cite{DL1}, for instance).

In 2004, motivated by the modularity problem  of the  algebraic
correspondences  in algebraic number theory, Clozel and Ullmo
\cite{CU} were led to study the  rigidity problems for  local
holomorphic isometric maps and even much more general
volume-preserving maps between bounded symmetric domains equipped
with their Bergman metrics. By reducing the modularity problem to
the rigidity problem for local holomorphic isometries, Clozel-Ullmo
proved that an algebraic correspondence in the quotient of a bounded
symmetric domain preserving the Bergman metric has to be a modular
correspondence in the case of the unit disc in the complex plane and
in the case of  bounded symmetric domains of rank $\geq 2$. Notice
that in the one dimensional setting, volume preserving maps are
identical to the metric preserving maps. Thus the Clozel-Ullmo
result also applies to the volume preserving algebraic
correspondences in the lowest dimensional case. Motivated by  the
work in \cite{CU},
 Mok  carried out a systematic study of the rigidity   problem
for local isometric embeddings  in a very general setting.
Mok in [Mo2-4] proved the total geodesy for a local holomorphic
isometric embedding between bounded symmetric domains $D$ and
$\Omega$ when either (i) the rank of each irreducible component of
$D$ is at least two  or (ii) $D=\BB^n$ and $\Omega=(\BB^n)^p$ for
$n\geq 2$. In a paper of Yuan-Zhang  \cite{YZ}, the total geodesy is
obtained in the case of $D=\BB^n$ and $\Omega=\BB^{N_1} \times\cdots
\times \BB^{N_p}$ with $n\geq 2$ and $N_l$ arbitrary for $1 \leq l
\leq p$. Earlier, Ng in \cite{Ng2} had established a similar result when
$p=2$ and $2 \le n\le N_1, N_2\le 2n-1.$ In a  paper of Yuan and the
second author of this paper \cite{HY1}, we established the rigidity
result for local holomorphic isometric embeddings from a Hermitian
symmetric space of compact type into the product of Hermitian
symmetric spaces of compact type with even negative conformal
factors where certain non-cancellation property for the conformal
factors holds. (This cancellation  condition turns out be the
necessary and sufficient condition for the rigidity to hold due to
the presence of negative conformal factors.) In a recent paper of
Ebenfelt \cite{E}, a certain classification, as well as its
connection with  problems in CR geometry, has been studied for local
isometric maps when the cancellation property fails to hold. The
recent paper of Yuan \cite{Y} studied the rigidity problem for local
holomorphic maps preserving the $(p,p)$-forms between Hermitian
symmetric spaces of non-compact type. At this point, we should also
mention  other  related studies for the rigidity  of holomorphic
mappings. Here, we  quote the papers by Chan-Xiao-Yuan \cite{CXY}, Dinh-Sibony \cite{DS}, Huang [Hu1-2], Ji \cite{Ji},
Kim-Zaitsev \cite{KZ}, Mok \cite{Mo1}\cite{Mo5}, Mok-Ng\cite{MN1}, Ng [Ng1-Ng3], Xiao-Yuan [XY1-2] and many
references therein, to name a few.

The work of Clozel and Ullmo has left open an important question of
understanding the modularity problem for
volume-preserving correspondences in the quotient of
 Hermitian symmetric spaces of
higher dimension equipped with their Bergman metrics. In 2012, Mok
and Ng answered, in the affirmative, the question of Clozel and
Ullmo in \cite{MN} by establishing the rigidity  property for local
holomorphic volume preserving maps from  an irreducible Hermitian
manifold of non-compact type into its Cartesian products.

The present paper continues the above mentioned investigations,
especially those in \cite{CU}, \cite{MN} and \cite{HY1}. Our main
purpose is to establish the Clozel-Ullmo and Mok-Ng results for
local measure preserving maps between Hermitian symmetric spaces of
compact type. Notice that in the   Riemann sphere setting, Theorem
\ref{mainthm} also follows from the isometric rigidity result
obtained in an earlier paper of the second author with Yuan
\cite{HY1}. However,  the basic  approach in this paper
fundamentally  differs from that in \cite{HY1}. The method used in
\cite{HY1} is to first obtain the result in the simplest projective
space setting and then use the minimal rational curves to reduce the
general case to the much simpler projective space case. On the other
hand, restrictions of volume preserving maps are no longer volume
preserving and thus the reduction method in \cite{HY1} can not be
applied here. The approach we use in this paper is  first to
establish general results under certain geometric and analytic
assumptions (i.e., Proposition (I)-(III)) and then verify that these assumptions are automatically
satisfied based on a case by case argument in terms of the type of
the Hermitian space.

\medskip
We now briefly describe  the organization  of the paper and the
basic ideas for the proof of Theorem \ref{mainthm}.
 The major part of the paper is devoted   to showing the  algebraicity for a certain  component $F_j$  in Theorem \ref{mainthm}
 with  total degree depending only on the geometry of $(M,\omega)$. For
this, we introduce the concept of Segre family for an embedded
projective subvariety. Notice that in  the previous work, Segre
varieties were only defined for a real submanifold in a complex
space through complexification. Our Segre family is defined by
slicing the minimal embedding  with a hyperplane in the ambient
projective space, associated with  points in its conjugate space.
The Segre family thus defined is  invariant under holomorphic
isometric transformations, whose defining function is closely
related to the complexification of the  potential function of the
canonical metric. The first step in our proof is to show that a
certain component $F_j$ preserves at least locally the Segre family.
The next difficult step is then to show that preservation of the
Segre foliation  gives the algebraicity of $F_j$. To obtain the
algebraicity of $F_j$, we need to study the size that the space of
the jets of the map $F_j$ along the Segre variety  directions.
Indeed, an important  part of the paper is to  show that the space
of the jets of an associated embedding map $r_F$ along the Segre
direction up to a certain order depending only on $M$ and its
minimal embedding spans the whole target tangent space. This is a
main reason we need to describe precisely what the minimal embedding
is for each $M$.
Once this is done, we  can then show  that the map, when restricted
to each Segre variety, stays in the field generated by rational
functions and the differentiations of their defining functions as
well as their inverse, and  thus must be algebraic by a modified
version of the Hurwitz theorem.
The uniform bound of the total degree of $F_j$ is obtained by the
fact that we need only a fixed number of steps to perform  algebraic
and differential operations to reproduce the map from the minimal
embedding functions. After obtaining the algebraicity, we further
show that $F_j$ extends to a birational self-map of the space by a
monodromy argument, the geometry of the Segre foliation, an
iteration argument and the classical Bezout theorem. Finally, a
simple argument shows that a birational map which preserves the
Segre foliation  is the restriction of  a holomorphic self-isometry
of the space. Once $F_j$ is proved to be an isometry,  we can delete
$F_j$ from the original equation and then apply an induction
argument to conclude the rigidity for other components.



 The organization of the paper is
as follows: In $\S 2$, we first introduce the Segre family for a
polarized projective variety. We then describe the canonical and
minimal embedding of the space into a complex projective space in
terms of the type of the space.  In $\S 3$, we derive a general
theorem for partially degenerate holomorphic embeddings which will
play a fundamental role in the later development. In $\S 4$, we
provide the algebraicity  for one of the components of the
holomorphic mapping $F$ under additional assumptions which include
the partial non-degeneracy condition introduced in $\S 3$, the
generic transversality of the Segre varieties and the irreducibility
of the Segre family. In $\S 5$, we show that the partial
non-degeneracy holds for  local biholomorphisms between any
irreducible Hermitian space of compact type. $\S 6$ is devoted to
proving the generic transversality for the intersection of  the
Segre varieties. We prove in $\S 7$ the irreducibility of the
potential functions pulled back to a complex Euclidean space, which
has consequences on the irreducibility of the Segre varieties and
the Segre families.
 The argument in $\S 5$-$\S 7$ varies as  the type of the
space  varies and thus has to be done case by case.

We  include several  Appendices for convenience of  the reader.
 In Appendix
I, we give the concrete   functions for a minimal holomorphic
embedding of a Hermitian symmetric space of exceptional type into a
projective space. In Appendix II, we continue to establish
Proposition (I) for the rest cases. In Appendix III, we provide the
verification on the transversality for the Segre varieties for the
remaining cases not covered in $\S 6$.

\medskip
{\bf Acknowledgement}: The authors  would like to thank A. Buch, J.
Lu, L. Manivel, X. Yang and Z. Zhang for many discussions during the
preparation of this work. In particular, the first author would like
to express his gratitude to R. Bryant for answering many of his
questions on Hermitian symmetric spaces through the mathoverflow
website.

\section{Irreducible Hermitian symmetric spaces and their Segre
varieties}

\subsection{Segre varieties of projective subvarieties}
Write $z=(z_1,\cdots,z_n,z_{n+1})$ for the coordinates of ${\mathbb
C}^{n+1}$ and $[z]=[z_1,\cdots,z_n,z_{n+1}]$ for the homogeneous
coordinates of $\mathbb{CP}^{n}$. For a polynomial $p(z)$, we define
$\overline{p}(z):=\overline{p(\overline{z})}$. For a connected
projective variety $V\subset \mathbb{CP}^{n}$, write $\mathcal{I}_V$
for  the ideal consisting of homogeneous polynomials in $z$ that
vanish on $V$. We define the conjugate variety ${V}^*$ of $V$ to be
the projective variety defined by ${\mathcal I}^*_V:=\{\bar{f}: f\in
{\mathcal I}_V\}$. Apparently the map $z\mapsto \Ol{z}$ defines a
diffeomorphism from $V$ to $ V^*$. When $\mathcal{I}_V$ has a basis
consisting of polynomials with real coefficients, $V^*=V$. Also if
$V$ is irreducible and has a smooth piece parametrized by a
neighborhood of the origin  of a complex Euclidean space through
polynomials with real coefficients, then $V^*=V$.

Next for $[\xi]\in V^*$, we define the Segre variety $Q_\xi$ of $V$
associated with $\xi$ by $Q_\xi=\{[z]\in V:
\sum_{j=1}^{n+1}z_j\xi_j=0\}$ which is a subvariety of
codimension one in $V$. Similarly, for $[z]\in V$, we define the
Segre variety $Q^*_{z}$ of $V^*$ associated with $z$ by
$Q^*_z=\{[\xi]\in V^*: \sum_{j=1}^{n+1}z_j\xi_j=0\}$. It is clear
that $[z]\in Q_\xi$ if and only if $[\xi]\in Q^*_z.$ The Segre
family of $V$  is defined to be  the projective variety ${\mathcal
M}:=\{([z],[\xi])\in V\times V^*, [z]\in Q_{\xi}\}$.

Now, we let $(M, \omega)$ be an irreducible Hermitian symmetric
space of compact type canonically  embedded in a certain minimal
projective space $\mathbb{CP}^N$, that we will describe in detail
later in this section. Then under this embedding,  its conjugate
space $M^*$ is just $M$ itself. Taking $\omega $ to be the natural
restriction of the Fubini-Study metric to $M$,  the holomorphic
isometric group of $M$ is then the restriction of a certain subgroup
of the unitary actions of the ambient space. Now, for  two points
$p_1,p_2\in M$, let $U$ be an $(N+1)\times (N+1)$ unitary matrix
such that $\sigma ([z])=[z]\cdot U$ is an isometry sending $p_1$
to $p_2$. Then $\sigma^*([\xi])=[\xi] \Ol{U}$ is an isometry of
$M^*$. By a straightforward verification, we see that $\sigma^*$
biholomorphically  sends $Q^*_{p_1}$ to $Q^*_{p_2}$. Similarly, for
any $q_1,q_2\in M^*$, $Q_{q_1}$ is unitary equivalent to $Q_{q_2}$.
In the canonical embeddings which we will describe later, the
hyperplane
 section at infinity of the manifold  is a Segre variety. Since the one at infinity is built up from Schubert
 cells and all Segre varieties are holomorphically equivalent to each other,  one deduces
 that each Segre variety of $M$ is irreducible.
 This
 fact will play a role in
 the proof of our main theorem.

\subsection{Canonical embeddings and explicit coordinate functions}
We now  describe  a special type of  canonical embedding of the Hermitian
symmetric space $M$  of compact type into $\mathbb{CP}^N$. This
embedding will play a crucial role in our computation leading to the
proof of Theorem \ref{mainthm}.  See \cite{He} for the classification of the irreducible Hermitian symmetric spaces of compact type. See also \cite{Lo1}, \cite{Lo2}  on the typical canonical embeddings of the Heritian symmetric spaces of compact type and the related theory of Hermitian positive Jordan triple system.

\medskip

${\clubsuit 1}.$ Grassmannians (spaces of type I): Write $G(p,q)$ for
the Grassmannian space consisting of  $p$ planes in
$\mathbb{C}^{p+q}. $ (Since $G(p,q)$ is biholomorphically equivalent
to $ G(q,p)$, we will assume $p\leq q$ in what follows).
There is a matrix representation of $G(p,q)$ as the equivalence
classes of $p\times(p+q)$ non-degenerate matrices under the matrix
multiplication from the left  by elements of $GL(p,\mathbb{C})$. A
Zariski open affine chart $\mathcal A$  for $G(p,q)$  is identified
with  $\mathbb{C}^{pq}$ with coordinates $Z$ for elements of the
form:
\begin{displaymath}\label{local affine}
\left(\begin{matrix}I_{p\times p}&Z
\end{matrix}\right)=
\left(\begin{matrix}1&0&0&\cdot\cdot\cdot&0&z_{11}&z_{12}&\cdot\cdot\cdot&z_{1q}\\
0&1&0&\cdot\cdot\cdot&0&z_{21}&z_{22}&\cdot\cdot\cdot&z_{2q}\\
&&&\cdot\cdot\cdot&&&\cdot\cdot\cdot&\\
0&0&0&\cdot\cdot\cdot&1&z_{p1}&z_{p2}&\cdot\cdot\cdot&z_{pq}\\
\end{matrix}\right),\ \hbox{where }\ Z \text{ is a } p\times q \text{ matrix}.
\end{displaymath}

The Pl\"ucker embedding
$G(p,q)\rightarrow\mathbb{CP}(\Lambda^p\mathbb{C}^{p+q})$ is given
by mapping the $p-$plane $\Lambda$ spanned by vectors
$v_1,...,v_p\in\mathbb{C}^{p+q}$ into the wedge product $v_1\wedge
v_2\wedge...\wedge v_p\in\wedge^{p}\mathbb{C}^{p+q}$. The action
induced by the multiplication through elements of $SU(p+q)$ from the
right induces a unitary action in the embedded ambient projective
space. In homogenous coordinates, the embedding is given by the
$p\times p$ minors of the $p\times (p+q)$ matrices (up to a sign).
   More specifically,  in the above local affine chart, we have the
   following (up to a sign in front of the components):
\begin{equation}\label{coord1}
Z\ra [1,
  Z(\begin{matrix}
 i_{1} & ... & i_{k} \\
 j_{1} & ... & j_{k}
   \end{matrix}),...]
\end{equation}
which is denoted for simplicity of notation, in what follows, by
$[1,r_{z}]=\left[1,\psi_{1},\psi_{2},...,\psi_{N} \right].$
Here and in what follows, $Z(\begin{matrix}
 i_{1} & ... & i_{k} \\
 j_{1} & ... & j_{k}
   \end{matrix}
)$ is the determinant of the submatrix of $Z$ formed by
its $i_{1}^{\text{th}},...,i_{k}^{\text{th}}$ rows and
$j_{1}^{\text{th}},...,j_{k}^{\text{th}}$ columns,
where the indices
run through  $$k=1,2,...,p,1\leq i_1<i_2<...<i_k\leq
p,1\leq j_1<j_2<...<j_k\leq q.$$ In particular when $k=1$,
$Z(\begin{matrix}
 i_{1} \\
 j_{1}
   \end{matrix})=z_{i_1j_1}.$
Notice that under such an embedding into the projective space,
$(G(p,q))^*=G(p,q).$ We thus have the same affine coordinates for
$(G(p,q))^*$:
\begin{displaymath}
\left(\begin{matrix}I_{p\times p}&\Xi
\end{matrix}\right)=
\left(\begin{matrix}1&0&0&\cdot\cdot\cdot&0&\xi_{11}&\xi_{12}&\cdot\cdot\cdot&\xi_{1q}\\
0&1&0&\cdot\cdot\cdot&0&\xi_{21}&\xi_{22}&\cdot\cdot\cdot&\xi_{2q}\\
&&&\cdot\cdot\cdot&&&\cdot\cdot\cdot&\\
0&0&0&\cdot\cdot\cdot&1&\xi_{p1}&\xi_{p2}&\cdot\cdot\cdot&\xi_{pq}\\
\end{matrix}\right),\ \ \Xi \text{ is a } p\times q \text{ matrix.}
\end{displaymath}
By the definition in $\S 2.1$, it follows that the restriction  of
the Segre family to the product of these Zariski open affine subsets
has the following canonical defining function:
\begin{equation}\label{1001} \rho(z,\xi)=1+\sum_
{\begin{subarray}{c}1\leq i_1<i_2<...<i_k\leq p,\\
 1\leq j_1<j_2<...<j_{k}\leq q\\
 k=1,...,p
 \end{subarray}}
Z(\begin{matrix}
 i_{1} & ... & i_{k} \\
 j_{1} & ... & j_{k}
   \end{matrix})\Xi
   (\begin{matrix}
  i_{1} & ... & i_{k}\\
 j_{1} & ... & j_{k}
   \end{matrix})
\end{equation}
Here $z=(z_{11},z_{12},...,z_{pq}),\
\xi=(\xi_{11},\xi_{12},...,\xi_{pq})$.
%
For simplicity of notation and terminology, we   call this
quasi-projective algebraic variety embedded in
$\mathbb{C}^{pq}\times\mathbb{C}^{pq}$, which is defined by
(\ref{1001}), the Segre family of $G(p,q)$. Our defining function $\rho(z, \xi)$ of the Segre family is closely related to the generic
norm of the corresponding Hermitian positive Jordan triple system(cf. \cite{Lo1}, \cite{Lo2}).

\medskip
$\clubsuit 2.$  Orthogonal Grassmannians (type II): Write
$G_{II}(n,n)$ for the submanifold of the Grassmannian $G(n,n)$
consisting of isotropic $n$-dimensional subspaces of
$\mathbb{C}^{2n}$.  Then $\tilde S\in G_{II}(n,n)$ if and only if
\begin{equation}\label{submanifold1}
\tilde S\left(
\begin{matrix}0&I_{n\times n}\\
I_{n\times n}&0
\end{matrix}\right)\tilde S^{T}=0.
\end{equation}
In the aforementioned  open affine piece of the Grassmannian
$G(n,n)$ with  $\tilde S=(I,S)$,
$\wt{S}\in G_{II}(n,n)$ if and only if  $S$ is an $n\times n$
antisymmetric matrix. We  identify this open affine chart ${\mathcal
A}$  of $G_{II}(n,n)$ with $\mathbb{C}^{\frac{n(n-1)}{2}}$ through
the holomorphic  coordinate map:
\begin{equation}\label{antisymmetric}
\left(
\begin{matrix}I_{n\times n}&Z
\end{matrix}
\right):= \left(
\begin{matrix}1&0&0&\cdot\cdot\cdot&0&0&z_{12}&\cdot\cdot\cdot&z_{1n}\\
0&1&0&\cdot\cdot\cdot&0&-z_{12}&0&\cdot\cdot\cdot&z_{2n}\\
&&&\cdot\cdot\cdot&&&\cdot\cdot\cdot&\\
0&0&0&\cdot\cdot\cdot&1&-z_{1n}&-z_{2n}&\cdot\cdot\cdot&0\\
\end{matrix}
\right) \ra \ (z_{12},\cdots z_{(n-1)n}).\
\end{equation}
Later in the paper we will sometimes use the notation $z_{ji}:=-z_{ij}$ if $j >i$ for this type II case. The Pl\"ucker embedding of $G(n,n)$ gives a 2-canonical embedding of
$G_{II}(n,n).$  Unfortunately this embedding is not good enough  for
our purposes later. Therefore, we will use a different  embedding in
this paper, which is  given by the spin representation of $O_{2n}$.
This embedding  is what is called a one-canonical embedding of
$G_{II}(n,n).$ We briefly describe this embedding as following. More
details can be found in [Chapter 12; PS].

Let $V$ be a real vector space of dimension $2n$ with a given inner
product, and let $\mathcal{K}(V)$ be the space consisting of all
orthogonal complex structures on V preserving this  inner product.
An element of $\mathcal{K}(V)$ is a linear  orthogonal
transformation $J:V\rightarrow V$ such that $J^2=-1.$  Any two
choices of $J$ are conjugate in the orthogonal group $O(V)=O_{2n}$,
and thus $\mathcal{K}(V)$ can be identified with the homogeneous
space $O_{2n}/U_n.$  On the other hand,  there is  a one-to-one
correspondence  assigning    the complex $J$ to a complex
$n$-dimensional isotropic subspace $W$ of
$V_{\mathbb{C}}(=V\bigotimes\mathbb{C}).$ $\mathcal{K}(V)$ has two
connected components $\mathcal{K}_{\pm}(V):$  Noticing that any
complex structure defines an orientation on $V$,  these two
components correspond to the two possible orientations on $V$. Write
one for $\mathcal{K}_{+}(V),$ which is actually our $G_{II}(n,n)$.

Now fix an isotropic $n$-dimensional subspace $W\subset V_{\mathbb
C}$ with the associated complex structure $J$  of $V_{\mathbb{C}}$
and pick a basis for V: $\{x_1,...,x_n,y_1,...,y_n\} $ with
$J(x_i)=y_i, J(y_i)=-x_i$. Then $W$ is spanned by
$\{x_i-\sqrt{-1}y_i\}_{i=1}^n. $ Define $\Ol W$ to be the space
spanned by $\{x_i+\sqrt{-1}y_i\}_{i=1}^n.$   As shown in [PS], there
is a holomorphic embedding $\mathcal{K}(V)\hookrightarrow
\mathbb{CP}(\Lambda(W))$, where $\Lambda(W)$ is the exterior algebra
of $W$.  This  embedding  is equivariant under the action of $O(V).$
Thus $\mathcal K_{+}(V)\hookrightarrow \mathbb{CP}(\Lambda(W))$ is
equivariant under $SO(V).$ Choose the open affine cell of $\mathcal
K_{+}(V)$ such that $\{Y\in \mathcal K_{+}(V)|Y\cap\Ol
W=\varnothing\}$. Then it can be identified with
$\eqref{antisymmetric}.$

We next describe the 1-canonical embedding by Pfaffians as
following:
Let $\Pi$ be the set of all partitions of $\{1,2,...,2n\}$ into
pairs  without regard to order.
An element $\alpha \in \Pi$ can be written as
$\alpha=\{ (i_{1}, j_{1}), (i_{2}, j_{2}),...,(i_{n},j_{n})\}$
with $i_{k}<j_{k}$ and $i_{1}<i_{2}<...<i_{n}.$ Let

$$\pi=\left[\begin{array}{cccccc}
             1 & 2 & 3 & 4 & ... & 2n \\
             i_{1} & j_{1} & i_{2} & j_{2} & ... & j_{n} \\
\end{array}\right]
$$
be the corresponding permutation. Given a partition $\alpha$ as
above and a $(2n)\times (2n)$ matrix $A=(a_{jk})$ , define
$$A_{\alpha}=\mathrm{sgn}(\pi)a_{i_{1}j_{1}}a_{i_{2}j_{2}}\cdots a_{i_{n}j_{n}}.$$
The Pfaffian of $A$ is then given by
$$\mathrm{pf}(A)=\sum_{\alpha \in \Pi} A_{\alpha}.$$
The Pfaffian of an $m \times m$ skew-symmetric matrix for $m$ odd is
defined to be zero.

Therefore in the coordinate system $\eqref{antisymmetric}$, the
embedding of ${\mathcal A}$ is given by
\begin{equation}\label{coord2}
[1,
,..., \mathrm{pf} (Z_{\sigma}),...].
\end{equation}
Write $S_{k}$ for the collection of all subsets of $\{1,...,n\}$
with $k$ elements.
 The  $\sigma$ in (\ref{coord2}) runs through all elements
of $S_{k}$ with  $2\le k\le n$ and $k$ even. For
$\sigma=\{i_1<\cdots<i_k\}$, $Z_\sigma$ is defined as the submatrix
$Z(\begin{matrix}
 i_{1} & ... & i_{k} \\
 i_{1} & ... & i_{k}
   \end{matrix}).$
 For instance, $ \left(\mathrm{pf} (Z_{\sigma})\right)_{\sigma \in
S_{2}}=(z_{12},...,z_{(n-1)n}).$ We also write  $\eqref{coord2}$ as
$[1,r_{z}]=\left[1,\psi_{1},\psi_{2},...,\psi_{N} \right]$ for
simplicity of notation. We choose the local coordinates for
$(G_{II}(n,n))^*$ in a similar way
\begin{equation}
\left(
\begin{matrix}I_{n\times n}&\Xi
\end{matrix}\right)=
\left(
\begin{matrix}1&0&0&\cdot\cdot\cdot&0&0&\xi_{12}&\cdot\cdot\cdot&\xi_{1n}\\
0&1&0&\cdot\cdot\cdot&0&-\xi_{12}&0&\cdot\cdot\cdot&\xi_{2n}\\
&&&\cdot\cdot\cdot&&&\cdot\cdot\cdot&\\
0&0&0&\cdot\cdot\cdot&1&-\xi_{1n}&-\xi_{2n}&\cdot\cdot\cdot&0\\
\end{matrix}
\right).
\end{equation}
The defining function for  the Segre family  (in the product of such
affine pieces) is given by
\begin{equation}\rho(z,\xi)=1+\sum_
{\begin{subarray}{c}\sigma \in S_{k},\\
2 \leq k \leq n ,2| k
 \end{subarray}}
{\rm{Pf}}(Z_{\sigma}){\rm{Pf}}
   (\Xi_{\sigma}).
\end{equation}


\medskip
$\clubsuit 3.$ Symplectic Grassmannians (type III): Write
$G_{III}(n,n)$ for the submanifold of the Grassmannian space
$G(n,n)$ defined as follows: Take the matrix representation of each
element of the Grassmannian $G(n,n)$ as an $n\times 2n$
non-degenerate matrix. Then $\wt A\in G_{III}(n,n),$ if and only if,
\begin{equation}\label{submanifold}
\wt A\left(
\begin{matrix}0&I_{n\times n}\\
-I_{n\times n}&0
\end{matrix}\right)\wt A^{T}=0.
\end{equation}
In the Zariski open affine piece  of the Grassmannian $G(n,n)$
defined before, we can take a representative matrix of the form:
$\wt A=(I, Z)$. Then we conclude that $\wt{A}\in G_{III}(n,n)$ if
and only if $Z$ is an $n\times n$ symmetric matrix. We identify this
Zariski open affine chart ${\mathcal A}$ of $G_{III}(n,n)$ with
$\mathbb{C}^{\frac{n(n+1)}{2}}$  through  the holomorphic coordinate
map:
\begin{displaymath}
\tilde A=\left(\begin{matrix}I_{n\times n}&Z
\end{matrix}\right):=
\left(\begin{matrix}1&0&0&\cdot\cdot\cdot&0&z_{11}&z_{12}&\cdot\cdot\cdot&z_{1n}\\
0&1&0&\cdot\cdot\cdot&0&z_{12}&z_{22}&\cdot\cdot\cdot&z_{2n}\\
&&&\cdot\cdot\cdot&&&\cdot\cdot\cdot&\\
0&0&0&\cdot\cdot\cdot&1&z_{1n}&z_{2n}&\cdot\cdot\cdot&z_{nn}\\
\end{matrix}\right)\ra\ \ (z_{11},\cdots, z_{nn}).
\end{displaymath}
Later in the paper we sometimes use the notation $z_{ji}:=z_{ij}$ if $j >i$ for this type III case. Through the Pl\"ucker embedding of the Grassmannian,
$G_{III}(n,n)$ is embedded into
$\mathbb{CP}(\Lambda^n\mathbb{C}^{2n})(\cong \mathbb{CP}^{N^*}).$
In the above local coordinates, we write down the embedding as (up
to a sign)
\begin{equation}
Z\ra [1,\cdots,
 Z(\begin{matrix}
 i_{1} & ... & i_{k} \\
 j_{1} & ... & j_{k}
   \end{matrix}),...]:=[1,\psi_1,\cdots,\psi_{N^*}].
\end{equation}
Choose the local affine open piece of $(G_{III}(n,n))^*$ consisting
of elements in the following form:
\begin{displaymath}
\left(\begin{matrix}I_{n\times n}&\Xi
\end{matrix}
\right)=
\left(\begin{matrix}1&0&0&\cdot\cdot\cdot&0&\xi_{11}&\xi_{12}&\cdot\cdot\cdot&\xi_{1n}\\
0&1&0&\cdot\cdot\cdot&0&\xi_{12}&\xi_{22}&\cdot\cdot\cdot&\xi_{2n}\\
&&&\cdot\cdot\cdot&&&\cdot\cdot\cdot&\\
0&0&0&\cdot\cdot\cdot&1&\xi_{1n}&\xi_{2n}&\cdot\cdot\cdot&\xi_{nn}.\\
\end{matrix}\right).
\end{displaymath}
The defining function of  Segre family in the product of such affine
open pieces is given by
\begin{equation}\label{typeIIIprerho}\rho(z,\xi)=1+\sum_
{\begin{subarray}{c}1\leq i_1<i_2<...<i_k\leq n,\\
 1\leq j_1<j_2<...<j_{k}\leq n\\
 k=1,...,n
 \end{subarray}}
Z(\begin{matrix}
 i_{1} & ... & i_{k} \\
 j_{1} & ... & j_{k}
   \end{matrix})\Xi
   (\begin{matrix}
  i_{1} & ... & i_{k}\\
 j_{1} & ... & j_{k}
   \end{matrix})
\end{equation}
However the Pl\"ucker embedding is not a useful canonical embedding
to us  for $G_{III}(n,n)$, due to the fact that $\{\psi_j\}$ is not
a linearly independent system. For instance, 
$$Z\left(
     \begin{array}{cc}
       1 & 2 \\
       3 & 4 \\
     \end{array}
   \right) + Z\left(\begin{array}{cc}
                             1 & 4 \\
                             2 & 3 \\
                           \end{array}
                         \right)  = Z\left(
                \begin{array}{cc}
                  1 & 3 \\
                  2 & 4 \\
                \end{array}
              \right).$$
This embedding can not serve our purposes here. We therefore derive
from this embedding a minimal embedding into a certain  projective
subspace
 in $\mathbb{CP}(\Lambda^n\mathbb{C}^{2n})(\cong \mathbb{CP}^{N^*})$.
   We denote this minimal projective subspace by $\mathcal H\cong\mathbb
   {CP}^N$, which is discussed in detail below.
  We notice that the embedding $G_{III}(n,n)
   \hookrightarrow\mathbb {CP}^N$ is equivariant under the transitive action of $Sp(n)$.

Following the notations we set up in  the Grassmannian case, we write
$[1,\psi_{1},\cdots\psi_{N^*}]$ for the map of the Pl\"ucker
embedding into $\mathbb{CP}^{N^*}.$ Write
$(\psi_{i_{1}},...,\psi_{i_{m_{k}}})$ for those components of degree
$k$ in $z$ among $\{ \psi_{j} \}_{j=1}^{N^*}.$ Here $1 \leq k \leq n,$ and $\{i_{1},...,i_{m_{k}}\}$
depends on $k.$ For instance, if $k=1,$ then
$$(\psi_{i_{1}},...,\psi_{i_{m_{1}}})=(z_{11},...,z_{nn}),$$
where $z_{ij}$ is repeated twice if $i \neq j.$
Let $\{\psi^{(k)}_{1},
\cdots,\psi^{(k)}_{m_{k}^*}\}$ be a maximally linearly independent subset
of $\{ \psi_{i_{1}},...,\psi_{i_{m_{k}}} \}$ over ${\mathbb R}$ (and thus
also over ${\mathbb C}$). For instance,
$$\{ \psi^{(1)}_{1},\cdots,\psi^{(1)}_{m_{1}^*} \}=\{z_{ij}\}_{i \leq j}.$$
Let $A_{k}$ be the $m_{k}^* \times m_{k}$ matrix
such that
$(\psi_{i_{1}},\cdots\psi_{i_{m_{k}}})=(\psi^{(k)}_1,\cdots\psi^{(k)}_{m_{k}^*})\cdot
A_{k}.$ Apparently $A_{k}$ has real entries and is of full rank.
Hence $A_{k}\cdot
A_{k}^t$ is positive definite. 

Then $\{\psi_1^*,\cdots,\psi^*_{N}\}: =\{ \psi^{(k)}_{1},
\cdots,\psi^{(k)}_{m_{k}^*} \}_{1 \leq k \leq n}$ forms a basis of
$\{\psi_{1},\cdots\psi_{N^*}\}$, where $N=m_{1}^* +...+ m_{n}^*$.
Moreover, if we write $A$ as the $(m_{1}^* +...+ m_{n}^*) \times
(m_{1}+...+m_{n})$ matrix:
$$A=\begin{pmatrix}
A_{1}\\
& \cdots \\
&& A_{n}
\end{pmatrix},
$$
Then $A$ has  full rank and we have a real orthogonal matrix $U$
such that
$$U=\begin{pmatrix}
U_{1}\\
& \cdots \\
&& U_{n}
\end{pmatrix}, \ \  U^t(A\cdot
A^t)U=
\begin{pmatrix}
\mu_{1}\\
& \cdots \\
&& \mu_{N}
\end{pmatrix} \ \text{ with  each}~ \mu_j>0.
$$
Here $U_{k}, 1 \leq k \leq n$, is an $m_{k}^* \times m_{k}^*$
orthogonal matrix. Now we define
\begin{displaymath}(\psi_{1}^1,...,\psi_{N_{1}}^1,\psi_{1}^2,...,\psi_{N_{2}}^2,...,\psi_{1}^{n-1},...,\psi_{N_{n-1}}^{n-1},
\psi^{n} ):=(\psi^*_1,\cdots\psi^*_{N})\cdot U\cdot\left(
\begin{matrix}{\sqrt{\mu_{1}}}&&&\\
&{\sqrt{\mu_{2}}}&&\\
&&\cdots&\\
&&&{\sqrt{\mu_{N}}}
\end{matrix}
\right).
\end{displaymath}
Here $N_{1}+...+N_{n-1}+N_{n}=N^*,$ where we set $N_{n}=1.$ We will
also sometimes write $\psi_{N_n}^n =\psi^n.$ As a direct
consequence,
\begin{equation}
\begin{aligned}(\psi_{1}^1,...,\psi_{N_{1}}^1,&\psi_{1}^2,...,\psi_{N_{2}}^2,...,\psi_{1}^{n-1},...,\psi_{N_{n-1}}^{n-1},
\psi^{n} )\cdot(\Ol{\psi_{1}^1},...,\Ol{\psi_{N_{1}}^1},\Ol{\psi_{1}^2},...,\Ol{\psi_{N_{2}}^2},...,\Ol{\psi_{1}^{n-1}},...,\Ol{\psi_{N_{n-1}}^{n-1}},
\Ol{\psi^{n}} )\\
&=(\psi_{1},\cdots,\psi_{N^*})\cdot(\Ol{\psi_{1}},\cdots,\Ol{\psi_{N^*}})={\rm det}(I+Z\bar{Z}^t)=\rho(z,\overline{z}).
\end{aligned}
\end{equation}


Moreover
$\{\psi_{1}^1,...,\psi_{N_{1}}^1,\psi_{1}^2,...,\psi_{N_{2}}^2,...,\psi_{1}^{n-1},...,\psi_{N_{n-1}}^{n-1},
\psi^{n}\}$ forms a linearly independent system; and
$\{\psi_{1}^k,...,\psi_{N_{k}}^k\}$ are polynomials in $z$ of degree
$k$ for  $k=1,...,n.$ Now our canonical embedding of the
aforementioned affine piece ${\mathcal A}$ of $G_{III}(n,n)$ is
taken as
$$z\in {\mathbb C}^{\frac{n(n+1)}{2}}
\ra [1,\psi_{1}^1,...,\psi_{N_{1}}^1,
\psi_{1}^2,...,\psi_{N_{2}}^2,...,\psi_{1}^{n-1},...,\psi_{N_{n-1}}^{n-1},
\psi^{n}].$$

For simplicity,  we will still  denote
$(\psi_{1}^1,...,\psi_{N_{1}}^1,\psi_{1}^2,...,\psi_{N_{2}}^2,...,\psi_{1}^{n-1},...,\psi_{N_{n-1}}^{n-1},
\psi^{n})$ by
\begin{equation}\label{3r}
r_{z}=\left(\psi_{1},\psi_{2},...,\psi_{N} \right)=\left(\psi_{1}^1,...,\psi_{N_{1}}^1,\psi_{1}^2,...,\psi_{N_{2}}^2,...,\psi_{1}^{n-1},...,\psi_{N_{n-1}}^{n-1},
\psi^{n}\right).
\end{equation}
Here, for instance,
$(\psi_{1},...,\psi_{\frac{n(n+1)}{2}})=(\psi_{1}^1,...,\psi_{N_{1}}^{1})=(a_{ij}z_{ij})_{1
\leq i \leq j \leq n}$, where $a_{ij}$ equals to $1$ if $i=j,$
equals to $\sqrt{2}$ if $i < j.$ Hence the defining function of the
Segre family, which is the same as $\eqref{typeIIIprerho}$, is given
by $\rho(z,\xi)=1+\sum_{i=1}^{N}\psi_i(z)\psi_i(\xi).$

\medskip

$\clubsuit 4.$  Hyperquadrics (type IV): Let $ Q^n$ be the
hypersurface in $\mathbb{CP}^{n+1}$ defined by
$$\left \{[x_{0},...,x_{n+1}] \in \mathbb{CP}^{n+1}: \sum_{i=1}^n x_{i}^2-2x_{0}x_{n+1}=0 \right\},$$
where $[x_{1},...,x_{n+2}]$ are the homogeneous coordinates for
$\mathbb{CP}^{n+1}.$
 It is invariant under the action of the  group $SO(n+2)$.  We mention that under the present embedding,
 the action is not the standard $SO(n+2)$ in $GL(n+2)$. However  it is   conjugate to the standard
 $SO(n+2)$ action
  by a certain element  $g\in U(n+2)$.
An Zariski open affine piece ${\mathcal A}\subset Q^n$
 identified with $\mathbb{C}^n$ is given by $(z_{1},...,z_{n})
\mapsto [1,\psi_1,...,\psi_{n+1}]= [1,z_{1},...,z_{n},
\frac{1}{2}\sum_{i=1}^n z_{i}^2],$ which will be denoted by
$[1,r_{z}]=\left[1,\psi_{1},\psi_{2},...,\psi_{n+1} \right]$. Choose
the same local chart for $(Q^n)^*:$ $(\xi_{1},...,\xi_{n})
\rightarrow [1,\xi_{1},...,\xi_{n}, \frac{1}{2}\sum_{i=1}^n
\xi_{i}^2]$. Then the defining function of the  Segre family
restricted to $\mathbb{C}^n\times\mathbb{C}^n\hookrightarrow
Q^n\times (Q^n)^*$ is given by

\begin{equation}
\rho(z,\xi)=1+\sum_{i=1}^n z_i\xi_i+\frac{1}{4}(\sum_{i=1}^n z_i^2)(\sum_{i=1}^n \xi_i^2)
\end{equation}

$\clubsuit 5.$ The exceptional manifold $M_{16}:= E_6\slash
SO(10)\times SO(2):$ As shown in \cite{IM1},\cite{IM2}, this exceptional
Hermitian symmetric space can be realize as the Cayley plane. Take
the exceptional $3\times3$ complex Jordan algebra
\begin{equation} \label{10000002}
\mathcal{J}_3(\mathbb{O})=\left\{\left(
\begin{matrix}c_1&x_3&\bar{x}_2\\
\bar{x}_3&c_2&x_1\\
x_2&\bar{x}_1&c_3
\end{matrix}\right):c_i\in\mathbb{C},x_i\in\mathbb{O}
\right\}\cong\mathbb{C}^{27}.
\end{equation}
Here $\mathbb{O}$ is the complexified algebra of octonions, which is
a complex vector space of dimension 8. Denote  a standard  basis of
$\mathbb{O}$ by $\{e_0,e_1,...,e_7\}$.   The multiplication rule in
terms of this basis is given in Appendix I. The conjugation operator
appeared in (\ref{10000002}) is  for octonions, which is defined as
follows: $\bar{x}=x_0e_1-x_1e_1-...-x_7e_7,$ if
$x=x_0e_0+x_1e_1+x_2e_2+...+x_7e_7,x_i\in\mathbb{C}.$
Moreover under this basis,  $\mathcal{J}_3(\mathbb{O})\cong\mathbb{C}^{27}$ is realized by identifying each matrix
\begin{equation*} 
X=\left(
\begin{matrix}\xi_1&\eta&\bar{\kappa}\\
\bar{\eta}_3&\xi_2&\tau\\
{\kappa}&\bar{\tau}&\xi_3
\end{matrix}\right)\in \mathcal{J}_3(\mathbb{O})
\end{equation*}
with the  point $(\xi_1,\xi_2,\xi_3,\eta_0,\eta_1,\dots,\eta_7,\kappa_0,\kappa_2,\dots,\kappa_7,\tau_0,\tau_1,\dots,\tau_7)\in \mathbb C^{27}$, where $\eta=\sum_{i=0}^7\eta_ie_i,\kappa=\sum_{i=0}^7\kappa_ie_i$ and $\tau=\sum_{i=0}^7\tau_ie_i.$

The Jordan multiplication is defined as $A\circ
B=\frac{1}{2}(AB+BA)$ for $A, B\in \mathcal{J}_3(\mathbb{O})$ . The
subgroup $SL$($\mathbb{O}$) of $GL$($\mathcal{J}_3$($\mathbb{O}$))
consisting of automorphisms preserving the determinant is the
adjoint group of type $E_6.$ The action of $E_6$ on the
projectivization $\mathbb{CP}\mathcal{J}_3$($\mathbb{O}$) has
exactly three orbits: the complement of the determinantal
hypersurface, the regular part of this hypersurface, and its
singular part which is the closed $E_6-$orbit.  The closed orbit is
the Cayley plane or the hermitian symmetric space of compact type
corresponding to $E_6$. It can be defined by the quadratic equation
$$X^2=\text{trace}(X)X,\quad\quad X\in\mathcal{J}_3(\mathbb{O}),$$  or as the closure of the affine
cell ${\mathcal A}$
\begin{displaymath}\mathbb{OP}_1^2=\left\{\left(
\begin{matrix}1&x&y\\
\bar{x}&x\bar{x}&y\bar{x}\\
\bar{y}&x\bar{y}&y\bar{y}
\end{matrix}\right):x,y\in\mathbb{O}
\right\}\cong\mathbb{C}^{16}
\end{displaymath}
in the local coordinates $(x_0,x_1,...,x_7,y_0,...,y_7)$.  The
precise formula for the canonical embedding map is given in Appendix
II.  We denote this embedding by $[1,r_{z}]=\left
[1,\psi_{1},\psi_{2},...,\psi_{N} \right ].$

To find the defining function for its Segre family over the product
of such  standard affine sets, we  choose local coordinates for the
conjugate Cayley plane to be
$(\kappa_0,\kappa_1,...,\kappa_7,\eta_0,\eta_1,...,\eta_7).$
 Then
\begin{equation}\label{e16rho}\rho(z,\xi)=1+\sum_{i=0}^7x_i\kappa_i+\sum_{i=0}^7y_i\eta_i+\sum_{i=0}^7A_i(x,y){A}_i(\kappa,\eta)
+B_0(x,y){B}_0(\kappa,\eta)+B_1(x,y) B_1(\kappa,\eta),
\end{equation}
where  $A_j, B_j$ are defined as in Appendix I,
$z=(x_0,...,x_7,y_0,...,y_7)$ and $
\xi=(\kappa_0,...,\kappa_7,\eta_0,...,\eta_7).$

\medskip

$\clubsuit 6.$ The other exceptional manifold $M_{27}=E_7\slash
E_6\times SO(2):$ As shown in \cite{CMP}, it can be realized as the
Freudenthal variety. Consider the Zorn algebra
$$\mathcal{Z}_2(\mathbb{O})=\mathbb{C}\bigoplus\mathcal{J}_3(\mathbb{O})\bigoplus\mathcal{J}_3(\mathbb{O})\bigoplus\mathbb{C}$$
One can prove that there exists an action of $E_7$ on that
$56-$dimensional vector space (see \cite{Fr}). The closed
$E_7-$orbit inside $\mathbb{CP}\mathcal{Z}_2(\mathbb{O})$ is the
Freudenthal variety
 $E_7\slash E_6\times
SO(2).$  An affine cell ${\mathcal A}$ of Freudenthal variety is
$[1,X,{\rm{Com}}(X),\det(X)]\in\mathbb{CP}\mathcal{Z}_2
 (\mathbb{O}).$ Here $X$ belongs to $\mathcal{J}_3({\mathbb O})$; ${\rm{Com}}(X)$ is the comatrix of $X$ such that $X{\rm{Com}}(X)=\det(X)I$ under the usual matrix multiplication rule. Notice that ${\rm{Com}}(X)=X\times X$, where $X\times X$ is the Freudenthal multiplication defined as follows (see \cite{O}):
 $$X\times X:=X^2-{\rm tr}(X)X+\frac{1}{2}({\rm tr}(X)^2-{\rm tr} (X^2)I.$$
 For  explicit expressions for $X\times X$ and $\det(X)$ in terms of the entries of $X$, see \cite{O} or Appendix I in this paper.

The embedding of  $E_7\slash E_6\times
SO(2)\hookrightarrow\mathbb{CP}^N$ in local coordinates $z$ is given
in Appendix I.  Choose  the local affine open piece for $(E_7\slash
E_6\times SO(2))^*$ with coordinates $$\xi=
(\xi_1,\xi_2,\xi_3,\eta_0,...,\eta_7,\kappa_0,...,\kappa_7,\tau_0,...,\tau_7).$$
We denote this embedding  by
$[1,r_{z}]=\left[1,\psi_{1},\psi_{2},...,\psi_{N} \right].$ The
defining function for the Segre family is then
$\rho(z,\xi)=1+r_{z}\cdot{r}_{\xi}$, where
\begin{equation}\label{e27rho}
\begin{aligned}
&r_{z}=(x_{1},x_{2},x_{3},y_{0},...,y_{7},t_{0},...,t_{7},w_{0},...,w_{7},A(z),B(z),C(z),D_{0}(z),...D_{7}(z),\\
&    \ \ \ \ \ E_{0}(z),...,E_{7}(z),F_{0}(z),...,F_{7}(z),G(z))\\
 &r_{\xi}=(\psi_{1}(\xi),\psi_{2}(\xi),...,\psi_{N}(\xi))=(\xi_1,\xi_2,\xi_3,\eta_0,...,\eta_7,\kappa_0,...,\kappa_7,\tau_0,...,\tau_7,\\
 &\ \ \ \ \ \ \ {A}(\xi),{B}(\xi),{C}(\xi),{D}_{0}(\xi),...,{D}_{7}(\xi),
{E}_{0}(\xi),...,{E}_{7}(\xi),{F}_{0}(\xi),...,{F}_{7}(\xi),{G}(\xi))
\end{aligned}
\end{equation}
Here see Appendix I for the definition of the functions appeared in
the formula.

\bigskip
Summarizing the above,  for each irreducible Hermitian symmetric
space of compact type $M$ of dimension $n$, we now have described  a
canonical embedding from $M$ into a projective space ${\mathbb
P}^N$, which restricted to a certain Zariski open affine piece
${\mathcal A}$ holomorphically equivalent to  ${\mathbb C}^n$ takes
the form: $z(\in {\mathbb C}^n)\mapsto [1,\kappa_1
z_1,\cdots,\kappa_i z_i,\cdots,\kappa_n z_n, O(z^2)]$. Here
$\kappa_i=1$ for all $i$ except in the case of type III where
$\kappa_i$ can be $1$ or $\sqrt 2$. This is the embedding we will
use in later discussions. Notice in our embedding, the conjugate
space $M^*$ is the same as $M$. For simplicity of notation, we will
also write ${\mathcal M}$ for the restriction of the Segre family of
$M$ restricted to ${\mathcal A}\times {\mathcal A}^*= {\mathbb
C}^n\times {\mathbb C}^n$. From this embedding and
the invariant property of Segre varieties, we immediately conclude
the following:

\begin{lemma}\label{unique}Assume $A$ and $B$ are two distinct points of
$ M$.  Then their associated Segre varieties are different, namely,
$Q_A\not =Q_B$.
\end{lemma}
{\it Proof of Lemma \ref{unique}}: Since the holomorphic isometric
group acts transitively on $M$, we can assume
$A=(0,0,...,0)\in\mathbb C^n\cong\mathcal A\subset M.$ Therefore
$Q_A$ is the  hyperplane section of $M\hookrightarrow\mathbb P^N$ at
infinity, namely, $Q_A=M\backslash \mathcal A.$ Now if $B\in\mathcal
A$, because $B\not=(0,0,...,0)$, there are non-trivial linear terms
in the defining function of $Q_B$.  This leads to the fact that the
defining function of $Q_B$ has to be a non-constant polynomial in
$\mathbb C[\xi_1,...,\xi_n]$. Therefore $Q_B\cap {\mathbb C^n}\not
=\emptyset$ and thus does not coincide with $Q_A.$ If $B\in\mathcal
M\backslash\mathcal A,$ by the  symmetric property of Segre
varieties, we have  $(0,...,0)\in Q_B$. Therefore $Q_B\not =Q_A.$ We
then arrive at the conclusion.\ \ \ $\endpf$

\medskip Finally, since in our setting,  $M^*=M$ and the Segre family
on $M$ and $M^*$ are the same. For simplicity of notation, we do not
distinguish, in what follows,  $Q^*$ and ${\mathcal M}^*$ from $Q$
and ${\mathcal M}$, respectively.

\subsection{Explicit expression of the volume forms}
From now on, we assume that $M$ is an irreducible Hermitian
symmetric space of compact type and we choose the canonical
embedding $M\hookrightarrow\mathbb{CP}^N$ as described in $\S 2.2$
according to its type. We denote the metric on $M$ induced from
Fubini-Study of $\mathbb{CP}^N$ by $\omega$, and the volume form by
$d\mu=\omega^n$ (up to a positive constant). Notice that the metric
we obtained is always invariant under the action of a certain
transitive subgroup $G\subset Aut(M)$ (which comes from the
restriction of a subgroup of the unitary group of the ambient
projective space). Hence by a theorem of Wolf \cite{W}, $\omega$ is
the unique $G$ invariant metric on $M$ up to a scale. We claim
$\omega$ must be K\"ahler-Einstein. Indeed, since the Ricci form
Ric($\omega$) of $\omega$ is  invariant under $G$,  for a small
$\epsilon,\ \omega+\epsilon\rm{Ric}(\omega)$ is thus also a $G$
invariant metric on $M$. By \cite{W}, it is a multiple of $\omega$,
and thus $\rm{Ric}(\omega)=\lambda\omega.$
 Write $d\mu$ as the product of  $V$ and the
standard Euclidean volume form over the affine subspace ${\mathcal
A}$, where $V$ is a positive function in $z.$ Since
${\rm{Ric}}(\omega)=- i \partial\bar\partial\log V$,
$-i\partial\bar\partial\log V=\lambda\omega$.
 Notice that $\lambda>0$.
In the local affine open piece ${\mathcal A}$ defined before,
$\omega=i\partial\bar\partial\log\rho(z,\bar z)$, where
$\rho(z,\xi)$ is the defining function for the associated   Segre
family. As we will see later ($\S 7$), $\rho(z,\xi)$ is an
irreducible polynomial in $(z,\xi)$.
Then we have
$$\partial\bar\partial\log (V \rho(z,\bar z)^{\lambda})=0.$$
Hence, $\log (V \rho(z,\bar z)^{\lambda})=\phi(z)+\Ol{\psi(z)},$
where both $\phi$ and $\psi$ are holomorphic functions. Therefore $V=\frac{e^{\phi(z)+\Ol{\psi(
z)}}}{\rho(z,\bar z)^{\lambda}}.$ Because $\rho(z,\xi)$ is an
irreducible polynomial, from the way $V$ is defined, $V$ must be a
rational function of the form $\frac{p(z,\Ol{z})}{\rho(z,\bar z)^m}$
with $p,\rho$ relatively prime to each other. Since $\phi,\psi$ are globally
defined, by a monodromy argument, it
is clear that $\lambda$ has to be an integer. Also both $e^{\phi(z)}$
and $e^{\Ol{\psi(\xi)}}$ must be rational functions. Again, since
$\phi, \psi$ are also globally defined, this forces $\phi,\psi$ to
be constant functions.
Therefore, we conclude that $V=c{\rho(z,\bar z)^{-\lambda}}$.
Here $\lambda$ is a certain positive integer and $c$ is a positive constant. Next by a
well-known result (see \cite{BaMa}), two K\"ahler-Einstein metrics
of $M$ are different by an automorphism of $M$ (up to a positive
scalar multiple). Therefore, to prove Theorem \ref{mainthm}, we can
assume, without loss of generality, that the K\"ahler-Einstein
metric in Theorem \ref{mainthm} is the metric obtained by
restricting the Fubini-Study metric to $M$ through the embedding
described  in this section.

\section{A basic property  for partially degenerate holomorphic maps}
In this section, we introduce a notion of degeneracy for holomorphic
maps and derive an important consequence, which  will be
fundamentally applied in the proof of our main theorem.

Let $\psi(z):=(\psi_{1}(z),...,\psi_{N}(z))$ be a vector-valued holomorphic function from a
 neighborhood $U$ of $0$ in $\mathbb{C}^m, m \geq 2$, into $\mathbb{C}^N, N > m,$ with $\psi(0)=0.$ Here we write
  $z=(z_{1},...,z_{m})$ for the coordinates of $\mathbb{C}^m.$ In the following, we will write
$\widetilde{z}=(z_{1},...,z_{m-1}),$ i.e., the vector  $z$ with  the
last component $z_{m}$ being dropped out.  Write
$\frac{\partial^{|\alpha|}}{\partial \widetilde{z}^{\alpha}}=
\frac{\partial^{|\alpha|}}{\partial z_{1}^{\alpha_{1}}...\partial
z_{m-1}^{\alpha_{m-1}}}$ for an $(m-1)-$multiindex $\alpha,$ where
$\alpha=(\alpha_{1},...,\alpha_{m-1}).$ Write
$$\frac{\partial^{|\alpha|}}{\partial \widetilde{z}^{\alpha}} \psi(z)=\left( \frac{\partial^{|\alpha|}}{\partial \widetilde{z}^{\alpha}} \psi_{1}(z),...,\frac{\partial^{|\alpha|}}{\partial \widetilde{z}^{\alpha}} \psi_{N}(z)\right).$$

We introduce the following definition.

\begin{definition}\label{df21}
Let $k \geq 0.$ For a point $p \in U,$ write $E_{k}(p)=\mathrm{Span}_{\mathbb{C}}\{
\frac{\partial^{|\alpha|}}{\partial \widetilde{z}^{\alpha}} \psi(z)|_{z=p}:
0 \leq |\alpha| \leq k \}.$
We write $r$ for  the greatest number such that for any neighborhood
$O$ of $0,$ there exists $p \in O$ with $\mathrm{dim}_{\mathbb{C}}
 E_{k}(p)=r.$\ \ \ $r$ is called the $k-$th  $\widetilde{z}-$rank of $\psi$ at $0,$
 which is
written   as $\mathrm{rank}_{k}(\psi, \widetilde{z}).\ ~F$ is called
$\widetilde{z}-$nondegenerate if $\mathrm{rank}_{k_{0}}(\psi,
\widetilde{z})=N$ for some $k_{0} \geq 1.$
\end{definition}

\begin{remark}
It is easy to see that $\mathrm{rank}_{k}(\psi,\wt{z})=r$ if and
only if the following matrix
$$\left(\begin{matrix}
 \frac{\partial^{|\alpha^0|}}{\partial \widetilde{z}^{\alpha^0}} \psi(z) \\
  ...\\
  ...\\
 \frac{\partial^{|\alpha^s|}}{\partial \widetilde{z}^{\alpha^s}} \psi(z)
  \end{matrix}\right)
$$
has an $r \times r$ submatrix with determinant not identically zero
for $z \in U$ for some multiindices $\{\alpha^{0},...,\alpha^{s}\}$
with all $0 \leq |\alpha^{j}| \leq k.$ Moreover, any $l \times l\
(l>r)$ submatrix of the matrix has identically zero determinant for
any choice of $\{\alpha^{0},...,\alpha^{s}\}$ with  $0 \leq
|\alpha^{j}| \leq k.$

In particular, $\psi$ is $\widetilde{z}-$nondegenerate if and only if there exist multiindices $\beta^{1},...,\beta^{N}$ such that
$$\left| \begin{matrix}
  \frac{\partial^{|\beta^{1}|}}{\partial \widetilde{z}^{\beta^{1}}} \psi_{1}(z)& ...  &  \frac{\partial^{|\beta^{1}|}}{\partial \widetilde{z}^{\beta^{1}}} \psi_{N}(z) \\
 ... & ... & ... \\
 \frac{\partial^{|\beta^{N}|}}{\partial \widetilde{z}^{\beta^{N}}}\psi_{1}(z) & ... & \frac{\partial^{|\beta^{N}|}}{\partial \widetilde{z}^{\beta^{N}}}\psi_{N}(z)
\end{matrix}
 \right|$$
is not identically zero. Moreover,  $\mathrm{rank}_{i+1}(\psi, \widetilde{z}) \geq \mathrm{rank}_{i}(\psi,\widetilde{z})$ for any $i \geq 0.$
\end{remark}

For the rest of this section, we further assume that the first $m$
components of $\psi$, i.e.,
$(\psi_{1},...,\psi_{m}): \mathbb{C}^m \rightarrow \mathbb{C}^m$
is a biholomorphic map in a neighborhood of $0 \in \mathbb{C}^m.$
Then we have,

\begin{lemma}\label{22-1} It holds that
$\mathrm{rank}_{0}(\psi, \widetilde{z})=1, \mathrm{rank}_{1}(\psi,
\widetilde{z})=m,$  and for $k \geq 1, \mathrm{rank}_{k}(\psi,
\widetilde{z}) \geq m.$
\end{lemma}

\medskip
{\it Proof of Lemma \ref{22-1}}: We first notice that it holds
trivially that $\mathrm{rank}_{0}(\psi, \widetilde{z})=1$, for $F$
is not identically zero. We now
 prove $\mathrm{rank}_{1}(\psi, \widetilde{z})=m.$ First notice that $\mathrm{rank}_{1}(\psi, \widetilde{z}) \leq m$
 as there are only $m$ distinct multiindices $\beta$ such that $|\beta| \leq 1.$ On the other hand,
  since $\psi$ has full rank at $0,$ we have,
$$
\left| \begin{matrix}
 \frac{\partial \psi_{1}}{\partial z_{1}} & ... &  \frac{\partial \psi_{m}}{\partial z_{1}} \\
 ... & ...  & ...  \\
 \frac{\partial \psi_{1}}{\partial z_{m}} & ... &  \frac{\partial \psi_{m}}{\partial z_{m}}
         \end{matrix}
 \right|(0) \neq 0.$$
 This together with the fact $\psi(0)=0$ implies that the $z_{m}$ derivative of
\begin{equation}\label{eq33}
\left| \begin{matrix}
 \psi_{1} &...  & \psi_{m} \\
 \frac{\partial \psi_{1}}{\partial z_{1}} & ... & \frac{\partial \psi_{m}}{\partial z_{1}} \\
 ... & ... & ... \\
 \frac{\partial \psi_{1}}{\partial z_{m-1}} & ... & \frac{\partial \psi_{m}}{\partial z_{m-1}}
\end{matrix}
\right|
\end{equation}
is nonzero at $p=0.$ Consequently, the quantity in (\ref{eq33}) is
not identically zero in $U.$ By the definition of the $\wt{z}$-rank,
we then arrive at the conclusion.\ \ \ $\endpf$

\medskip
We now prove the following degeneracy theorem in terms of its
$\wt{z}$-rank, which will be  used to derive Theorem \ref{thm211}.
\begin{theorem}\label{thm24}
Let $\psi=(\psi_{1},...,\psi_{m},\psi_{m+1},...,\psi_{N})$ be a
holomorphic map from a neighborhood of $0 \in \mathbb{C}^m$ into
$\mathbb{C}^N$ with $\psi(0)=0$. Recall that $\widetilde{z}=(z_{1},...,z_{m-1}),$ i.e., the vector  $z$ with  the
last component $z_{m}$ being dropped out.  Assume  that
$(\psi_{1},...,\psi_{m})$ is a biholomorphic map from a neighborhood
of $0 \in \mathbb{C}^m$ into a neighborhood of $0 \in \mathbb{C}^m.$
Suppose
\begin{equation}\label{eq34}
{\rm rank}_{N-m+1}(\psi, \widetilde{z})< N.
\end{equation}
Then there exist $N$ holomorphic functions
$g_{1}(z_{m}),...,g_{N}(z_{m})$ near $0$ in the $z_{m}-$Gauss plane
with $\{ g_{1}(0),...,g_{N}(0)\}$ not all zero such that the
following holds for any $(z_{1},...,z_{m})$ near $0.$
\begin{equation}\label{eq35}
\sum_{i=1}^N g_{i}(z_{m})\psi_{i}(z_{1},...,z_{m}) \equiv 0.
\end{equation}
In particular, one can make one of the $\{g_{i}\}_{i=1}^N $ to be identically
one.
\end{theorem}

The geometric intuition for the theorem is as follows: The space of
1-jets has dimension $m$ by Lemma \ref{22-1}. We expect that at
least one more dimension is increased when we go from the space of
$k$-jets to the space of $(k+1)-$ jets until we reach the maximum
possible value $N$. The theorem  says that if   this process fails,
namely, the assumption in \eqref{eq34} holds,   we then  end up with
a function relationship as in  \eqref{eq35}.

\medskip
 {\it Proof of Theorem \ref{thm24}}: We consider the following
set,
$$\mathcal{S}=\{l\geq 1: \mathrm{rank}_{l}(\psi, \widetilde{z}) \leq l+m-2\}.$$
Note that $1 \notin \mathcal{S}$, for $\mathrm{rank}_{1}(F)=m.$ We
claim that $\mathcal{S}$ is not empty. Indeed, we have  $1+ N-m \in
\mathcal{S}$ by (\ref{eq34}). Now write $t'$ for  the minimum number
in $\mathcal{S}.$ Then $2 \leq t' \leq 1+N-m.$ Moreover, by the
choice of $t',$
\begin{equation}
\mathrm{rank}_{t'}(\psi, \widetilde{z}) \leq t'+m-2,\
\mathrm{rank}_{t'-1}(\psi, \widetilde{z})\ \geq t'+ m -2.
\end{equation}
This yields that
\begin{equation}
\mathrm{rank}_{t'}(\psi, \widetilde{z})=\mathrm{rank}_{t'-1}(\psi, \widetilde{z})=t'+m-2.
\end{equation}
We  write  $t:=t'-1, \   n: =t'+m -2 .$ Here we note $t \geq 1, m \leq n \leq N-1.$ Then there exist
multiindices $\{\gamma^{1},...,\gamma^{n}\}$ with each $|\gamma^{i}|
\leq t$ and $j_{1},...,j_{n}$ such that
\begin{equation}\label{eq38}
\Delta(\gamma^{1},...,\gamma^{n}|j_{1},...,j_{n}):=\left| \begin{matrix}
 \frac{\partial^{|\gamma^{1}|} \psi_{j_{1}}}{\partial \tilde{z}^{\gamma^{1}}} & ... & \frac{\partial^{|\gamma^{1}|} \psi_{j_{n}}}{\partial \tilde{z}^{\gamma^{1}}} \\
 ... & ...  & ... \\
\frac{\partial^{|\gamma^{n}|} \psi_{j_{1}}}{\partial \tilde{z}^{\gamma^{n}}}  & ... &  \frac{\partial^{|\gamma^{n}|} \psi_{j_{n}}}{\partial \tilde{z}^{\gamma^{n}}}
\end{matrix}\right|~\text{is not identically zero in}~U.
\end{equation}
Since $\mathrm{rank}_{1}(\psi, \widetilde{z})=m,$  we can choose
$(\gamma^{1},...,\gamma^{n}|j_{1},...,j_{n})$ such that
$$\gamma^{1}=(0,..,0),
\gamma^{2}=(1,0,...,0),...,\gamma^{m}=(0,...,0,1).$$ For any
$\alpha^{1},...,\alpha^{n+1}$ with $|\alpha^{i}| \leq t+1,$ and
$l_{1},...,l_{n+1}$, we have
\begin{equation}\label{eq39}
\Delta(\alpha^{1},...,\alpha^{n+1}|l_{1},...,l_{n+1})= \left|\begin{matrix}
 \frac{\partial^{|\alpha^{1}|} \psi_{l_{1}}}{\partial \tilde{z}^{\alpha^{1}}} & ...  & \frac{\partial^{|\alpha^{1}|} \psi_{l_{n}}}{\partial \tilde{z}^{\alpha^{1}}} & \frac{\partial^{|\alpha^{1}|} \psi_{l_{n+1}}}{\partial \tilde{z}^{\alpha^{1}}} \\
 ... & ... & ... & ... \\
 ... & ... & ...  & ... \\
 \frac{\partial^{|\alpha^{n+1}|} \psi_{l_{1}}}{\partial \tilde{z}^{\alpha^{n+1}}} & ... &  \frac{\partial^{|\alpha^{n+1}|} \psi_{l_{n}}}{\partial \tilde{z}^{\alpha^{n+1}}} & \frac{\partial^{|\alpha^{n+1}|} \psi_{l_{n+1}}}{\partial \tilde{z}^{\alpha^{n+1}}}
\end{matrix}\right| \equiv 0~\text{in}~U.
\end{equation}

We write $\varGamma$ for the collection of
$(\gamma^{1},...,\gamma^{n}|j_{1},...,j_{n}), j_{1}<...<j_{n},$ with
$\gamma^{1}=(0,..,0)$ and with $\eqref{eq38}$ being held. We
associate each $(\gamma^{1},...,\gamma^{n}|j_{1},...,j_{n})$ with an
integer $s(\gamma^{1},...,\gamma^{n}|j_{1},...,j_{n}):=s_{0}$ where
$s_{0}$ is the least number $s \geq 0$ such that
$$\frac{\partial^{s_{1}+...+s_{m-1}+s} \Delta(\gamma^{1},...,\gamma^{n}|j_{1},...,j_{n})}{\partial z^{s_{1}}_{1} \partial z^{s_{2}}_{2}...\partial z^{s_{m-1}}_{m-1}\partial z^s_{m}}(0) \neq 0.$$
for some integers $s_{1},...,s_{m-1}.$ Then $s(\gamma^{1},...,\gamma^{n}|j_{1},...,j_{n}) \geq 0$ for any $(\gamma^{1},...,\gamma^{n}|j_{1},...,j_{n}) \in \varGamma.$

Let $(\beta^{1},...,\beta^{n}|i_{1},...,i_{n}) \in \varGamma, i_{1}
< ...< i_{n} $ be  indices with the least
$s(\gamma^{1},...,\gamma^{n}|j_{1},...,j_{n})$ among all
$(\gamma^{1},...,\gamma^{n}|j_{1},...,j_{n}) \in \varGamma.$

\bigskip

\bigskip

We write $\{i_{n+1},...,i_{N}\}=\{1,...,N\} \setminus \{i_{1},..,i_{n}\},$
where $i_{n+1} <...< i_{N}.$
Write $\tilde{U}=\{z \in U: \Delta(\beta^{1},...,\beta^{n}|i_{1},...,i_{n}) \neq 0\}.$
We then have the following:

\bigskip

\begin{lemma}\label{lemma27}
Fix $j \in \{i_{n+1},...,i_{N}\}.$ Let $i \in \{i_{1},..,i_{n}\}.$
Write $\{i'_{1},...,i'_{n-1}\}= \{i_{1},...,i_{n}\} \setminus
\{i\}.$ There exists a holomorphic function $g^{j}_{i}(z_{m})$ in
$\tilde{U}$ which only depends on $z_{m}$ such that the following
holds for $z \in \tilde{U}:$
\begin{equation}\label{eq40}
\left|\begin{matrix}
 \frac{\partial^{|\beta^{1}|} \psi_{i'_{1}}}{\partial \tilde{z}^{\beta^{1}}} & ... & \frac{\partial^{|\beta^{1}|} \psi_{i'_{n-1}}}{\partial \tilde{z}^{\beta^{1}}} & \frac{\partial^{|\beta^{1}|} \psi_{j}}{\partial \tilde{z}^{\beta^{1}}}  \\
  ... & ... & ... &  ...  \\
 ... & ... & ... & ...  \\
 \frac{\partial^{|\beta^{n}|} \psi_{i'_{1}}}{\partial \tilde{z}^{\beta^{n}}} & ...  & \frac{\partial^{|\beta^{n}|} \psi_{i'_{n-1}}}{\partial \tilde{z}^{\beta^{n}}}  & \frac{\partial^{|\beta^{n}|} \psi_{j}}{\partial \tilde{z}^{\beta^{n}}}
\end{matrix}\right|(z)
=g^{j}_{i}(z_{m})\left|\begin{matrix}
 \frac{\partial^{|\beta^{1}|} \psi_{i'_{1}}}{\partial \tilde{z}^{\beta^{1}}} & ... & \frac{\partial^{|\beta^{1}|} \psi_{i'_{n-1}}}{\partial \tilde{z}^{\beta^{1}}} &
 \frac{\partial^{|\beta^{1}|} \psi_{i}}{\partial \tilde{z}^{\beta^{1}}}  \\
 ... & ... & ... & ...  \\
 ... & ...  & ...  & ...  \\
 \frac{\partial^{|\beta^{n}|} \psi_{i'_{1}}}{\partial \tilde{z}^{\beta^{n}}} & ... & \frac{\partial^{|\beta^{n}|} \psi_{i'_{n-1}}}{\partial \tilde{z}^{\beta^{n}}} &
\frac{\partial^{|\beta^{n}|} \psi_{i}}{\partial \tilde{z}^{\beta^{n}}}
\end{matrix}
\right|(z),
\end{equation}
or equivalently,
\begin{equation}\label{eqnpsi27}
\left|\begin{matrix}
 \frac{\partial^{|\beta^{1}|} \psi_{i'_{1}}}{\partial \tilde{z}^{\beta^{1}}} & ... & \frac{\partial^{|\beta^{1}|} \psi_{i'_{n-1}}}{\partial \tilde{z}^{\beta^{1}}}
 & \frac{\partial^{|\beta^{1}|} (\psi_{j}-g^{j}_{i}(z_{m})\psi_{i})}{\partial \tilde{z}^{\beta^{1}}}  \\
  ... & ... & ... &  ...  \\
 ... & ... & ... & ...  \\
 \frac{\partial^{|\beta^{n}|} \psi_{i'_{1}}}{\partial \tilde{z}^{\beta^{n}}} & ...  &
\frac{\partial^{|\beta^{n}|} \psi_{i'_{n-1}}}{\partial
\tilde{z}^{\beta^{n}}}  & \frac{\partial^{|\beta^{n}|}
(\psi_{j}-g^j_{i}(z_{m})\psi_{i})}{\partial \tilde{z}^{\beta^{n}}}
\end{matrix}\right| \equiv 0.
\end{equation}
\end{lemma}

\bigskip

{\it Proof of Lemma \ref{lemma27}}:
For simplicity of notation, we  write $ \frac{\partial }{\partial
\tilde{z}^{\beta^{i}}}$ for $ \frac{\partial^{|\beta^{i}|}
}{\partial \tilde{z}^{\beta^{i}}}$, and for $\mu= i$ or $j,$ write the matrix

$$V_{\mu}:=\left(\begin{matrix}
 \frac{\partial \psi_{i'_{1}}}{\partial \tilde{z}^{\beta^{1}}} & ... & \frac{\partial \psi_{i'_{n-1}}}{\partial \tilde{z}^{\beta^{1}}} & \frac{\partial \psi_{\mu}}{\partial \tilde{z}^{\beta^{1}}}  \\
  ... & ... & ... &  ...  \\
 ... & ... & ... & ...  \\
 \frac{\partial \psi_{i'_{1}}}{\partial \tilde{z}^{\beta^{n}}} & ...  & \frac{\partial \psi_{i'_{n-1}}}{\partial \tilde{z}^{\beta^{n}}}  & \frac{\partial \psi_{\mu}}{\partial \tilde{z}^{\beta^{n}}}
\end{matrix}\right)=\left[\begin{matrix} {\bf v}_{\mu}^1 \\ \vdots \\ {\bf v}_{\mu}^n    \end{matrix} \right],$$
where ${\bf v}_{\mu}^1, \cdots, {\bf v}_{\mu}^n$ are the row vectors of $V_{\mu}.$ To prove (\ref{eq40}), one just
needs to show that, for each $1 \leq \nu \leq m-1,$
\begin{equation}\label{eq42}
\frac{\partial}{\partial z_{\nu}} \frac{\mathrm{det}(V_j)}{\mathrm{det}(V_i)}\equiv 0~\text{in}~\tilde{U}.
\end{equation}
Indeed, by the quotient rule, the numerator of the left-hand side of (\ref{eq42}) equals to
\begin{equation}\nonumber
\begin{split}
\mathrm{det} &\left( \begin{matrix} \mathrm{det}(V_i) & \mathrm{det}(V_j) \\ \frac{\partial}{\partial z_{\nu}} \mathrm{det}(V_i) &  \frac{\partial}{ \partial z_{\nu}} \mathrm{det}(V_j) \end{matrix} \right) \\
= \mathrm{det} & \left( \begin{matrix} \mathrm{det}(V_i) & \mathrm{det}(V_j) \\ \mathrm{det}\left[\begin{matrix}  \frac{\partial}{\partial z_{\nu}}{\bf v}_{i}^1 \\ {\bf v}_{i}^2 \\ \vdots \\ {\bf v}_{i}^n    \end{matrix} \right]  & \mathrm{det}\left[\begin{matrix}  \frac{\partial}{\partial z_{\nu}}{\bf v}_{j}^1 \\ {\bf v}_{j}^2 \\ \vdots \\ {\bf v}_{j}^n    \end{matrix} \right]  \end{matrix} \right) + \cdots + 
\mathrm{det} \left( \begin{matrix} \mathrm{det}(V_i) & \mathrm{det}(V_j) \\ \mathrm{det}\left[\begin{matrix}  {\bf v}_{i}^1 \\ \vdots \\ {\bf v}_{i}^{n-1} \\ \frac{\partial}{\partial z_{\nu}}{\bf v}_{i}^n    \end{matrix} \right]  & \mathrm{det}\left[\begin{matrix}  {\bf v}_{j}^1 \\ \vdots \\ {\bf v}_{j}^{n-1} \\ \frac{\partial}{\partial z_{\nu}} {\bf v}_{j}^n    \end{matrix} \right]  \end{matrix} \right).
\end{split}
\end{equation}
By (\ref{eq39}) and Lemma 4.4 in [BX], each term on the
right-hand side of the equation above equals $0.$  For instance, the
last term above equals to
\begin{equation}
\left|
    \begin{array}{cc}
      \left|\begin{matrix}
 \frac{\partial \psi_{i'_{1}}}{\partial \tilde{z}^{\beta^{1}}} & ... & \frac{\partial \psi_{i'_{n-1}}}{\partial \tilde{z}^{\beta^{1}}} &
 \frac{\partial \psi_{i}}{\partial \tilde{z}^{\beta^{1}}}  \\
 ... & ... & ... & ...  \\
 ... & ...  & ...  & ...  \\
 \frac{\partial \psi_{i'_{1}}}{\partial \tilde{z}^{\beta^{n}}} & ... & \frac{\partial \psi_{i'_{n-1}}}{\partial \tilde{z}^{\beta^{n}}} &
\frac{\partial \psi_{i}}{\partial \tilde{z}^{\beta^{n}}}
\end{matrix}
\right| & \left|\begin{matrix}
 \frac{\partial \psi_{i'_{1}}}{\partial \tilde{z}^{\beta^{1}}} & ... & \frac{\partial \psi_{i'_{n-1}}}{\partial \tilde{z}^{\beta^{1}}} & \frac{\partial \psi_{j}}{\partial \tilde{z}^{\beta^{1}}}  \\
  ... & ... & ... &  ...  \\
 ... & ... & ... & ...  \\
 \frac{\partial \psi_{i'_{1}}}{\partial \tilde{z}^{\beta^{n}}} & ...  & \frac{\partial \psi_{i'_{n-1}}}{\partial \tilde{z}^{\beta^{n}}}  & \frac{\partial \psi_{j}}{\partial \tilde{z}^{\beta^{n}}}
\end{matrix}\right| \\
    \left| \begin{matrix}
  \frac{\partial \psi_{i'_{1}}}{\partial \tilde{z}^{\beta^{1}}} & ... & \frac{\partial \psi_{i'_{n-1}}}{\partial \tilde{z}^{\beta^{1}}} &
 \frac{\partial \psi_{i}}{\partial \tilde{z}^{\beta^{1}}}  \\
 ... & ... & ... & ...  \\
 \frac{\partial \psi_{i'_{1}}}{\partial \tilde{z}^{\beta^{n-1}}} & ...  & \frac{\partial \psi_{i'_{n-1}}}{\partial \tilde{z}^{\beta^{n-1}}} & \frac{\partial \psi_{i}}{\partial \tilde{z}^{\beta^{n-1}}}  \\
 \frac{\partial}{\partial z_{\nu}}(\frac{\partial \psi_{i'_{1}}}{\partial \tilde{z}^{\beta^{n}}}) & ... & \frac{\partial}{\partial z_{\nu}}(\frac{\partial \psi_{i'_{n-1}}}{\partial \tilde{z}^{\beta^{n}}}) &
\frac{\partial}{\partial z_{\nu}}(\frac{\partial \psi_{i}}{\partial \tilde{z}^{\beta^{n}}})
\end{matrix}
\right| & \left|\begin{matrix}
 \frac{\partial \psi_{i'_{1}}}{\partial \tilde{z}^{\beta^{1}}} & ... & \frac{\partial \psi_{i'_{n-1}}}{\partial \tilde{z}^{\beta^{1}}} & \frac{\partial \psi_{j}}{\partial \tilde{z}^{\beta^{1}}}  \\
  ... & ... & ... &  ...  \\
 \frac{\partial \psi_{i'_{1}}}{\partial \tilde{z}^{\beta^{n-1}}} & ... & \frac{\partial \psi_{i'_{n-1}}}{\partial \tilde{z}^{\beta^{n-1}}} & \frac{\partial \psi_{j}}{\partial \tilde{z}^{\beta^{n-1}}}  \\
 \frac{\partial}{\partial z_{\nu}}(\frac{\partial \psi_{i'_{1}}}{\partial \tilde{z}^{\beta^{n}}}) & ...  & \frac{\partial}{\partial z_{\nu}} (\frac{\partial \psi_{i'_{n-1}}}{\partial \tilde{z}^{\beta^{n}}})  & \frac{\partial}{\partial z_{\nu}}(\frac{\partial \psi_{j}}{\partial \tilde{z}^{\beta^{n}}})
\end{matrix}\right| \\
    \end{array}
  \right|.
\end{equation}
It is a multiple of the following determinant  (by Lemme 4.4 in [BX]):
\begin{equation}
\left|
  \begin{array}{ccccc}
    \frac{\partial \psi_{i'_{1}}}{\partial \tilde{z}^{\beta^{1}}} & ... & \frac{\partial \psi_{i'_{n-1}}}{\partial \tilde{z}^{\beta^{1}}} & \frac{\partial \psi_{i}}{\partial \tilde{z}^{\beta^{1}}} & \frac{\partial \psi_{j}}{\partial \tilde{z}^{\beta^{1}}} \\
    ... & ... & ... & ... & ... \\
    \frac{\partial \psi_{i'_{1}}}{\partial \tilde{z}^{\beta^{n}}} & ... &  \frac{\partial \psi_{i'_{n-1}}}{\partial \tilde{z}^{\beta^{n}}}  & \frac{\partial \psi_{i}}{\partial \tilde{Z}^{\beta^{n}}} & \frac{\partial \psi_{j}}{\partial \tilde{z}^{\beta^{n}}} \\
    \frac{\partial \psi_{i'_{1}}}{\partial \tilde{z}^{\beta^{n+1}}} & ... & \frac{\partial \psi_{i'_{n-1}}}{\partial \tilde{z}^{\beta^{n+1}}}  & \frac{\partial \psi_{i}}{\partial \tilde{z}^{\beta^{n+1}}} & \frac{\partial \psi_{j}}{\partial \tilde{z}^{\beta^{n+1}}} \\
  \end{array}
\right|,
\end{equation}
where $\frac{\partial}{\partial
\tilde{z}^{\beta^{n+1}}}=\frac{\partial}{\partial
z_{\nu}}(\frac{\partial} {\partial \tilde{z}^{\beta^{n}}}),$ which
is identically zero by (\ref{eq39}). This establishes Lemma
\ref{lemma27}. $\endpf$
\medskip

The extendability of $g_{i}^j(z_{m})$ will be needed for our later
argument:

\begin{lemma}\label{lemmaext}
For any $i,j$ as above, the holomorphic function $g^j_{i}(z_{m})$
can be extended holomorphically to a neighborhood of $0$ in  the
$z_{m}-$plane.
\end{lemma}

{\it Proof of Lemma \ref{lemmaext}}: First, $g^j_{i}$ is defined on
the projection $\pi_{m}(\tilde{U})$ of $\tilde{U},$ where $\pi_{m}$
is the natural projection of $(z_{1},...,z_{m})$ to its last
component $z_{m}.$ If $0 \in \pi_{m}(\tilde{U}),$  the claim follows
trivially. Now assume that $0 \notin \pi_{m}(\tilde{U}).$ If we
write $s=s(\beta_{1},...,\beta_{n}|i_{1},...,i_{n}),$ by its
definition, then there exists $(a_{1},...,a_{m-1}) \in
\mathbb{C}^{m-1}$ close to $0,$ such that
\begin{equation}\label{eq43}
\left|\begin{matrix}
 \frac{\partial^{|\beta^{1}|} \psi_{i'_{1}}}{\partial \tilde{z}^{\beta^{1}}} & ... & \frac{\partial^{|\beta^1|} \psi_{i'_{n-1}}}{\partial \tilde{z}^{\beta^{1}}} &
 \frac{\partial^{|\beta^1|} \psi_{i}}{\partial \tilde{z}^{\beta^{1}}}  \\
 ... & ... & ... & ...  \\
 ... & ...  & ...  & ...  \\
 \frac{\partial^{|\beta^n|} \psi_{i'_{1}}}{\partial \tilde{z}^{\beta^{n}}} & ... & \frac{\partial^{|\beta^n|} \psi_{i'_{n-1}}}{\partial \tilde{z}^{\beta^{n}}} &
\frac{\partial^{|\beta^n|} \psi_{i}}{\partial \tilde{z}^{\beta^{n}}}
\end{matrix}
\right|(a_{1},...,a_{m-1},z_{m})=c z^{s}_{m} + o(|z_{m}|^{s}),~ c \neq 0.
\end{equation}

Then there exists $r > 0$ small enough such that for any $0 < |z_{m}|< r, (a_{1},...,a_{m-1},z_{m}) \in \tilde{U}.$ That is, at any of such points,
equation (\ref{eq43}) is not zero.

We now substitute  $(a_{1},...,a_{m-1},z_{m}), 0 < |z_{m}| <r,$ into
the equation (\ref{eq40}), and compare the vanishing order as
$z_{m}\ra 0$:
\begin{equation}
c_{1}z^{s'}_{m}+ o(|z_{m}|^{s'})=
g^j_{i}(z_{m})(cz^s_{m}+o(|z_{m}|^s)),~c \neq 0.
\end{equation}
for some $s' \geq 0.$  Note that $0 \leq s \leq s'$  by the definition of $s$ and the choice of $(\beta_{1},...,\beta_{n}|i_{1},...,i_{n}).$
The holomorphic  extendability across $0$ of $g^j_{i}(z_{m})$ then follows easily.
$\endpf$
\medskip

We next make the following observation:
\medskip
\begin{claim}\label{claim1}
 For each fixed  $j \in \{i_{n+1},...,i_{N}\}$ and any $i'_{1} <...<i'_{n-1}$ with $\{i'_{1},...,i'_{n-1}\} \subset \{i_{1},...,i_{n}\},$ we have:
\begin{equation}\label{eqnpsi1}
\left| \begin{matrix}
  \frac{\partial^{|\beta^1|} \psi_{i'_{1}}}{\partial \tilde{z}^{\beta^{1}}} & ... &  \frac{\partial^{|\beta^1|} \psi_{i'_{n-1}}}{\partial \tilde{z}^{\beta^{1}}} &
 \frac{\partial^{|\beta^1|} (\psi_{j}-\sum_{k=1}^{n} g^j_{i_{k}}\psi_{i_{k}}) }{\partial \tilde{z}^{\beta^{1}}}  \\
 ... & ... & ... & ...  \\
 ... & ... & ...  & ...  \\
 \frac{\partial^{|\beta^n|} \psi_{i'_{1}}}{\partial \tilde{z}^{\beta^{n}}} & ... & \frac{\partial^{|\beta^n|} \psi_{i'_{n-1}}}{\partial \tilde{z}^{\beta^{n}}} &
\frac{\partial^{|\beta^n|} (\psi_{j}-\sum_{k=1}^{n} g^j_{i_{k}}\psi_{i_{k}}) }{\partial \tilde{z}^{\beta^{n}}}
\end{matrix}\right|(z) \equiv 0, \forall z \in \tilde{U}.
\end{equation}
\end{claim}
{\it Proof of Claim \ref{claim1}}:
 Note that for each $i'_{l}, 1 \leq l \leq n-1,$ the
following trivially holds:
\begin{equation}\label{eqnpsi2}
\left| \begin{matrix}
  \frac{\partial^{|\beta^1|} \psi_{i'_{1}}}{\partial \tilde{z}^{\beta^{1}}} & ... &  \frac{\partial^{|\beta^1|} \psi_{i'_{n-1}}}{\partial \tilde{z}^{\beta^{1}}} &
 \frac{\partial^{|\beta^1|} (g^j_{i'_{l}}\psi_{i'_{l}}) }{\partial \tilde{z}^{\beta^{1}}}  \\
 ... & ... & ... & ...  \\
 ... & ... & ...  & ...  \\
 \frac{\partial^{|\beta^n|} \psi_{i'_{1}}}{\partial \tilde{z}^{\beta^{n}}} & ... & \frac{\partial^{|\beta^n|} \psi_{i'_{n-1}}}{\partial \tilde{z}^{\beta^{n}}} &
\frac{\partial^{|\beta^n|} (g^j_{i'_{l}}\psi_{i'_{l}}) }{\partial \tilde{z}^{\beta^{n}}}
\end{matrix}\right|(z) \equiv 0,
\end{equation}
for the last column in the matrix is a multiple of  one of the first
$(n-1)$ columns. Then (\ref{eqnpsi1}) is an immediate consequence of (\ref{eqnpsi27}) and
(\ref{eqnpsi2}).  $\endpf$
\medskip

\begin{lemma}\label{lemma29}
For each fixed $j \in \{i_{n+1},...,i_{N}\},$
we have $\psi_{j}(z)- \sum_{k=1}^n g^j_{i_{k}}(z_{m}) \psi_{i_{k}}(z)\equiv 0$
 for any $z \in \tilde{U},$ and thus it
holds also for all $z \in U.$
\end{lemma}
{\it Proof of Lemma \ref{lemma29}}:
This can be concluded easily from the following Lemma \ref{BX47} and Claim \ref{claim1}. Here one needs to use the fact that $\beta^1=(0,...,0).$
$\endpf$

\begin{lemma}\label{BX47}
(\cite{BX}, Lemma 4.7) Let ${\bf{b}}_{1},\cdots,{\bf{b}}_{n}$ and
${\bf{a}}$ be  $n$-dimensional column vectors
 with elements in $\mathbb{C}$, and let $B=({\bf{b}}_{1},\cdots,{\bf{b}}_{n})$ denote the $n\times n$ matrix.
   Assume that $\mathrm{det}B \neq 0$ and
 $\mathrm{det}({\bf{b}}_{i_{1}},{\bf{b}}_{i_{2}},\cdots,{\bf{b}}_{i_{n-1}},{\bf{a}})=0$ for
 any $1 \leq i_{1} < i_{2} < \cdots < i_{n-1} \leq n.$ Then $\bf{a}=0.$
\end{lemma}
Theorem \ref{thm24} now follows easily from Lemma \ref{lemma29}.\
$\endpf$
\medskip

If we further assume that $\psi_{i}(z), m+1 \leq i \leq N,$ vanishes
at least to the second order, then we have the following, which
plays a crucial role in our proof of Theorem \ref{mainthm}.
\begin{theorem}\label{thm211}
Let $\psi=(\psi_{1},...,\psi_{m},\psi_{m+1},...,\psi_{N})$ be a
holomorphic map from a neighborhood of $0 \in \mathbb{C}^m$ into
$\mathbb{C}^N$ with $\psi(0)=0$. Assume that
$(\psi_{1},...,\psi_{m})$ is a biholomorphic map from a neighborhood
of $0 \in \mathbb{C}^m$ into a neighborhood of $0 \in \mathbb{C}^N.$
Assume that $\psi_{j}(z)=O(|z|^2)$ for $m+1 \leq j \leq N.$ Suppose
that ${\rm rank}_{N-m+1}(\psi)< N.$ Then there exist
$a_{m+1},...,a_{N} \in \mathbb{C}$ that are not all zero such that
\begin{equation}\label{eq47}
\sum_{i=m+1}^N a_{j}\psi_{j}(z_{1},...,z_{m-1},0) \equiv 0,
\end{equation}
for all $(z_{1},...,z_{m-1})$ near $0.$
\end{theorem}
{\it Proof of Theorem \ref{thm211}}:
We first have the following:

\medskip
\begin{claim}\label{22-2} For each $1 \leq i \leq m,$  $ g_{i}(0)=0.$
\end{claim}

\medskip

{\it Proof of Claim \ref{22-2}}: Suppose not. Write ${\bf
c}:=(g_{1}(0),...,g_{m}(0)) \neq 0.$ Then
$(g_{1}(z_m),...,g_{m}(z_m))$
$={\bf c}+ O(|z_m|). $ The fact that
$\psi_{i}(z)=O(|z|^2), i \geq m+1,$ implies
\begin{equation}\label{eq48}
\sum_{i=1}^{m}g_{i}(z_m)\psi_{i}(z)=O(|z|^2).
\end{equation}
Notice that (the Jacobian of) $(\psi_{1},...,\psi_{m})$ is of full rank at $0.$ 
Hence
\begin{equation}
\left( \begin{matrix}
 \frac{\partial \psi_{1}}{\partial z_{1}}(0) & ... &  \frac{\partial \psi_{m}}{\partial z_{1}}(0) \\
 ... & ...  & ...  \\
 \frac{\partial \psi_{1}}{\partial z_{m}}(0) & ... &  \frac{\partial \psi_{m}}{\partial z_{m}}(0)
         \end{matrix}
 \right) {\bf c}^t \neq 0.
\end{equation}
This is a contradiction to (\ref{eq48}).
$\endpf$

\medskip
Finally, letting $z_{m}=0$ in equation (\ref{eq35}), we  obtain
(\ref{eq47}). By claim \ref{22-2}, $ (g_{m+1}(0),...,g_{N}(0))\not =0.$
This establishes Theorem \ref{thm211}.\ \ $\endpf$

\section{Proof of the main theorem assuming three extra propositions}

 In this section,  we give a proof of our main theorem under several
extra assumptions (i.e., Propositions (I)-(III)), which will be verified one by one  in the later
sections.

 Let $M\subset\mathbb{CP}^N$ be an irreducible   Hermitian symmetric space  of compact
 type, which has been canonically (and isometrically) embedded in the complex projective space through the way
  described in $\S 2$. In this section, we write $n$ as the complex dimension of $M.$ We also
 have
 on $M$ an affine  open piece ${\mathcal A}$ that is   biholomorphically equivalent to the
 complex Euclidean space of the same dimension,  such that $M\setminus{\mathcal A}$ is a codimension one complex subvariety
 of $M$.
 We identify  the coordinates of ${\mathcal A}$ by the parametrization map with
  $z=(z_1,...,z_n)\in\mathbb{C}^n$ through what is
 described in  $\S 2$, which we wrote  as $[1,\psi_1,...,\psi_N]$,
where $\psi_1,...,\psi_N$  are polynomial maps in $(z_1,...,z_n)$
with $\psi_j=\kappa_jz_j$, where $\kappa_j=1$ or $\sqrt 2$, for $j=1,\cdots, n$.
 We also write  $\Ol{F}(\xi)$ for $\Ol{F(\Ol{\xi}})$ for $\xi=(\xi_1,...,\xi_n)\in {\mathbb C}^n.$
 We still use  $\rho(z,\xi)$ for the
defining function of the Segre family of $M$ restricted to
${\mathcal A}\times {\mathcal A}^*$, which will be  canonically
identified with $\mathbb{C}^n\times\mathbb{C}^n$. Since the
coefficients of $\psi_1,...,\psi_N$ are all real, $\Ol{\psi}=\psi$
and ${\mathcal A}^*={\mathcal A}$. Hence, we have
\begin{equation}\label{unirho}
  \rho(z,\xi)=1+\sum_{i=1}^N\psi_i(z)\psi_i(\xi).
\end{equation}
Recall the standard metric $\omega$ of $M$ on ${\mathcal A}$  is
given by
\begin{equation}\label{unimetric}
  \omega=i \partial \overline{\partial}
\mathrm{log }(\rho(z,\bar{z})).
\end{equation}
 The volume form $d\mu=c_n\omega^n$  associated to  $\omega$,
by $\S 2$ ,  is now  given in ${\mathcal A}$ by the multiplication
of  $V$ with the standard Euclidean volume form, where
\begin{equation}\label{univolume}
 V =\frac{c}{(\rho(z,\bar{z}))^\lambda}
\end{equation}
with  $c > 0$ and  ${\lambda}$  a certain positive integer depending
on $M$. For instance, ${\lambda}=p+q$ when $X=G(p,q)$ \cite{G}. Here
$c_n$ is a certain positive constant depending only on $n$.





\begin{theorem}\label{main}
Let ${\mathcal A}\subset M $ be as above equipped  with the standard
metric $\omega.$ Let $F_j,j=1,...,m,$ be a holomorphic map from $U
\subset{\mathcal A}$ into $M$,
where $U$ is a connected open neighborhood of $\mathcal A$. Assume
that $F_j^*(d\mu)\not \equiv 0$ for each $j$  and assume that
\begin{equation}\label{volume preserving equation}
  d\mu= \sum_{j=1}^{m}\lambda_jF_j^{*}(d\mu),
\end{equation}
for certain positive constants  $\lambda_j >0$ with $j=1,\cdots, m$.
Then for any $j\in\{1,2,...,m\}$, $F_{j}$ extends to  a holomorphic
isometry of $(M,\omega)$.
\end{theorem}

For convenience of our discussions, we first fix some notations: In
what follows, we identify $\mathcal A$ with ${\mathbb C}^n$ having
$z=(z_1,\cdots,z_n)$ as its coordinates.
On $U\subset{\mathcal A}\subset M$ and after shrinking $U$ if
needed, we write the holomorphic map $F_j$, for $j=1,...,m$, from
$U\ra {\mathcal A}={\mathbb C}^n$, as follows:
\begin{equation}\label{rational functions}
 F_j=(F_{j,1},F_{j,2},...,F_{j,n}),\ j=1,...,m.
\end{equation}
Still write the holomorphic  embedding  from ${\cal A}$ into
$\mathbb{CP}^N$ as $[1,\psi_1,\cdots,\psi_N]$. We define
$\mathcal{F}_j(z)=({\mathcal F}_{j,1},...,{\mathcal
F}_{j,N})=(\psi_1(F_j),\psi_2(F_j),...,\psi_N(F_j))$ for $
j=1,...,m.$ Finally, all Segre varieties and Segre families are
restricted to ${\mathcal A}={\mathbb C}^n$.


The main purpose of this section is to  give a proof of  Theorem
$\ref{main}$, assuming the following three propositions hold. These propositions
will be separately established in terms of the type of $M$ in $\S 5,\S
6$ and $\S 7$. This then completes the proof of our main theorem.

\medskip

\textbf{Proposition} (I): Write
$\mathcal{L}_{i}=\frac{\partial}{\partial z_{i}}-
\frac{\frac{\partial \rho}{\partial z_{i}}(z,\xi)}{\frac{\partial
\rho}{\partial z_n}(z,\xi)}\frac{\partial}{\partial z_n}, 1 \leq i
\leq n-1,$ which are holomorphic vector fields (whenever defined)
tangent to the Segre family ${\mathcal M}$ of
$M\hookrightarrow\mathbb {CP}^N$ restricted to ${\mathcal
A}\times{\mathcal A}^*=\mathbb C^n\times\mathbb C^n$ defined by
$\rho(z, \xi)=0.$ Under the notations we set up above,
for any local biholomorphic map $F=(f_1,\cdots,f_n): U\ra
{\mathbb C}^n$ with $F(0)=0$,
there are $
z^0\in U,\xi^0\in Q_{z^0},\beta^1,...,\beta^N$, such that
\begin{equation}
\frac{\partial \rho}{\partial z_n}(z^0,\xi^0)\not =0,\ \ \
\Lambda(\beta^1,...,\beta^N)(z^0,\xi^0):=\left|\begin{matrix}
 \mathcal{L}^{\beta^1}\mathcal F_{1} & ...  & \mathcal{L}^{\beta^1}\mathcal F_{N} \\
... & ...  & ...  \\
\mathcal{L}^{\beta^N}\mathcal F_{1}  & ...  &
\mathcal{L}^{\beta^N}\mathcal F_{N}
\end{matrix}\right|(z^ 0,\xi^0)\not =0.
\end{equation}
Here $\beta^l=(k^l_1,...,k^l_{n-1}),$ $k^l_1,...,k^l_{n-1}$ are
non-negative integers, for $l=1,2,...,N;$ $\beta^1=(0,0,...,0);$
$\mathcal{L}^{\beta^l}=\mathcal{L}_{1}^{k^l_1}\mathcal{L}_2^{k^l_2}\mathcal{L}_3^{k^l_3}...\mathcal{L}_{n-1}^{k^l_{n-1}};$
$\mathcal{F}(z)=({\mathcal F}_{1},...,{\mathcal
F}_{N})=(\psi_1(F),\psi_2(F),...,\psi_N(F))$.
 Moreover, $s_l:=\sum_{i=1}^{n-1}k^l_i \ (l=1,...,N)$ is a
non-negative integer bounded from above by a universal constant
depending only on $(M,\omega)$.  Also, in what follows, when we like
to emphasize the dependence of $\Lambda(\beta^1,...,\beta^N)$ on
$F$, we also write it as $\Lambda_{F}(\beta^1,...,\beta^N)$.

\medskip

\textbf{Proposition} (II):  Suppose that $\xi^0\in\mathbb{C}^n$ with
$\xi^0\neq(0,0,...,0)$. Then for a generic smooth point $z^0$ on the
Segre variety $Q_{\xi^0}$ and a small neighborhood $U\owns z^0$,
there is a $z^1 \in U\cap Q_{\xi^0}$ such that  $Q_{z^0}$ and $
Q_{z^1}$ both are smooth at $\xi^0$ and intersect transversally at
$\xi^0$, too. Moreover,
we can find a biholomorphic parametrization near $\xi^0$:
$(\xi_1,\xi_2,...,\xi_n)=\mathcal{G}(\tilde\xi_1,\tilde\xi_2,...,\tilde\xi_n)$
with $(\tilde\xi_1,\tilde\xi_2,...,\tilde\xi_n)\in U_1\times
U_2\times...\times U_n\subset\mathbb C^n$, where  $U_1$ and $U_2$
are  small neighborhoods of $1\in {\mathbb C}$, and $U_j$ for $j\ge
3$ are small neighborhoods of $0\in {\mathbb C}$
such that (i). ${\mathcal G}(1,1,0,\cdots,0)=\xi_0$, (ii).
$\mathcal{G}(\{\tilde\xi_1=1\}\times U_2\times...\times U_n)\subset
Q_{z^0},\mathcal{G}(U_1\times\{\tilde\xi_2=1\}\times
U_3\times...\times U_n)\subset Q_{z^1},$ and (iii).
$\mathcal{G}(\{\tilde\xi_1=t\}\times U_2\times...\times U_n)\
\hbox{or }\ \mathcal{G}(U_1\times\{\tilde\xi_2=s\}\times
U_3\times...\times U_n),s\in U_1,t\in U_2$  is an  open piece  of a
certain  Segre variety for each fixed $t$ and $s$. Moreover
$\mathcal{G}$ consists of algebraic functions with total degree
bounded by a constant depending only  on the manifold $M$.

\medskip

\textbf{Proposition}\ (III):  For any $\xi\neq0(z\not = 0,\
\rm{respectively}) \in {\mathbb C}^n$, $\rho(z,\xi)$ is an
irreducible polynomial in $z$ (and in $\xi$, respectively). (In
particular, $Q^*_\xi$ and $Q_{z}$ are irreducible.) Moreover, if $U$
is a connected open set in $\mathbb{C}^n$, then the Segre family
$\mathcal{M}$ restricted to $U\times\mathbb{C}^n$ is an irreducible
complex subvariety and thus its regular points form a connected
complex submanifold.
In particular, $\mathcal{M}$  is an irreducible complex subvariety
of $\mathbb{C}^n\times\mathbb{C}^n$.




\medskip
The rest of this section is splitted into several subsections. In the
first subsection, we discuss a partial algebraicity for a certain
component  $F_{j_0}$  in Theorem $\ref{main}$. In $\S 4.2$, we show
$F_{j_0}$ is algebraic. In $\S 4.3$, we further prove the
rationality of $F_{j_0}$. $\S 4.4 $ is devoted to proving that
$F_{j_0}$ extends to a birational map from $M$ to itself
and extends to a holomorphic isometry, which can be used,  through
an induction argument, to prove Theorem $\ref{main}$ assuming Propositions (I)-(III).

\subsection{An algebraicity lemma}




We use the notations we have set up so far. We now proceed to the
proof Theorem $\ref{main}$ under the hypothesis that Propositions  (I)-(III) hold.

\medskip
Denote by $J_f(z)$ the determinant of the complex Jacobian matrix of a holomorphic map $f: B \rightarrow \mathbb{C}^n,$ where $B \subset \mathbb{C}^n$ is an open subset and $z=(z_1, \cdots, z_n) \in B.$ For any holomorphic map $g(\xi)$ from an open subset of $\mathbb{C}^n$ to $\mathbb{C}^m,$ where $\xi \in \mathbb{C}^n,$ we define $\overline{g}(\xi):=\overline{g(\overline{\xi})}.$

Now from $\eqref{unirho}\eqref{unimetric}\eqref{univolume}\eqref{volume preserving
 equation}$,
we obtain
\begin{equation}\label{orirational}
\sum_{j=1}^{m}\lambda_j\frac{|J_{F_j}(z)|^2}{(1+\sum_{i=1}^N\psi_i(F_{j}(z))\psi_i(\overline F_{j}
(\bar z)))^{\lambda}}
= \frac{1}{(1+\sum_{i=1}^N\psi_i(z)\psi_i(\bar z))^{\lambda}},\ \
z=(z_1,...,z_n)\in U.
\end{equation}
Recall that  $F_j=(F_{j,1},F_{j,2},...,F_{j,n}),j=1,...,n.$
Complexifying $\eqref{orirational}$, we have
\begin{equation}\label{unicomplexify}
\sum_{j=1}^{m}\lambda_j\frac{J_{F_j}(z)\overline{J_{F_j}}(\xi)}{(1+\sum_{i=1}^N\psi_i(F_{j}(z))\psi_i(\overline F_{j}(\xi)))^{\lambda}}=
\frac{1}{(1+\sum_{i=1}^N\psi_i(z)\psi_i(\xi))^{\lambda}},\ \
(z,\xi)\in U\times {\rm conj}(U).
\end{equation}
Here ${\rm conj}(U)=:\{z:\ov{z}\in U\}$. Using  the
transitive action of the holomorphic isometric group of $(M,\omega)$
on $M$, we assume that $0\in U$, $F_j(0)=0\in {\mathcal A}$ and
$J_{F_j}(0)\not = 0$ for each $j$. Also, letting $U=B_{r}(0)$ for a
sufficiently small $r>0,$ we have ${\rm conj}(U)=U$. Hence, we will
assume that $\eqref{unicomplexify}$ holds for $(z, \xi) \in U \times
U$.

 We will need the following algebraicity lemma.
\begin{lemma}\label{HYlemma}
Let $F_j's$ be as in Theorem $\ref{main}$. Then there exist  Nash
algebraic maps
$$\widehat{F}_{1}(z,X_{1},...,X_{m}),...,\widehat{F}_{m}(z,X_{1},...,X_{m})$$
holomorphic in $(z,X_{1},...,X_{m})$ near
$(0,\overline{J_{F_{1}}}(0),...,\overline{J_{{F}_{m}}}(0))\in\mathbb
C^n\times\mathbb C^m$ such that
\begin{equation}
\overline{F}_{j}(z)=\widehat{F}_{j}(z, \overline{J_{F_{1}}}(z),...,\overline{J_{F_{m}}}(z)),j=1,...,m
\end{equation}
for $z=(z_1,...,z_n)$ near $0.$
\end{lemma}

{\it Proof of Lemma \ref{HYlemma}}:
Recall that $\psi_i=\kappa_iz_i$, where $\kappa_i=1$ or $\sqrt 2$, for $i=1,\cdots,n$ and $\psi_i=O(|z|^2)$ is
a  polynomial of $z$ for each $n+1\le i\le N$.
We obtain  from (\ref{unicomplexify}) the following:
\begin{align*}
\sum_{j=1}^{m}\lambda_j\big(J_{F_{j}}(z)\overline{J_{F_{j}}}(\xi)-\lambda(\sum_{i=1}^n
(J_{F_{j}}(z) \kappa_iF_{j,i}(z))(\overline{J_{F_{j}}}(\xi)\kappa_i
\overline{F}_{j,i}(\xi)))+ P_{j}(z,\overline{F_{j}}(\xi),
\overline{J_{F_{j}}}(\xi))\big)
\end{align*}
\begin{align}\label{simeq5}
=\frac{1}{(1+\sum_{i=1}^N\psi_i(z)\psi_i(\xi))^{\lambda}}.
\end{align}
Here each
$P_{j}(z,\overline{F_{j}}(\xi), \overline{J_{F_{j}}}(\xi))$
is a rational function in $z,\overline{F_{j}}(\xi)$ and $ \overline{J_{F_{j}}}(\xi)$.

We now set $X_{j}=J_{F_{j}}, 1 \leq j \leq m.$ Set $Y_{j}, 1 \leq j
\leq m,$ to be  the vectors:
$$Y_{j}=(Y_{j1},...,Y_{jn}):=(\kappa_1J_{F_{j}}F_{j,1},...,\kappa_nJ_{F_{j}}F_{j,n}).$$ Then equation (\ref{simeq5}) can be rewritten as
\begin{equation}\label{simeq6}
\sum_{j=1}^m\lambda_j\left( X_{j}(z)\overline{X}_{j}(\xi)-{\lambda}
Y_{j}(z)\cdot \overline{Y}_{j}(\xi) + Q_{j}(z,
\overline{X}_{j}(\xi),
\overline{Y}_{j}(\xi))\right)=\frac{1}{(1+\sum_{i=1}^N\psi_i(z)\psi_i(\xi))^{\lambda}}
\end{equation}
over $U \times U.$ Here each $Q_{j}$ with  $1 \leq j \leq m$ is
rational in $\overline X_{j}, \overline Y_{j}.$ Moreover, each
$Q_{j}, 1 \leq j \leq m$, has no terms of the form $\overline
X_{j}^k\overline Y_{js}^l$ with $l \leq 1$ for any $s \geq 1$ in its
Taylor expansion at $(\overline{X_j}(0),\overline{Y_j}(0))$.

We write $D^{\alpha}=\frac{\partial^{|\alpha|}}{\partial
z_{1}^{\alpha_{1}}...\partial z_{n}^{\alpha_{n}}}$ for an
$n-$multiindex $\alpha=(\alpha_{1},...,\alpha_{n}).$  Taking
differentiation in (\ref{simeq6}), we obtain, for each multiindex
$\alpha,$ the following:
\begin{align*}
&\sum_{j=1}^m\left((D^{\alpha}X_{j}(z))\overline{X}_{j}(\xi)-{\lambda}
(D^{\alpha}Y_{j}(z))\cdot \overline{Y}_{j}(\xi) + D^{\alpha}Q_{j}(z,
\overline{X}_{j}(\xi), \overline{Y}_{j}(\xi)) \right)\\
&\,\,\,\,\,\,\,\,\,\,\,\,\,\,\,\,\,\,\,\,\,\,\,\,\,\,\,\,\,\,\,\,\,\,\,\,\,\,
\,\,\,\,\,\,\,\,\,\,\,\,\,\,\,\,\,\,\,\,\,\,\,\,\,\,\,\,\,\,\,\,\,\,\,\,\,\,\,
=D^{\alpha}\big(\frac{1}{(1+\sum_{i=1}^N\psi_i(z)\psi_i(\xi))^{\lambda}}\big).
\end{align*}
Again each $D^{\alpha}Q_{j}, 1 \leq j \leq m,$ is rational in
$(\overline{X}_{j},\overline{Y}_{j})$ and has no terms of the form
$\overline X_{j}^k\overline Y_{js}^l$ with $l \leq 1$ and $s \geq 1$
in its Taylor expansion at $(\overline{X_j}(0),\overline{Y_j}(0))$.
Applying a similar argument  as in [Proposition 3.1, \cite{HY1}], we
can algebraically solve for $\ov{F_j}$ to complete the proof of the
lemma. \  $\endpf$
\medskip

Let $\mathcal{R}$ be the field of rational functions in $z=(z_1,...,z_n).$ Consider the field extension

$$\mathcal{E}=\mathcal{R}(\overline{J_{F_{1}}}(z),...,\overline{J_{F_{m}}}(z)).$$
Let $K$ be the transcendental degree of the field extension
$\mathcal{E}/\mathcal{R}.$ If $K=0,$ then each of
$\{\overline{J_{F_{1}}},...,\overline{J_{F_{m}}}\}$ is Nash
algebraic. As a consequence of  Lemma \ref{HYlemma}, each $F_{j}$
with $1 \leq j \leq m$ is Nash algebraic. Otherwise, by re-ordering
the indices if necessary, we let
$\mathcal{G}=\{\overline{J_{F_{1}}},...,\overline{J_{F_{K}}}\}$ be
the maximal algebraic independent subset of
$\{\overline{J_{F_{1}}},...,\overline{J_{F_{m}}}\}.$ It follows that
the transcendental degree of $\mathcal{E}/\mathcal{R}(\mathcal{G})$
is zero.  For any $l >K,$ there exists a minimal polynomial
$P_{l}(z,X_{1},...,X_{K},X)$ such that
$P_{l}(z,\overline{J_{F_{1}}}(z),...,\overline{J_{F_{K}}}(z),\overline{J_{F_{l}}}(z))
\equiv 0.$ Moreover,

$$\frac{\partial P_{l}(z,X_{1},...,X_{K},X)}{\partial X}(z, \overline{J_{F_{1}}}(z),...,\overline{J_{F_{K}}}(z),\overline{J_{F_{l}}}(z))\not \equiv 0$$
in a small neighborhood $V$ of $0,$ for otherwise, $P_{l}$ cannot be
a minimal polynomial of $\overline{J_{F_{l}}}(z).$ Now the union of
the vanishing set  of the partial derivative with respect to  $X$ in
the above equation  for each $l$  forms a proper local complex
analytic variety near $0.$ Applying the algebraic version of the
implicit function theorem, there exists a small connected open
subset $U_{0} \subset U,$ with $0 \in \overline{U}_{0}$ and a
holomorphic algebraic function $\widehat{h}_{l}, l > K,$ in a
certain neighborhood $\widehat{U}_{0}$ of $\{(z,
\overline{J_{F_{1}}}(z),...,\overline{J_{F_{K}}}(z)): z \in U_{0}
\}$ in $\mathbb{C}^n \times \mathbb{C}^K,$ such that
$$\overline{J_{F_{l}}}(z)=\widehat{h}_{l}(z,\overline{J_{F_{1}}}(z),...,\overline{J_{F_{K}}}(z)),$$
for any $z \in U_{0}$. (We can assume here $U_{0}$ is the projection
of $\widehat{U}_{0}$). Substitute this into
$$\widehat{F}_{i}(z, \overline{J_{F_{1}}}(z),...,\overline{J_{F_{m}}}(z)),$$
and still denote it, for simplicity
of notation, by $\widehat{F_{j}}(z, \overline{J_{F_{1}}}(z),...,\overline{J_{F_{K}}}(z))$ with
$$\widehat{F_{j}}(z, \overline{J_{F_{1}}}(z),...,\overline{J_{F_{K}}}(z))=\widehat{F}_{j}(z, \overline{J_{F_{1}}}(z),...,\overline{J_{F_{m}}}(z))~\text{for}~z \in U_{0}.$$
In the following, for simplicity of notation,  we also write for $j
\leq K,$
$$\widehat{h}_{j}(z,\overline{J_{F_{1}}}(z),...,\overline{J_{F_{K}}}(z))=
\overline{J_{F_{j}}}(z)~\text{or}~\widehat{h}_{j}(z,X_{1},...,X_{K})=X_{j}.$$

Now we replace $\overline{F_{j}}(\xi)$ by $\widehat{F_{j}}(\xi,\overline{J_{F_{1}}}(\xi),...,\overline{J_{F_{K}}}(\xi)),$ and replace $\overline{J_{F_{j}}}(\xi)$ by $\widehat{h}_{j}(\xi,\overline{J_{F_{1}}}(\xi),...,\overline{J_{F_{K}}}(\xi)),$ for
$1 \leq j \leq m,$ in $\eqref{unicomplexify}.$
Furthermore, we write $X=(X_{1},...,X_{K}),$ and replace
$\overline{J_{F_{j}}}(\xi)$ by $X_{j}$ for $1 \leq j \leq K$ in
$$\widehat{F_{j}}(\xi,\overline{J_{F_{1}}}(\xi),...,\overline{J_{F_{K}}}(\xi)),
\widehat{h}_{j}(\xi,\overline{J_{F_{1}}}(\xi),...,\overline{J_{F_{K}}}(\xi)),
1 \leq j \leq m.$$ We  define
a new function $\Phi$ as follows:
\begin{equation}
\Phi(z,\xi,X):= \sum_{j=1}^m{\lambda_j}\frac{J_{F_{j}}(z)\widehat{h}_{j}(\xi,X)}{(1+\sum_{i=1}^N\psi_i(F_{j}(z))\psi_i(\widehat F_{j}(\xi,X)))^{\lambda}}-
\frac{1}{(1+\sum_{i=1}^N\psi_i(z)\psi_i(\xi))^{\lambda}}.
\end{equation}

\begin{lemma}\label{HYlemma2}
Shrinking $U$ if necessary, we have $\Phi(z,\xi,X) \equiv 0,$ i.e.,
\begin{equation}\label{equ10}
\sum_{j=1}^m{\lambda_j}\frac{J_{F_{j}}(z)\widehat{h}_{j}(\xi,X)}{(1+\sum_{i=1}^N\psi_i(F_{j}(z))\psi_i(\widehat F_{j}(\xi,X)))^{\lambda}}=
\frac{1}{(1+\sum_{i=1}^N\psi_i(z)\psi_i(\xi))^{\lambda}}.
\end{equation}
or,
\begin{equation}\label{equ11}
\begin{aligned}
(1+\sum_{i=1}^N\psi_i(z)\psi_i(\xi))^{\lambda}& \sum_{j=1}^m \left({\lambda_j}J_{F_{j}}(z)\widehat{h}_{j}(\xi,X)\prod_{1 \leq k \leq m, k \neq j}(1+\sum_{i=1}^N\psi_i(F_{k}(z))\psi_i(\widehat F_{k}(\xi,X)))^{\lambda}\right)\\
&=\prod_{1 \leq j \leq m} (1+\sum_{i=1}^N\psi_i(F_{j}(z))\psi_i(\widehat F_{j}(\xi,X)))^{\lambda}
\end{aligned}
\end{equation}
for  $z \in U$ and $(\xi, X) \in \widehat{U}_{0}.$
\end{lemma}
{\it Proof of Lemma \ref{HYlemma2}}:  Suppose not. Notice $\Phi$ is
Nash algebraic in $(\xi, X)$ for each fixed $z\in U$, by Lemma
\ref{HYlemma}. For a generic fixed $z=z_{0}$ near $0,$ since
$\Phi(z,\xi,X) \not \equiv 0,$ there exist polynomials
$A_{l}(\xi,X)$ for $0 \leq l \leq N$ with $A_{0}(\xi,X) \not \equiv
0$ such that
$$\sum_{0 \leq l \leq N} A_{l}(\xi,X)\Phi^{l}(z,\xi,X) \equiv 0.$$
As $\Phi(z_{0}, \xi,
\overline{J_{F_{1}}}(\xi),...,\overline{J_{F_{K}}}(\xi)) \equiv 0$
for $\xi \in U_{0},$ then it follows that $A_{0}(\xi,
\overline{J_{F_{1}}}(\xi),...,\overline{J_{F_{K}}}(\xi)) \equiv 0$
for $\xi \in U_{0}.$ This is a contradiction to the assumption that
$\{\overline{J_{F_{1}}}(\xi),...,\overline{J_{F_{K}}}(\xi)\}$ is an
algebraic independent set.\  $\endpf$
\medskip

Now that $\widehat{F}_{j}(\xi,X), 1 \leq j \leq m$, is algebraic in
its variables, if $\widehat{F}_{j}, 1 \leq j \leq m$, is independent
of $X,$ then $F_{j}$ is algebraic by Lemma \ref{HYlemma}. This fact
motivates the remaining work in this section.

\subsection{Algebraicity and rationality with uniformly  bounded degree}

In this subsection, we prove the algebraicity and rationality for at
least one of the $F_j's$. We start with the following:

\begin{lemma}\label{analyticity of bad locus}Let $F_j(z),j\in\{1,...,m\}$, be a local holomorphic map defined
 on a neighborhood of $0\in U$ as in $\eqref{unicomplexify}$.  Suppose that there exist $z^0\in U$ and $\xi^0\in
Q_{z^0}$ such that $\Lambda(\beta^1,...,\beta^N)(z^0,\xi^0)$ is well
defined and  non-zero with $\beta^1=(0,0,...,0).$
Then there is an analytic variety $W\subsetneq U$ such that 
when $z\in U\backslash W$, $\Lambda(\beta^1,...,\beta^N)(z,\xi)$ is
a rational function in $\xi$ over $Q_z$ and
    $\Lambda(\beta^1,...,\beta^N)(z,\xi)\not\equiv 0$ on $ Q_z$.
\end{lemma}
{\it Proof of Lemma \ref{analyticity of bad locus}}: By the
assumption, $\frac{\partial \rho}{\partial z_n}(z_0,\xi_0)\not =0$
and
\begin{equation}
\Lambda(\beta^1,...,\beta^N)(z^0,\xi^0)=\left|\begin{matrix}
 \mathcal{L}^{\beta^1}\mathcal F_{j,1} & ...  & \mathcal{L}^{\beta^1}\mathcal F_{j,N} \\
... & ...  & ...  \\
\mathcal{L}^{\beta^N}\mathcal F_{j,1}  & ...  &
\mathcal{L}^{\beta^N}\mathcal F_{j,N}
\end{matrix}\right|(z^ 0,\xi^0)
\end{equation} is   non-zero with $\beta^1=(0,0,...,0).$

 By the definition, ${\mathcal L}_{i}=\frac{\partial}{\partial z_{i}}-
\frac{\frac{\partial \rho}{\partial z_{i}}(z,\xi)}{\frac{\partial
\rho}{\partial z_n}(z,\xi)} \frac{\partial}{\partial z_n}$ and
$\mathcal{L}^{\beta^l}=\mathcal{L}_1^{k^l_1}\mathcal{L}_2^{k^l_2}\mathcal{L}_3^{k^l_3}...\mathcal{L}_{n-1}^{k^l_{n-1}}$
for $\beta^l=(k^l_1,...,k^l_{n-1}),$ $k^l_1,...,k^l_{n-1}$.
 Hence $\Lambda(\beta^1,...,\beta^N)(z,\xi)$ can be written in
the form
$\Lambda(\beta^1,...,\beta^N)(z,\xi)=\frac{\mathcal{G}_1(z,\xi)}{\mathcal{G}_2(z,\xi)}$.
Here $\mathcal{G}_1(z,\xi)=\sum_{|I|=0}^{M_1}\Phi_{I}(z)\xi^I,
\mathcal{G}_2(z,\xi)=\sum_{|J|=0}^{M_2}\Psi_{J}(z)\xi^J,$ with
 $\Phi_{I}$ and $\Psi_{J}$ being holomorphic functions defined over
$U\subset\mathbb{C}^n.$ In fact, $\mathcal{G}_2(z,\xi)$ is simply
taken as a certain sufficiently large power of
$\rho_{z_n}:=\frac{\partial \rho}{\partial z_n}$.

By our assumption, we have $\mathcal{G}_1,\mathcal{G}_2$  not equal
to zero at $(z^0,\xi^0).$ Hence, $\mathcal{G}_1,\mathcal{G}_2$ are
not zero elements in $\mathcal{O}(U)[\xi_1,...,\xi_n]$, the
polynomial ring of $\xi$ with coefficients from the holomorphic
function space over $U$.

By Proposition (III), the defining function of the Segre family
$\rho$ can be written  in the form
$\rho(z,\xi)=\sum_{|\a|=0}^{M_3}\Theta_k(z)\xi^\a,$ which is an
irreducible polynomial in $(z,\xi).$ And for each fixed $z$,  by
Proposition (III), we also have $\rho(z,\xi)$   irreducible as a
polynomial of $\xi$ only.

Then the set of  $z\in U$ where
$\Lambda(\beta^1,...,\beta^N)(z,\xi)$ is  undefined over  $Q_{z}$ is
a subset of  $z\in U$ where $\mathcal{G}_2(z,\xi)$, as a polynomial
of $\xi$, contains the factor $\rho(z,\xi)$ as a polynomial in
$\xi$.
We denote the latter set by $W_2$. Similarly,
 the set of $z\in
U$  with $\Lambda(\beta^1,...,\beta^N)(z,\xi)\equiv 0$ over $ Q_{z}$
is a subset of $z\in U$ where $\mathcal{G}_1(z,\xi)$, as a
polynomial of $\xi$, contains a factor $\rho(z,\xi)$, which we
denote by $W_1$.

Notice that $\rho(z,\xi)\in\mathcal{O}(U)[\xi_1,...,\xi_n]$ depends
on each $\xi_j$ for $1\le j\le n$.
 Also notice that
$\mathcal{G}_2(z,\xi)$, as a certain power of $\rho_{z_n}(z,\xi)$,
depends on  $\xi_n$.

We next characterize $W_2$ by the resultant $R_2$ of
$\mathcal{G}_2(z,\xi)$ and $\rho(z,\xi)$ as polynomials in $\xi_n$.
We rewrite  $\mathcal{G}_2$ and $\rho$ as polynomials of $\xi_n$ as
follows:
$$\mathcal{G}_2=\sum_{i=0}^ka_i(z,\xi_1,...,\xi_{n-1})\xi_n^i,\ \rho=\sum_{j=0}^lb_j(z,\xi_1,...,\xi_{n-1})\xi_n^j.$$
Here  the leading terms $a_k,b_l\not\equiv 0$ with $k,l\geq 1.$
We  write the resultant as
$R_2(z,\xi_1,...\xi_{n-1})=\sum_{I}c_I(z){\xi'}^I,$ where $c_I's$
are holomorphic functions of $z\in U$.

For those points $z\in W_2$, $R_2(z,\cdot)\equiv 0$ as a polynomial
of $\xi_1,...,\xi_{n-1}.$ Then $W_2$ is contained in the complex
analytic set $\wt W_2:=\{c_I=0,\forall I\}.$ If $\wt{ W_2}=U$, then
we can find non-zero polynomials
$f,g\in\mathcal{O}(U)[\xi_1,...,\xi_{n-1}][\xi_n]$ such that
$f\rho+g\mathcal{G}_2\equiv 0,$ where the degree of $g$ in $\xi_n$
is less than  the degree of  $\rho$ in $\xi_n$. Hence
$\{\mathcal{G}_2=0\}\cup\{g=0\}\supset\{\rho=0\}\cap
(U\times\mathbb{C}^n). $ Again by the irreducibility of
$\{\rho=0\}\cap (U\times\mathbb{C}^n)$, since $\{g=0\}$ is a thin
set in $\{\rho=0\}\cap (U\times\mathbb{C}^n), \mathcal{G}_2$
vanishes on $\{\rho=0\}\cap (U\times\mathbb{C}^n).$ This contradicts
$\mathcal{G}_2(z^0,\xi^0)\neq 0.$ Hence $W_2\subset\wt W_2$ and $\wt
W_2$ is a proper complex analytic subset of $U.$

By a similar argument,  we  can prove that $W_1$ is contained in
$\wt W_1$
 that is also a proper analytic set of $U$.  Let $W=\wt W_1\cup \wt W_2.$
Then when $z\in U\backslash W$, $\Lambda(\beta^1,...,\beta^N)(z,\xi)$
is well-defined over $Q_{z}$ as a rational function in $\xi$  and
    $\Lambda(\beta^1,...,\beta^N)(z,\xi)\not\equiv 0$ on $ Q_z$.\
$\endpf$

\medskip

\begin{lemma}\label{continuation of hat x}
Let $ \psi(\xi,X)$
be a non-zero Nash-algebraic function in $(\xi,X)
=(\xi_1,...,\xi_n,X_1,...,X_m)\\
\in {\mathbb C}^n\times {\mathbb
C}^m.$ Write $E$ for a proper complex analytic variety of ${\mathbb
C}^n\times {\mathbb C}^m$ that contains the branch locus of $\psi$ and the
zeros of the leading coefficient in the minimal polynomial of $\psi$.
Then there exists a proper analytic set $W_{1}$ in $\mathbb{C}^n$
such that
$$\{\xi|\ \exists X^0,
(\xi,X^0)\not \in E\}
\supset\mathbb{C}^n\backslash W_1.$$

\end{lemma}
{\it Proof of Lemma \ref{continuation of hat x}}:
Since $\psi$ is algebraic, there is an irreducible polynomial
$\Phi(\xi,X;Y) =\sum_{i=0}^{k}\phi_i(\xi,X)Y^i $ such that
$\Phi(\xi,X,\psi(\xi,X))\equiv 0.$
 If $k=1$ then $\psi$ is a rational function and thus  $E$ is just the poles and points of
 indeterminancy. The proof is then obvious and
  we hence assume $k\geq 2$.

Define $\Psi(\xi,X,Y)=\frac{\partial\Phi}{\partial Y}.$ Since $k\geq
2,$ the degree of  $\Psi$ in $Y$ is at least one.  Consider
$\Phi,\Psi$ as polynomials in $Y,$ and write $R(\xi,X)$ for their
resultant. Then
 the branch  locus is contained in $\{(\xi,X)|R(\xi,X)=0\}$.
Notice that $R\not \equiv 0,$ for  $\Phi$ is irreducible.
Write  $R=\sum_{I}r_{I}(\xi)X^I$  with some  $r_{I}\neq 0$.  Write $\phi_k(\xi,X)=\sum\phi_{k,i}(\xi)X^i$ and
$W_1=\{r_I(\xi)=0\ , \forall I\}\cup \{ \phi_{k,i}(\xi)=0\ , \forall
\ i\}$, which is a proper complex analytic set in $\mathbb{C}^n$.
Then $\{\xi|\ \exists X^0, (\xi,X^0)\not \in E\}
\supset\mathbb{C}^n\backslash W_1.$\ \ $\endpf$

\medskip
Let $E$ be a proper complex analytic variety containing the union of
the branch loci of $\widehat h_j,\widehat F_{j}$ for $j=1,\cdots,m$
and the zeros of the leading coefficients in their minimal
polynomials. For any point $(z^0,\xi^0,X^0)\in
U\times((\mathbb{C}^n\times\mathbb{C}^K)\backslash E),$ we can find
a smooth Jordan curve $\gamma$ in
$U\times((\mathbb{C}^n\times\mathbb{C}^K)\backslash E)$ connecting
$(z^0,\xi^0,X^0)$ with a certain point in $U\times
(\widehat{U}_{0}\sm E).$ We can holomorphically continue the
following  equation along $\gamma$:
\begin{equation}\label{proequation}
\begin{aligned}
(
\rho(z,\xi))^{\lambda}& \sum_{j=1}^m \left( {\lambda_j}J_{F_{j}}(z)\widehat{h}_{j}(\xi,X)\prod_{1 \leq k \leq m, k \neq j}(1+\sum_{i=1}^N\psi_i(F_{k}(z))\psi_i(\widehat F_{k}(\xi,X)))^{\lambda}\right)\\
&=\prod_{1 \leq j \leq m} (1+\sum_{i=1}^N\psi_i(F_{j}(z))\psi_i(\widehat F_{j}(\xi,X)))^{\lambda},\ \ \ \ z \in U,\ (\xi, X) \in \widehat{U}_{0},
\end{aligned}
\end{equation}
to a neighborhood of $(z^0,\xi^0,X^0)$.
For our later discussions, we further define
$$\mathcal{M}_{\text {sing,z}}=\{(z,\xi): \frac{\partial
\rho}{\partial z_{j}}=0,\forall j \},\mathcal{M}_{\text{reg,z}}=\mathcal{M}\backslash\mathcal{M}_{\rm{sing,z }} ;$$
$$\mathcal{M}_{\rm{SING}}=\{(z,\xi): \frac{\partial
\rho}{\partial \xi_{j}}=0,\forall j \}\cup\{(z,\xi): \frac{\partial
\rho}{\partial z_{j}}=0,\forall j \},\ \ \
\mathcal{M}_{\rm{REG}}=\mathcal{M}\backslash\mathcal{M}_{\rm{SING}} ;$$
 $\hbox{Pr}_z: \mathbb C^{2n}\rightarrow\mathbb C^n\ \ \ (z,\xi)\mapsto(z) \text{ and Pr}
_{\xi}:\mathbb C^{2n}\rightarrow\mathbb C^n\ \ \
(z,\xi)\mapsto(\xi).$

Notice that $\mathcal{M}_{\rm{REG}}$ is a Zariski open subset of
$\mathcal{M}$ and  the restrictions of Pr$_z$,Pr$_{\xi}$ to
$\mathcal{M}_{\rm{REG}}$ are open mappings. Also,
 for
$(z^0,\xi^0)\in \mathcal{M}_{\rm{REG}},$ $Q_{z^0}$ is smooth at
$\xi^0,$ and $Q_{\xi^0}$ is smooth at $z^0.$ By Proposition (III),
 $\mathcal{M}_{\rm{reg,z}}\cap (Q_{\xi^0},\xi^0)$
is Zariski open in $(Q_{\xi^0},\xi^0)$.

\begin{lemma}\label{verygeneral}
With the notations we have set up so far, there exists a point
$(z^0,\xi^0,X^0)\in(U\times\mathbb C^n\times
 \mathbb C^K)$ with
$(z^0,\xi^0)\in\mathcal{M}_{\rm{REG}}\cap (U\times\mathbb C^n)$ and
$(\xi^0,X^0)\not\in E$. Moreover, for each $j=1,...,m$, we can find
$\beta_j^1,...,\beta_j^N$ \text{ with }$\beta_j^1=(0,...,0)$ \text{
such that } $\Lambda_{F_j}(\beta_j^1,...,\beta_j^N)(z^0,\xi^0)\neq
0$.
\end{lemma}
{\it Proof of Lemma \ref{verygeneral}}: This is  an easy consequence
of  Propositions (I) (III), Lemma $\ref{analyticity of bad locus}$ and
the Zariski openness of ${\mathcal M}_{\rm REG}$ in ${\mathcal M}$.
\ \ $\endpf$

\medskip
Let  $(z^0,\xi^0,X^0)$ be chosen as in Lemma \ref{verygeneral}. We
then analytically continue the equation $\eqref{proequation}$ to a
neighborhood of the point $(z^0,\xi^0,X^0)$
through a Jordan curve $\gamma$ described above. We denote one of
such neighborhoods by $V_1\times V_2\times V_3,$ where $V_1, V_2$
and $V_3$ are chosen to be  a small neighborhood of $z^0,\xi^0$, and
$X^0$, respectively. It is clear, after shrinking $V_1,V_2,V_3$ if
needed, that there exists a $j_0\in\{1,...,m\}$ such that
$$1+\sum_{i=1}^N\psi_i(F_{j_0}(z))\psi_i(\widehat
{F_{j_0}}(\xi,X))=0,\ \ \text {for }(z,\xi)\in {\mathcal M}\cap
(V_1\times V_2), X\in V_3,
$$
We next proceed to prove the algbraicity for $F_{j_0}(z).$
\begin{theorem}\label{algebraicity}
 $\wh{F_{j_0}}(\xi,X)$, for $\xi\in V_2,\ X\in V_3$, is independent of $X$ and is thus a
  Nash algebraic function of $\xi$.
  Hence $F_{j_0}$ is an algebraic function of $z$. Moreover, the algebraic
total degree of $\wh{F_{j_0}}(\xi,X)=\Ol{F_{j_0}}(\xi)$, and thus of
$F_{j_0}(z)$, is uniformly bounded by a constant  depending only on
the manifold $(X,\omega)$ and the described canonical embedding.
\end{theorem}

Before  proceeding  to the proof, we state  a slightly modified
version of a classical result of Hurwitz. We first give the
following definition:
\begin{definition}Suppose $F$ is an algebraic function defined on $\xi\in\mathbb{C}^n$.
The total degree of $F$ is defined to be   the total degree of its
minimum polynomial. Namely, let $P(\xi;X)$ be an irreducible minimum
polynomial of $F$, the total degree of $F$ is defined as the degree
of $P(\xi;X)$ as a polynomial in $(\xi,X)$.
\end{definition}
We next state some simple facts  about algebraic functions, whose
proof is more or less standard. (See, for instance,\ \cite{Fa}):

\begin{lemma}\label{james-002}
1.\ Suppose $\phi_1,\phi_2$ are algebraic functions defined in
$\xi\in U\subset\mathbb{C}^n$ with total degree bounded by $N$. Then
$\phi_1\pm\phi_2,\phi_1\phi_2,1/\phi_1$ (if $\phi_1\not\equiv 0$)
are algebraic functions and their degrees are bounded above by a
constant  depending only on $N, \ n$.

 2. Suppose
$\phi_1(z_1,...,z_n)$ is an algebraic function of total degree
bounded by $N$, and suppose that $\psi_1(\xi_1,...,\xi_m),...,$
$\psi_n(\xi_1,...,\xi_m)$ are algebraic functions with total degree
bounded by $N$ as well. Let
$$A_0=(\xi_1^0,\xi_2^0,...,\xi_m^0)\in\mathbb{C}^m,$$ where
$\psi_1,...,$ $\psi_n$ are holomorphic functions in a neighborhood
of $A_0$ and let $\phi_1$ be a holomorphic function in a
neighborhood $U\subset\mathbb{C}^n$ of
$(\psi_1(A_0),\psi_2(A_0),...,\
 \ \psi_m(A_0))$.  Then the composition
$\Phi(\xi_1,...,\xi_m)=\phi_1(\psi_1(\xi_1,...,\xi_m),\psi_2(\xi_1,...,\xi_m),\psi_3(\xi_1,...,\xi_m),...,\psi_n(\xi_1,...,\xi_m))$
is an algebraic function with total degree bounded by a constant
$C(N,n,m)$ depending only on $(N,\ n,\ m)$.

\noindent 3. Suppose
$P_1(z_1,z_2,...,z_m,\xi_1,\xi_2,...,\xi_n),...,P_n(z_1,z_2,...,z_m,\xi_1,\xi_2,...,\xi_n)$
are algebraic functions with total degrees bounded from above by N
which   are holomorphic in a neighborhood $U\times
V\subset\mathbb{C}^m\times\mathbb{C}^n$ of
$A_0=(z_1^0,...,z_m^0,\xi_1^0,...,\xi_n^0)$. Suppose that
\begin{displaymath}\begin{cases}P_1(z_1,z_2,...,z_m,\xi_1,...,\xi_n)=0\\P_2(z_1,z_2,...,z_m,\xi_1,...,\xi_n)=0\\...\\P_n(z_1,z_2,...,z_m,\xi_1,...,\xi_n)=0\
\end{cases}
\end{displaymath}
has a solution at
$A_0=(z^0,\xi^0)=(z_1^0,...,z_m^0,\xi_1^0,...,\xi_n^0)$ and
$\frac{\partial(P_1,P_2,...,P_n)}{\partial(\xi_1,\xi_2,...,\xi_n)}(z_1^0,z_2^0,...,z_m^0,\xi_1^0,...,\xi_n^0)\neq
0.$  Then we can solve
$\xi_1=\phi_1(z_1,z_2,...,z_m)$,$\xi_2=\phi_2(z_1,z_2,...,z_m)$,...,$\xi_n=\phi_n(z_1,z_2,...,z_m)$
with $\phi_j(z^0)=\xi^0$  in a neighborhood of
$z^0\in\tilde{U}\subset U\subset\mathbb{C}^m$, where
$\phi_1,...,\phi_n$ are algebraic functions with total degree
bounded by $C(N,n,m).$

\end{lemma}
We now state the following modified  version of  the classical
Hurwitz theorem with a controlled total degree \cite{BM}.
\begin{theorem} \label {james-003} Let ${F}(s,t,\xi_1,\xi_2,...,\xi_m)$ be holomorphic  over
 $U\times V\times W\subset\mathbb{C}^{m+2}$. Suppose that
 for any fixed $s\in U\subset {\mathbb C}$, ${F}$ is an algebraic function in $(t,\xi_1,...,\xi_m)$ with its total
 degree uniformly
 bounded by N; and for any  fixed $t\in V\subset {\mathbb C}$, ${F}$ is an algebraic function of
 $(s,\xi_1,...,\xi_m)$
  with its total degree uniformly  bounded by $N$.
  Then ${F}$ is an algebraic function with total degree bounded by a constant depending only on $(m,N)$.
\end{theorem}
The proof of Theorem \ref{james-003} is more or less the same as in
the classical setting \cite{BM}. (See, for example, the Ph. D. thesis of
the first author \cite{Fa})

\medskip

{\it Proof of Theorem \ref{algebraicity}}: By the choice of $(z^0,\xi^0,X^0)$,
 \ there exist $\beta_{j_0}^1,...,\beta_{j_0}^N$ such that
\begin{equation}
\Lambda_{F_{j_0}}(\beta_{j_0}^1,...,\beta_{j_0}^N)(z^0,\xi^0)=\left|\begin{matrix}
 \mathcal{L}^{\beta_{j_0}^1}\mathcal F_{j_0,1} & ...  & \mathcal{L}^{\beta_{j_0}^N}\mathcal F_{j_0,N} \\
... & ...  & ...  \\
\mathcal{L}^{\beta_{j_0}^N}\mathcal F_{j_0,1}  & ...  & \mathcal{L}^{\beta_{j_0}^N}\mathcal F_{j_0,N}
\end{matrix}\right|(z^0,\xi^0)\neq 0.
\end{equation}

We can also assume that $(z_0,\xi_0)$ satisfies the  assumption in
Proposition (II) after a slight perturbation of $z_0$ inside
$Q_{\xi_0}$  if needed. By Proposition (II), we can find $z^1\in
V_1\cap Q_{\xi^0}$ such that $Q_{z^0}$
 intersects $ Q_{z^1}$ transversally at $\xi^0.$  Moreover  there exists a neighborhood
 $B$ of $\xi^0$ and a biholomorphic parametrization of $B:$
  $(\xi_1,\xi_2,...,\xi_n)=\mathcal{G}(\tilde\xi_1,\tilde\xi_2,...,\tilde\xi_n)$
  with
$(\tilde\xi_1,\tilde\xi_2,...,\tilde\xi_n)\in U_1\times
U_2\times...\times U_n\subset\mathbb C^n$. Here $U_1, U_2$ are
as in Proposition (II).
Moreover, $\mathcal{G}(\{\tilde\xi_1=1\}\times U_2\times...\times
U_n)\subset Q_{z^0},\mathcal{G}(U_1\times\{\tilde\xi_2=1\}\times
U_3\times...\times U_n)\subset Q_{z^1}.$ Also, for   $s\in U_1,t\in
U_2$, $\mathcal{G}(\{\tilde\xi_0=t\}\times U_2\times...\times
U_n),\mathcal{G}(U_1\times\{\tilde\xi_1=s\}\times U_3\times...\times
U_n)$  are open pieces of  certain Segre varieties. Here
$\mathcal{G}$ consists of algebraic functions with total algebraic
degree uniformly bounded by $M$ and the canonical embedding.
Consider the equation:
\begin{equation}\label{j0equation}
1+{\mathcal F}_{j_0}(z)\cdot
\widehat{{\mathcal F}_{j_0}}(\xi,X)=0,\ (z,\xi,X)\in V_1\times
V_2\times V_3, (z,\xi)\in {\mathcal M}.
\end{equation}
 Since the holomorphic  vector fields
$\{\mathcal{L}_i\}_{i=1}^{n-1}$ are tangent to the Segre family, we
have
\begin{equation}\label{linear system}
\left(\begin{matrix}
 \mathcal{L}^{\beta_{j_0}^1}\mathcal F_{j_0,1}(z,\xi) & ...  & \mathcal{L}^{\beta_{j_0}^1}\mathcal
  F_{j_0,N}(z,\xi) \\
... & ...  & ...  \\
\mathcal{L}^{\beta_{j_0}^N}\mathcal F_{j_0,1}(z,\xi)  & ...  &
\mathcal{L}^{\beta_{j_0}^N}\mathcal F_{j_0,N}(z,\xi)
\end{matrix}\right)
\left(\begin{matrix}\widehat {\mathcal{F}_{j_0,1}}(\xi,X)\\
...\\
\widehat {\mathcal{F}_{j_0,N}}(\xi,X)

\end{matrix}\right)=
\left(
\begin{matrix}-1\\
\cdots\\
0
\end{matrix}
\right),
\end{equation}
where $ (z,\xi)(\approx (z^0, \xi^0))\in {\mathcal M},\ X\approx
X^0.$

By the Cramer's rule, we conclude that
$\{\widehat{\mathcal{F}_{j_0,l}}(\xi,X)\}_{l=1}^N$ are rational
functions of $\xi$ with a uniformly bounded  degree on an open piece
of each Segre variety $Q_z$ for $z\approx z^0$.
By the previous   modified  Hurwitz Theorem (Theorem \ref{james-003}), we conclude the
algebraicity of $\widehat{\mathcal{F}_{j_0,l}}(\xi,X)$ for
$l=1,...,N.$. Since in $\eqref{linear system}$ the matrix $\left(  \mathcal{L}^{\beta_{j_0}^{\mu}}\mathcal F_{j_0,\nu}(z,\xi)    \right)_{1 \leq \mu, \nu \leq N}$ and the
right hand side are independent of $X$, these functions must also be
independent of the $X$-variables. Moreover,
by Lemma $\ref{james-002}$ and Theorem $\ref{james-003}$, the total
algebraic degree of
$\Ol{F}_{j_0,l}(\xi)=\widehat{\mathcal{F}_{j_0,l}}(\xi,X)$, for
$l=1,...,n$, is  uniformly bounded.  Since $\Ol{F}$ is obtained by
holomorphically continuing the conjugation function $\Ol{F}$ of $F$,
we conclude the  algebraicity  of $F_{j_0,l}$ for each $1\le l\le
n$. Also  the total algebraic degree of each $F_{j_0,l}$ is  bounded
by a constant  depending only on $(M,\omega)$. \ \ $\endpf$

\medskip
\begin{theorem}\label{rationality} Under the notations we have just set up,
 $F_{j_0}$ is a rational map,
 whose degree   depends only  on the canonical embedding $M\hookrightarrow\mathbb{CP}^N$.
\end{theorem}
For the proof Theorem $\ref{rationality}$, we first recall the
following Lemma of  [HZ]:
\begin{lemma}\label{LemmaHZ} (Lemma 3.7 in [HZ])
Let $U \subset \mathbb{C}^{n}$ be a simply connected open subset and
$\mathcal{S} \subset U$ be a closed complex analytic subset of
codimension one. Then for $p \in U \setminus \mathcal{S},$ the
fundamental group $\pi_{1}(U \setminus \mathcal{S},p)$ is generated
by loops obtained by concatenating (Jordan) paths
$\gamma_{1},\gamma_{2},\gamma_{3},$ where $\gamma_{1}$ connects $p$
with a point arbitrarily close to a smooth point $q_{0} \in
\mathcal{S},$ $\gamma_{2}$ is a loop around $\mathcal{S}$ near
$q_{0}$ and $\gamma_{3}$ is $\gamma_{1}$ reversed.
\end{lemma}

{\it Proof of Theorem \ref{rationality}}: We give a proof for the
rationality of  $F_{j_0}$. Once this is done, we then conclude that
the degree of $F_{j_0}$ is
 uniformly bounded,
 for we
know the total algebraic degree of $F$ is uniformly bounded by
Theorem $\ref{algebraicity}$.

Suppose that  $F_{j_0}$ and thus $\Ol {F_{j_0}}$   is not rational.
Write  $E\subset {\mathbb C}^n$ for a proper complex analytic
variety containing the branch locus of $F_{j_0}, \ov {F_{j_0}}$ and
the zeros of the leading coefficients of the minimal polynomials of
their components. We first notice that for $A\not =B\in {\mathbb
C}^n,\ Q^*_A\not= Q^*_B$, by Lemma \ref{unique}. Hence, for any
proper complex analytic variety $V^1,V^2\subset {\mathbb C}^n$ and
any point $(a,b)\in {\mathcal M}$, we can find $(a^1,b^1)\approx
(a,b)$ such that $a^1\in Q_{b^1}\sm V^1$ and $b^1\not \in V^2$.

 We choose $(z^0,\xi^0)$ as
above and assume further  that $z^0,\xi^0\not \in E$ (after a small
perturbation if needed). We choose a sufficiently small neighborhood
$W$ of $(z^0,\xi^0)$ in ${\mathcal M}_{\text{REG}}$ such that  for
each $(z^1,\xi^1)\in W$, we can find, by Lemma \ref{LemmaHZ}, a loop of the form
$\gamma=\gamma_{1}^{-1}\circ\gamma_{2}\circ\gamma_{1}$ in ${\mathbb
C}^n\sm E$ with $\gamma(0)=\gamma(1)=\xi^1,\gamma_1(1)=q$. Here
$\gamma_1$ is a simple curve connecting $\xi^1$ to $q$ with $q$ in a
small ball $B_p$ centered at a certain smooth point $p$ of $E$
 such that the fundamental
group of $B_p\sm E$ is generated by $\gamma_2$; and $\gamma_1^{-1}$
is the reverse curve of $\gamma_1$. Moreover, when $\Ol {F_{j_0}}$
is holomorphically continued along $\gamma$, we end up with a
different branch $\Ol {F_{j_0}}^{*}$ of $\Ol {F_{j_0}}$ near
$\xi^1$. We pick $p$ such that there is an $X_p \not \in E$ with
$(X_p, p)\in {\mathcal M}_{\text{reg,z}}$. (This  follows from
Proposition (III) and Lemma \ref{unique} as mentioned above.)
Take a certain small neighborhood ${\mathcal W}$ of $(X_p,p)$ in
${\mathcal M}_{\text{reg,z}}.$
We assume, without loss of generality, that the piece $\mathcal{W}$ of  ${\mathcal M}_{\text{reg,z}}$ is defined by a
holomorphic function of the form $z_1=\phi(z_2,\cdots,z_n,\xi)$. In particular, writing $X_p=(z_1^p, \cdots, z_n^p),$ we have $z_1^p=\phi(z_2^p, \cdots, z_n^p, p).$
Make $B_p$ sufficiently small such that it is  compactly contained
in the image of the projection of ${\mathcal W}$ into the
$\xi$-space. Write
$X_q=(\phi(z_2^p,\cdots,z^p_n,q),z_2^p,\cdots,z_n^p)$ and
 define the loop
 $\gamma_2^*(t)=(\phi(z_2^p,\cdots,z_n^p,\gamma_2(t)),z^p_2,\cdots,z^p_n)$.
 Then $\gamma_2^*$ is a loop whose  base point is at $X_q$. Also, we
 have
 $(\gamma^*_2(t),\gamma_2(t))\in {\mathcal M}$.

Notice that $X_p\not\in E$. After shrinking $B_p$ if needed, we
assume that $\gamma_2^*$ stays sufficiently close to $X_p$ and  is homotopically trivial in ${\mathbb
C}^n\sm E$.

Now we slightly thicken $\gamma_1$ to get a simply connected domain
$U_{1}$ of $\mathbb{C}^n \sm E$. Since ${\mathcal M}$ is irreducible
over ${\mathbb C}^n\times U_{1}$, we can find a smooth simple curve
$\wt{\gamma}_1=({\gamma_1}^*,\wh\gamma_1)$ in $\mathcal{M}\sm
((E\times\mathbb C^n)\cup(\mathbb C^n\times E))$ connecting
$(z^1,\xi^1)$ to $(X_q,q)$. Then $\wh{\gamma_1}$ is homotopic to
${\gamma_1}$ relatively to $\{\xi^1,q\}$ and ${\gamma_1}^*(1)=X_q$.
Now replace $\gamma$ by its homotopically equivalent loop
${\wh\gamma}_1^{-1}\circ \gamma_2\circ {\wh\gamma}_1$ and define
$\gamma^*={ \gamma^*}_1^{-1}\circ \gamma^*_2\circ {\gamma^*}_1$.
Define $\Gamma=(\gamma^*, \gamma)$. Then the image of $\Gamma$ lies inside ${\mathcal
M}\sm ((E\times\mathbb C^n)\cup(\mathbb C^n\times E))$. Continuing
Equation $\eqref{j0equation}$ along $\Gamma$ and noticing that it is
independent of $X$ now, we get both
$$ 1+{\mathcal F}_{j_0}(z)\cdot \Ol{\mathcal
F}_{j_0}(\xi)=0\ \rm{and }\ 1+{\mathcal F}_{j_0}(z)\cdot
\Ol{\mathcal F}^*_{j_0}(\xi)=0\ \forall (z,\xi)\in {\mathcal M}\cap
((V_1\sm E)\times
(V_2\sm E)).$$
Now as before, applying the uniqueness for the solution
of the linear system \eqref{linear system} (with an invertible
coefficient matrix), we then conclude that ${\Ol {F_{j_0}}}^*\equiv
{\Ol {F_{j_0}}}$. This is a contradiction. \ \ $\endpf$

\subsection{Isometric extension of  $F$}

For simplicity of notation, in the rest of this section, we denote
the map $F_{j_0}$ just by $F$.  Now that all components of $F$ are
rational functions, it is easy to verify that $F$ gives rise to a
rational map $M\dashrightarrow M$. By the Hironaka theorem (see [H]
and [K]), we have a (connected) complex manifold $Y$ of the same
dimension, holomorphic maps $\tau: Y\ra M$, $\sigma: Y\ra M$, and a
proper complex analytic variety $E_1$ of $M$ such that $\sigma:\
Y\sm \sigma^{-1}(E_1)\ra M\sm E_1$ is biholomorphic; $F: M\sm E_1\ra
M$ is well-defined; and for any $p\in Y\sm \sigma^{-1}(E_1)$,
$F(\sigma(p))=\tau(p)$.

Let $E_2$ be a proper complex analytic subvariety of $M$ containing
$E_1$, $M\sm {\mathcal A}$ and let $E_3\subset Y$ be the proper
subvariety where $\tau$ fails to be biholomorphic.
 Write
$E^*=\tau(\sigma^{-1}(E_2)\cup E_3)\cup (M\sm {\mathcal A})$ and
$E=\sigma(\tau^{-1}(E^*))$. Then $F: {\mathcal A}\sm E\ra {\mathcal
A}\sm E^*$ is a holomorphic covering map. We first prove

\begin{lemma}\label{biho}: Under the above notation,  $F: {\mathcal A}\sm E\ra {\mathcal
A}\sm E^*$ is a biholomorphic map.
\end{lemma}
{\it Proof of Lemma \ref{biho}}: We first notice that since $F$ is
biholomorphic near $0$ with $F(0)=0$. We can assume that $0\not\in
E$. Consider $F^2=F\circ F$. Then $\Ol{F^2}=\Ol{F}^2$. Since
$(F,\Ol{F})$ maps ${\mathcal M}$ into ${\mathcal M}$ whenever it is
defined, it is easy to see that $(F,\Ol{F})\circ
(F,\Ol{F})=(F^2,\Ol{F}^2)$ also maps ${\mathcal M}$ into ${\mathcal
M}$ at the points where it is well-defined. Hence, we can repeat a
similar  argument for $F$ to conclude that $F^2$, as a rational map,
also has its degree bounded by a constant independent of $F^2$.
Similarly, we can
conclude that for any positive integer $m$, $F^m$ is a rational
 map with degree bounded by a constant independent of $m$ and $F$.
Now, as for $F$, we can find complex anaytic subvarieties $E^{(m)},
\ E^{*(m)}$ of ${\mathbb C}^n$ such that $F^m$ is a holomorphic
covering map from ${\mathcal A}\sm E^{(m)}\ra {\mathcal A}\sm
E^{*(m)}$. Suppose $F:  {\mathcal A}\sm E\ra {\mathcal A}\sm E^*$ is
a $k$ to $1$ covering map. It is easy to see that  $F^m:\ {\mathcal
A}\sm E^{(m)}\ra {\mathcal A}\sm E^{*(m)}$ is a $k^m$ to $1$
covering map. However, since the degree $F^m$ is independent of $m$,
we conclude that $k=1$ by the following Bezout theorem:
\begin{theorem}{\rm{(\cite{S})}} The number of isolated solutions to a system of polynomial equations
$$f_1(x_1,...,x_n)=f_2(x_1,...,x_n)=...=f_n(x_1,...,x_n)=0$$
is bounded by $d_1d_2\cdot\cdot\cdot d_n$, where $d_i=\deg f_i.$
\end{theorem}
This proves the lemma. $\endpf$

\medskip

Now we prove that  $F$ extends to  a global holomorphic isometry of
$(M,\omega)$.
\begin{theorem}\label{isom} $F: (U,\omega|_{U})\ra (M,\omega)$ extends to  a global holomorphic  isometry of $(M,\omega)$.
\end{theorem}
{\it Proof of Theorem \ref{isom}}: By what we just achieved,
we then have two proper
complex analytic varieties $W_1,\ W_2$ of ${\mathbb C}^n$ such that
$F:{\mathbb C}^n\backslash W_1\rightarrow {\mathbb C}^n\backslash
W_2$ is biholomorphic. Similarly we have two proper complex analytic
subvarieties $W_1^*,\ W_2^*$ of ${\mathbb C}^n$ such that $\Ol{F}:
{\mathbb C}^n\sm W^*_1 \rightarrow  {\mathbb C}^n\sm  W_2^*$ is a
biholomorphic map. Hence
$$\mathfrak{F}=(F,\Ol{F}):
{\mathbb C}^n\sm W_1\times{\mathbb C}^n\sm W_1^*\rightarrow {\mathbb
C}^n\sm W_2\times{\mathbb C}^n\sm W_2^*$$ is  biholomorphic.
Let $\rho$ be the defining function of the Segre family as described
before. Since $\rho$ is irreducible as a polynomial in $(z,\xi)$,
$\mathcal M$ is an irreducible complex analytic variety of
${\mathcal A}$. Since $\mathfrak{F}$ maps a certain  open piece of
${\mathcal M}$ into an open piece of ${\mathcal M}$, by the
uniqueness of holomorphic functions, we see that
$\mathfrak{F}=(F,\Ol{F})$ also gives a biholomorphic map from
$({\mathbb C}^n\sm W_1\times{\mathbb C}^n\sm W_1^*)\cap{\mathcal M}$
to $ ({\mathbb C}^n\sm W_2\times{\mathbb C}^n\sm W_2^*)\cap{\mathcal
M}$. Hence $\rho_F=\rho(F(z),\Ol{F}(\xi))$ defines the same  subvariety as
$\rho$ does over ${\mathbb C}^n\sm W_1\times{\mathbb C}^n\sm W_1^*$.
Since $\rho_F$ is a rational function in $(z,\xi)$ with denominator
coming from the factors of the denominators of $F(z)$ and
$\Ol{F}(\xi)$, we can write
\begin{equation}\label{polysegre}
\rho_F(z,\xi)=(\rho(z,\xi))^l\frac{
P_1^{i_1}(z,\xi)P_2^{i_2}(z,\xi)\cdot\cdot\cdot
P_{\tau}^{i_{\tau}}(z,\xi)}{Q_1^{j_1}(z)\cdot\cdot\cdot
Q_{\mu}^{j_{\mu}}(z)R_1^{k_1}(\xi)\cdot\cdot\cdot
R_{\nu}^{k_{\nu}}(\xi)}
\end{equation}
 Here the zeros of $Q_j(z)$ and $R_j(\xi)$
stay in $W_1$ and $W_1^*$, respectively. All polynomials are
irreducible and prime to each other. By what we just mentioned
$P_j(z,\xi)$ can not have any zeros in ${\mathbb C}^n\sm
W_1\times{\mathbb C}^n\sm W_1^*$, for otherwise it must have $\rho$
as its factor by the irreducibility of $\rho$. Hence the zeros of
$P_j(z,\xi)$ must stay in $ (W_1\times{\mathbb C}^n)\cup ({\mathbb
C}^n\times W_1^*)$. From this, it follows easily that
$P_j(z,\xi)=P_{j,1}(z)$ or $P_j(z,\xi)=P_{j,2}(\xi)$. Namely,
$P_j(z,\xi)$ depends  either on $z$ or on $\xi$. Since
$\mathfrak{F}$ is  biholomorphic, we see that $l=1$. Thus replacing
$\xi$ by $\bar z$ and taking
 $i\partial\bar{\partial}\log$ to $\eqref{polysegre}$,
 we have $i\partial\bar{\partial}\log \rho_F(z,\bar{z})=i\partial\bar{\partial}\log \rho(z,\bar{z})$.
This shows that $F^*(\omega)=\omega$, or $F$ is a local isometry.
Now, by the Calabi Theorem (see \cite{Ca}), $F$ extends to a global holomorphic
isometry of $(M,\omega)$. This proves Theorem $\ref{isom}$. $\endpf$

\bigskip

We now are ready to give a proof of  Theorem $\ref{main}$. By what
we have obtained, there is a component $F_j$  for $F$ in Theorem
\ref{main} that extends to a holomorphic isometry to $(M,\omega)$.
Hence $F_j^*(d\mu)=d\mu$. Notice $\lambda_j<1$ due to the positivity
of all terms in the right hand side of the equation (\ref{volume
preserving equation}). After a cancellation, we reduce the theorem to
the case with only $(m-1)$- maps. Then by an induction argument, we
complete the proof of Theorem \ref{main}.  \ $\endpf$

\section{Partial non-degeneracy: Proof of Proposition (I)}
In this  section, we prove Proposition (I) for
irreducible compact Hermitian spaces of compact type. Since the
argument differs as its  type varies, we do it  on a case by case
base. For convenience of the reader, we give a detailed proof
here  for the Grassmannians and Hyperquadrics.  We will  include the
rest arguments in Appendix II.

\subsection{Spaces of type I}

With the same notations that we have set up  in $\S 2,$ $Z$ is a
$p\times q$ matrix ($p\leq q$); $Z(\begin{matrix}
 i_{1} & ... & i_{k} \\
 j_{1} & ... & j_{k}
   \end{matrix}
)$
is the determinant of the submatrix of $Z$ formed by its $i_{1}^{\text{th}},...,i_{k}^{\text{th}}$ rows and $j_{1}^{\text{th}},...,j_{k}^{\text{th}}$ columns;
${z}=(z_{11},...,z_{1q},z_{21},...,z_{2q},...,$ $z_{p1},...,z_{pq})$ is the coordinates of $\mathbb C^{pq}{\cong}\mathcal A\subset G(p,q).$
Let $0\in U$ be a small neighborhood of $0$ in ${\mathbb C}^{pq}$
and $F$ be a biholomorphic map defined over $U$ with $F(0)=0$. For
convenience of our discussions, we represent the map
$F:U\rightarrow{\mathcal A}$  as a holomorphic
 matrix-valued map:
$$F= \left( \begin{matrix}
 f_{11} & ... & f_{1q}  \\
 ... & ...  & ... \\
 f_{p1} & ...  & f_{pq}
\end{matrix}\right).$$
Similar to $Z(\begin{matrix}
 i_{1} & ... & i_{k} \\
 j_{1} & ... & j_{k}
   \end{matrix}
)$,
 $F(\begin{matrix}
 i_{1} & ... & i_{k} \\
 j_{1} & ... & j_{k}
   \end{matrix}
)$ denotes the determinant of  the submatrix formed by the
$i_{1}^{\text{th}},...,i_{k}^{\text{th}}$ rows and
$j_{1}^{\text{th}},...,j_{k}^{\text{th}}$ columns of the matrix $F$.
Recall in  $\eqref{coord1}, r_z$ is defined as 
\begin{displaymath}
(\psi_{1},\psi_{2},...,\psi_{N})=(\cdots,Z(\begin{matrix}
 i_{1} & ... & i_{k} \\
 j_{1} & ... & j_{k}
\end{matrix}),\cdots),1\leq i_{1}<...<i_{k} \leq p, 1 \leq j_{1}<...<j_{k} \leq q,1 \leq k \leq p.
\end{displaymath}
Similarly, we define:
\begin{displaymath}
r_{F}:=(\cdots,F(\begin{matrix}
 i_{1} & ... & i_{k} \\
 j_{1} & ... & j_{k}
\end{matrix}),\cdots),1\leq i_{1}<...<i_{k} \leq p, 1 \leq j_{1}<...<j_{k} \leq q,1 \leq k \leq
p.
\end{displaymath}
Notice that $r_{F}=(\psi_{1}(F(z)),..., \psi_{N}(F(z)))$.
We define $$\tilde{z}:=(z_{11},...,z_{1q},z_{21},...,
z_{2q},...,z_{p1},...,z_{p(q-1)}),$$ i.e., $\wt{z}$ is obtained from
$z$ by dropping the last component $z_{pq}.$ Write
$\frac{\partial^{|\alpha|}}{\partial \tilde{z}^{\alpha}}=
\frac{\partial^{|\alpha|}}{\partial
 z_{11}^{\alpha_{11}}...\partial z_{p(q-1)}^{\alpha_{p(q-1)}}}$ for
any $(pq-1)-$multiindex $\alpha,$ where 
$\alpha=(\alpha_{11},...,\alpha_{1p},\alpha_{21},...,\alpha_{2q},...,
\alpha_{p1},...,\alpha_{p(q-1)}).$


We apply the notion of the  partial degeneracy defined in Definition
\ref{df21} of $\S 3$ by letting $\psi=r_{F}$ and letting
$\widetilde{z}$ be as just defined with $m=pq.$  We  prove the following proposition:
\begin{proposition}\label{thm36}
$r_{F}$ are $\widetilde{z}-$nondegenerate near $0.$ More precisely,
$\mathrm{rank}_{1+N-pq}(r_{F}, \widetilde{z})=N.$
\end{proposition}

{\it Proof of Proposition \ref{thm36}}: If $p=1,q=n \geq 1$ i.e., the Hermitian symmetric space $M=\mathbb{P}^n,$ then it follows from Lemma \ref{22-1} that $\mathrm{rank}_1(r_F, \widetilde{z})=N=n.$ In the following we assume $p \geq 2.$

Suppose the conclusion is not true. Namely, assume
that
$\mathrm{rank}_{1+N-pq}(r_{F}, \widetilde{z})<N.$ Since  the
hypothesis of  Theorem \ref{thm211} is satisfied,  we see that there
exist $c_{pq+1},...,c_{N} \in \mathbb{C}$ which are not all zero
such that
\begin{equation}\label{eq58}
\sum_{i=pq+1}^N c_{i}\psi_{i}(F)(z_{11},...,z_{pq-1},0)\equiv 0.
\end{equation}

The next step is to show that (\ref{eq58}) cannot hold in the
setting of Proposition \ref{thm36}.  This is obvious if we can prove
the following:
\begin{lemma}\label{t1}

Let

$$H= \left( \begin{matrix}
 h_{11} & ...  &  h_{1p}  \\
 ... & ... & ... \\
 h_{p1} & ... & h_{pq}
\end{matrix}\right),  $$
be a vector-valued holomorphic function in a neighborhood $U$ of $0$
in $\tilde{z}=(z_{11},...,z_{p(q-1)}) \in \mathbb{C}^{pq-1}$ with
$H(0)=0.$ Assume that $H$ is of full rank at $0.$ Set
\begin{equation}
(\phi_{1},...,\phi_{m}):=r_{H}=\left(\left(~H\left(\begin{matrix}
 i_{1} & ... & i_{k} \\
 j_{1} & ... & j_{k}
\end{matrix}\right)~\right)_{1\leq i_{1}<...<i_{k} \leq p, 1 \leq j_{1}<...<j_{k} \leq q} \right)_{2 \leq k \leq p}.
\end{equation}
Here
$$m=\left(
      \begin{array}{c}
        p \\
        2 \\
      \end{array}
    \right)\left(
             \begin{array}{c}
               q \\
               2 \\
             \end{array}
           \right)+...+\left(
                         \begin{array}{c}
                           p \\
                           p \\
                         \end{array}
                       \right)\left(
                                \begin{array}{c}
                                  q \\
                                  p \\
                                \end{array}
                              \right).
$$
Let $a_{1},...,a_{m}$ be complex numbers such that
\begin{equation}\label{eq60}
\sum_{i=1}^{m} a_{i}\phi_{i}(\tilde{z})\equiv 0~ \text{for all}~\tilde{z} \in U.
\end{equation}
Then $a_{i}=0$ for each $1 \leq i \leq m.$
\end{lemma}
{\it Proof of Lemma \ref{t1}}:
We start with the simple case $p=q=2,$ in which  $m=1.$ Then by the
assumption (\ref{eq60}), $a_{1}\phi_{1}=0.$ Here
$$\phi_{1}=\left|\begin{matrix}
 h_{11} & h_{12} \\
 h_{21} & h_{22}
                \end{matrix}\right|.
$$
Note that $H=(h_{11},h_{12},h_{21},h_{22})$ is of full rank at $0.$
We assume, without loss of generality, that
$\tilde{H}:=(h_{11},h_{12},h_{21})$ is a local  biholomorphic map
from $\mathbb{C}^3$ to $\mathbb{C}^3.$ After an appropriate
biholomorphic change of coordinates preserving $0$, we can assume
$h_{11}=z_{11}, h_{12}=z_{12},h_{21}=z_{21},$ and still write the
last component as $h_{22}.$ Then we have
$$a_{1}\phi_{1}=a_{1}(z_{11}h_{22}-z_{12}z_{21}) \equiv 0,$$
which easily yields that $a_{1}=0.$

We then prove the lemma for the case of $p=2, q=3$, in which $m=3.$
As before, without loss of generality, we assume that
$\tilde{H}:=(h_{11},h_{12}, h_{13}, h_{21},h_{22})$ is a local biholomorphic map near $0$ from
$\mathbb{C}^5$ to $\mathbb{C}^5.$ After an appropriate biholomorphic
change of coordinates, we assume that
$\tilde{H}=(z_{11},...,z_{22}).$ By  (\ref{eq60}), we have
\begin{equation}
a_{1}\phi_{1}+...+a_{3}\phi_{3}=a_{1}\left| \begin{matrix}
 z_{11} & z_{12}  \\
 z_{21} & z_{22}
\end{matrix}\right|+a_{2}\left| \begin{matrix}
 z_{11} & z_{13}  \\
 z_{21} & h_{23}
\end{matrix}\right| + a_{3}\left| \begin{matrix}
 z_{12} & z_{13}  \\
 z_{22} & h_{23}
\end{matrix}\right|.
\end{equation}

The conclusion can be easily proved by checking the coefficients in
the Taylor expansion at $0$. Indeed, the quadratic terms
$z_{13}z_{21}, z_{13}z_{22}$ only appear once in the last two
determinants. This implies $a_{2}=a_{3}=0.$ Then trivially
$a_{1}=0.$

We also prove the case $p=q=3.$ In this case $m=10.$  As before,
without loss of generality, we assume that
$\tilde{H}=(h_{11},...,h_{32})$ is a biholomorphic map from
$\mathbb{C}^8$ to $\mathbb{C}^8.$ After an appropariate
biholomorphic change of coordinates, we can assume that
$\tilde{H}=(z_{11},...,z_{32}).$ Then by assumption, we have
\begin{equation}
\begin{split}
 a_{1} &\phi_{1}+...+a_{10}\phi_{10}=\\
 a_{1}& \left| \begin{matrix}
 z_{11} & z_{12}  \\
 z_{21} & z_{22}
\end{matrix}\right|+ a_{2}\left| \begin{matrix}
 z_{11} & z_{13}  \\
 z_{21} & z_{23}
\end{matrix}\right|+a_{3}\left| \begin{matrix}
 z_{12} & z_{13}  \\
 z_{22} & z_{23}
\end{matrix}\right|+ a_{4}\left| \begin{matrix}
 z_{11} & z_{12}  \\
 z_{31} & z_{32}
\end{matrix}\right|+ a_{5}\left| \begin{matrix}
 z_{11} & z_{13}  \\
 z_{31} & h_{33}
\end{matrix}\right| +
a_{6}\left|\begin{matrix}
 z_{12} & z_{13}  \\
 z_{32} & h_{33}
\end{matrix}\right|            \\
& + a_{7}\left| \begin{matrix}
 z_{21} & z_{22}  \\
 z_{31} & z_{32}
\end{matrix}\right|+ a_{8}\left| \begin{matrix}
 z_{21} & z_{23}  \\
 z_{31} & h_{33}
\end{matrix}\right|+a_{9}\left| \begin{matrix}
 z_{22} & z_{23}  \\
 z_{32} & h_{33}
\end{matrix}\right| + a_{10}\left| \begin{matrix}
 z_{11} & z_{12} & z_{13}  \\
 z_{21} & z_{22}  & z_{23}  \\
 z_{31} & z_{32}  & h_{33}
\end{matrix}\right| =0.
\end{split}
\end{equation}
We then check the coefficients for each term in its Taylor expansion
at $0$. First it is easy to note that $a_{5}=a_{6}=a_{8}=a_{9}=0$ by
checking the coefficients of quadratic terms
$$z_{13}z_{31},z_{13}z_{32}, z_{23}z_{31},z_{23}z_{32},$$ respectively.
Then by checking the coefficients of other quadratic terms, we see
that $a_{1}=a_{2}=a_{3}=a_{4}=a_{7}=0.$ Finally we check the
coefficient of  the cubic term $z_{13}z_{22}z_{31}$ to obtain that
$a_{10}=0.$

We now prove  the general case: $q \geq p \geq 2.$ As before, we
assume  without loss of generality that
 $\tilde{H}=(h_{11},...,h_{p(q-1)})$ is a biholomorphic map
from $\mathbb{C}^{pq-1}$ to $\mathbb{C}^{pq-1}.$ Furthermore, we
have $\tilde{H}=(z_{11},...,z_{p(q-1)})$ after an appropriate
biholomorphic change of coordinates. We again first consider the
coefficients of the quadratic terms in (\ref{eq60}). For that, we
consider the $2 \times 2$ submatrix involving $h_{pq},$ i.e.,
$H\left( \begin{matrix}
 l & p \\
 k & q
\end{matrix}
 \right), 1 \leq l < p, 1 \leq k < q. $ Note that $z_{lq}z_{pk}$ only appears in this $2 \times 2$ determinant, which
 yields that the coefficient $a_{i}$ associated to this
$2 \times 2$ determinant is $0,$ for any $1 \leq i < p, 1 \leq j <
q.$ Then by checking the coefficients of other quadratic terms, we
see that all coefficients $a_{i}'s$ that are associated to $2 \times
2$ determinants $H\left( \begin{matrix}
 l_{1} & l_{2}  \\
 k_{1} & k_{2}
\end{matrix} \right), 1 \leq l_{1}, l_{2} \leq p, 1 \leq k_{1}, k_{2} \leq q,$ are $0.$

We then consider the coefficients of cubic terms in (\ref{eq60}). We first look at those $3 \times 3$ submatrix involving $h_{pq}.$ i.e., $H\left( \begin{matrix}
 l_{1} & l_{2}  & p  \\
 k_{1} & k_{2} & q
\end{matrix}\right), 1 \leq l_{1} < l_{2} < p, 1 \leq k_{1}< k_{2} < q.$ Note that $z_{l_{1}q}z_{l_{2}k_{2}}z_{pk_{1}}$ only appears in this $3 \times 3$ matrix,
which yields that the $a_{i}$ associated to this $3 \times 3$ determinant is $0.$ Furthermore, we see that all coefficients $a_{i}$'s that are associated to
$3 \times 3$ determinants are $0.$

Now the conclusion can be proved inductively. Indeed, assume that we
have proved that all coefficients $a_{i}$'s that are associated with
the determinants of order up to $\mu \times \mu, 3 \leq \mu < p$ are
zero. Then we will prove that the coefficients associated with $(\mu
+1) \times (\mu+1)$ determinants are also $0.$ For this we consider
all such determinants which involve $h_{pq},$ i.e.,
$H\left(\begin{matrix}
 l_{1} & ... & l_{\mu} & p  \\
 k_{1} & ... & k_{\mu} & q
\end{matrix}\right)$where $1 \leq l_{1} < ...< l_{\mu} < p, 1 \leq k_{1}<...<k_{\mu}<q.$ We conclude the $a_i$ associated to it is $0$ by noting that
 $z_{l_{1}q}z_{l_{2}k_{\mu}}...z_{l_{\mu}k_{2}}z_{pk_{1}}$ only appears in this $(\mu+1) \times (\mu+1)$ determinant. Then we can show all coefficients that
are associated with other $(\mu+1) \times (\mu+1)$ determinants,
i.e.,
$$H\left( \begin{matrix}
 l_{1} & ...  & l_{\mu}  & l_{\mu +1}  \\
 k_{1} & ...  & k_{\mu}  & k_{\mu+1}
          \end{matrix}
 \right), 1 \leq l_{1}<...<l_{\mu+1} \leq p, 1 \leq k_{1}<...<k_{\mu+1} \leq q, (l_{\mu+1},k_{\mu+1}) \neq (p,q).$$
are $0$ by checking a term of the form
$z_{l_{1}k_{1}}...z_{l_{\mu+1}k_{\mu+1}}$ that only appears once in
the Taylor expansion of the left hand side of (\ref{eq58}). This
proves the lemma.\ \ $\endpf$

We thus get a contradiction to the equation (\ref{eq58}). This
establishes Proposition \ref{thm36}. \ \ $\endpf$

\medskip

\begin{remark}\label{rk38}
Let $F$ be as in  Proposition \ref{thm36}. There exist multiindices
${\beta}^1,...,{\beta}^N$ with  $|{\beta}^j| \leq 1+N-pq$ and
$$z^0=\left( \begin{matrix}
 z^0_{11} & ...  & z^0_{1q}  \\
 ... & ... & ...  \\
 z^0_{p1} & ...  & z^0_{pq}
\end{matrix}\right) \not =0$$
such that $z^0$ is near $0$ and
\begin{equation}
\Delta({\beta}^1,...,{\beta}^N):=\left| \begin{matrix}
\frac{\partial^{|\beta^1|} (\psi_{1}(F))}{\partial \tilde{z}^{{\beta}^1}} & ...  & \frac{\partial^{|\beta^1|} (\psi_{N}(F))}{\partial \tilde{z}^{{\beta}^1}} \\
... & ... & ... \\
\frac{\partial^{|\beta^N|} (\psi_{1}(F))}{\partial
\tilde{z}^{{\beta}^N}} & ... & \frac{\partial^{|\beta^N|}
(\psi_{N}(F))}{\partial \tilde{z}^{{\beta}^N}}
\end{matrix}\right|(z^0) \neq 0.
\end{equation}
Perturbing $z^0$ if necessary, we can thus  assume that $z^0_{pq}
\neq 0.$ Moreover, we can replace one of the
$\beta^1,...,\beta^N$ by $\beta=(0,...,0)$, because
$(\psi_1(F),...,\psi_N(F))$ are not identically zero (See also the proof of Theorem \ref{thm24}).  Without lost
of generality, we can assume that $\beta^1=(0,...,0).$
\end{remark}

The defining function of the Segre family  now is
\begin{equation}
\rho(z,\xi)= 1+ \sum_{k=1}^p\left( \sum_{1 \leq i_{1}<i_{2}<...< i_{k}\leq p,\\ 1 \leq j_{1}< j_{2}<...< j_{k} \leq q}
 Z(\begin{matrix}
 i_{1} & ... & i_{k} \\
 j_{1} & ... & j_{k}
   \end{matrix}
)
 \Xi(\begin{matrix}
 i_{1} & ... & i_{k} \\
 j_{1} & ... & j_{k}
   \end{matrix}
)\right)
\end{equation}
It is a complex manifold for any fixed $\xi$ close enough to the point
$$\xi^0=\left( \begin{matrix}
 0 & ...  & 0 & 0  \\
 0 & ...  & 0 & 0  \\
 0 & ... & 0  & \xi_{pq}^0
         \end{matrix}
 \right) \in \mathbb{C}^{pq}, \quad \xi^0_{pq}=-\frac{1}{z^0_{pq}}.$$
Write for each $1 \leq i \leq p, 1 \leq j \leq q, (i,j) \neq (p,q),$
\begin{equation}
\mathcal{L}_{ij}=\frac{\partial}{\partial z_{ij}}- \frac{\frac{\partial \rho}{\partial z_{ij}}(z,\xi)}{\frac{\partial \rho}{\partial z_{pq}}(z,\xi)}\frac{\partial}{\partial z_{pq}},
\end{equation}
which is a well-defined holomorphic tangent vector field along
$\mathcal{M}$ near $(z^0, \xi^0).$ Here we note that $\frac{\partial
\rho}{\partial z_{pq}}(z,\xi)$ is nonzero near $(z^0, \xi^0).$ For
any $(pq-1)$-multiindex $\beta=(\beta_{11},...,\beta_{p(q-1)})$, we
write
$$\mathcal{L}^{\beta}=\mathcal{L}_{11}^{\beta_{11}}...\mathcal{L}_{p(q-1)}^{\beta_{p(q-1)}}.$$
Now we define for any $N$ collection of $(pq-1)-$multiindices $\{\beta^1,...,\beta^N\},$
\begin{equation}
\Lambda(\beta^1,...,\beta^N)(z,\xi):=\left|\begin{matrix}
 \mathcal{L}^{\beta^1}(\psi_{1}(F)) & ...  & \mathcal{L}^{\beta^1}(\psi_{N}(F)) \\
... & ...  & ...  \\
\mathcal{L}^{\beta^N}(\psi_{1}(F))  & ...  & \mathcal{L}^{\beta^N}(\psi_{N}(F))
\end{matrix}\right|(z,\xi).
\end{equation}

\begin{theorem}\label{lemma39}
There exists multiindices $\{{\beta}^1,...,{\beta}^N\},$ such that
\begin{equation}\label{beta non vanishing}
\Lambda({\beta}^1,...,{\beta}^N)(z,\xi) \neq 0,
\end{equation}
for $(z,\xi)$ in a small neighborhood of $(z^0, \xi^0).$ Moreover, we can require $\beta^1=(0,...,0).$
\end{theorem}
{\it Proof of Theorem \ref{lemma39}}:
First we observe that $\mathcal{L}_{ij}$ evaluating at $(z^0,\xi^0)$
is just $\frac{\partial}{\partial z_{ij}}.$ More generally, for any
$(pq-1)-$multiindex $\beta,$ by an easy computation,
$\mathcal{L}^{\beta}$ evaluating at $(z^0, \xi^0)$ coincides with
$\frac{\partial}{\partial \tilde{z}^{\beta}}.$ Therefore, we can
just choose the same $\beta^1,...,\beta^N$ as in Remark \ref{rk38}.
$\endpf$

\subsection{Spaces of type IV}
In this subsection, we consider the hyperquadric case $M=Q^n$. This
case is more subtle because the tangent vector  fields of its Segre
family are more complicated.
Recall that $Q^n$ is defined by
$$\left \{[z_{0},...,z_{n+1}] \in \mathbb{CP}^{n+1}: \sum_{i=1}^n z_{i}^2-2z_{0}z_{n+1}=0 \right\},$$
where $[z_{0},...,z_{n+1}]$ is the homogeneous coordinates of
$\mathbb{CP}^{n+1}.$ The previously described  minimal  embedding
$\mathbb{C}^n ( {\mathcal A})\rightarrow  Q^n$ is given by
$$z:=(z_{1},...,z_{n}) \mapsto [1,\psi_1(z),...,\psi_{n+1}(z)]= [1,z_{1},...,z_{n},  \frac{1}{2}\sum_{i=1}^n z_{i}^2].$$
The defining function of the Segre family over ${\mathcal A}\times
{\mathcal A}$ is $\rho(z,\xi)=1+r_{z}\cdot r_{\xi},\text{ where }$

\begin{equation}
r_{z}=(z_{1},...,z_{n},\frac{1}{2}\sum_{i=1}^n
z_{i}^2),r_{\xi}=(\xi_{1},...,\xi_{n},\frac{1}{2}\sum_{i=1}^n
\xi_{i}^2).
\end{equation}
Let $F$ be a local biholomorphic map at $0$ with $F(0)=0$.  We write
\begin{equation}\label{eq90}
F=(f_{1},...,f_{n}),\ \ r_{F}=(f_{1},...,f_{n}, \frac{1}{2}
\sum_{i=1}^n f_{i}^2).
\end{equation}
Notice that $$r_{z}=(\psi_{1}(z),...,\psi_{n+1}(z)),
r_{F}=(\psi_{1}(F),...,\psi_{n+1}(F)).$$ We will need the following
lemma:
\begin{lemma}\label{lemma61}
For each fixed $\mu_{1},...,\mu_{n-1}$ with $(\sum_{i=1}^{n-1}
\mu_{i}^2)+1=0$ and each fixed $(z_{1},...,z_{n})$ with
$(\sum_{i=1}^{n-1}\mu_{i}z_{i})+z_{n} \neq 0,$ we can find
$(\xi_{1},...,\xi_{n})$ such that
\begin{equation}\label{eqn91}
1+z_{1}\xi_{1}+...+z_{n}\xi_{n}=0;~~ \sum_{i=1}^{n} (\xi_{i})^2
=0,~~\xi_{j}=\mu_{j}\xi_{n}, 1 \leq j \leq n-1,~~\xi_{n} \neq 0.
\end{equation}
\end{lemma}
{\it Proof of Lemma \ref{lemma61}}: We just need to set
$$\xi_{n}=\frac{-1}{(\sum_{i=1}^{n-1} \mu_{i}z_{i})+z_{n}},~~\xi_{j}=\mu_{j}\xi_{n}, 1 \leq j \leq n-1.$$
Then it is easy to verify that (\ref{eqn91}) is satisfied. $\endpf$

Recall that in the type I case,   the   vector fields
$\frac{\partial}{\partial \widetilde{z}^{\alpha}}$ in
$\mathbb{C}^{pq}$ are  tangent vector fields of the particular
hyperplane $\{z_{pq}=0\}.$ We can  formulate the result in  $\S 3$
in a more abstract way and extend it to a more general setting. For
instance, it can be generalized to the  complex hyperplane case. We
briefly discuss this in more details as  follows:

First fix $\mu_{1},...,\mu_{n-1}$ with $(\sum_{i=1}^{n-1}
\mu_{i}^2)+1=0.$ Take the complex hyperplane $\mathbb{H}:
z_{n}+\sum_{i=1}^{n-1}\mu_{i}z_{i}=0$ in $(z_{1},...,z_{n}) \in
\mathbb{C}^{n}.$ Write
$$L_{1}=\frac{\partial}{\partial z_{1}}-\mu_{1}\frac{\partial}{\partial z_{n}},...,
L_{n-1}=\frac{\partial}{\partial
z_{n-1}}-\mu_{n-1}\frac{\partial}{\partial z_{n}}.$$

Then $\{L_{i}\}_{i=1}^{n-1}$ forms a basis of the tangent vector
fields of $\mathbb{H}.$ For any multiindex
$\alpha=(\alpha_{1},..,\alpha_{n-1}),$ we write
$L^{\alpha}=L_{1}^{\alpha_{1}}... L_{n-1}^{\alpha_{n-1}}.$ We define
$L-$rank and  $L-$nondegeneracy as in Definition \ref{df21} by using
$r_{F}$ in (\ref{eq90}) and by using  $L^{\alpha}$ instead of
$\widetilde{z}^{\alpha}$ with $m=n.$ We write the $k$th $L$-rank
defined in
 this setting as $\mathrm{rank}_{k}(r_{F},L)$  We now need to prove
 the following
\begin{proposition}\label{thm62}
$\mathrm{rank}_{2}(r_{F},L)=n+1.$
\end{proposition}
{\it Proof of Proposition \ref{thm62}}: Suppose not. By applying the
same argument as in Section 3 and a linear change of coordinates, we
can first obtain a modified version of Theorem \ref{thm211}:
\begin{lemma}
There exist $n+1$ holomorphic functions $g_{1}(w),...,g_{n+1}(w)$
which are defined near $0$ on the $w-$plane with $\{ g_{1}(0),...,g_{n+1}(0)\}$ not
all zero such that the following holds for all $z \in U.$
\begin{equation}
\sum_{i=1}^{n+1}
g_{i}(z_{n}+\mu_{1}z_{1}+...+\mu_{n-1}z_{n-1})\psi_{i}(F(z)) \equiv
0.
\end{equation}
\end{lemma}

Then one shows with a similar argument  as in Section 3, by the fact
that $F$ has full rank at $0,$ that $g_{1}(0)=0,...,g_{n}(0)=0.$
Hence we obtain,
\begin{lemma}
There exists a non-zero constant $c \in \mathbb{C}$  such that
\begin{equation}\label{eqn93}
c\psi_{n+1}(F(z))=\frac{c}{2}\sum_{i=1}^n f_{i}^2 (z) \equiv 0,
\end{equation}
for all $z \in U$ when restricted on
$z_{n}+\sum_{i=1}^{n-1}\mu_{i}z_{i}=0.$
\end{lemma}

We then just need to show that (\ref{eqn93}) cannot hold by applying
the following lemma and a linear change of coordinates.

\begin{lemma}\label{lemma681}
Let $H=(h_{1},...,h_{n})$ be a vector-valued holomorphic function in
a neighborhood $U$ of $0$ in $\tilde{z}= (z_{1},...,z_{n-1})\in
\mathbb{C}^{n-1}$ with $H(0)=0.$ Assume that $H$ has full rank at
$0.$ Assume that $a$ is a complex number such that,
\begin{equation}\label{eqn94}
a \sum_{i=1}^n h_{i}^2 (\tilde{z}) \equiv 0,
\end{equation}
Then $a=0.$
\end{lemma}
{\it Proof of Lemma \ref{lemma681}}: Seeking a contradiction,
suppose not.  Notice that $H$ has full rank at $0.$ We assume,
without loss of generality, that $(h_{1},...,h_{n-1})$ gives a local
biholomorphic map near $0$ from $\mathbb{C}^{n-1}$ to
$\mathbb{C}^{n-1}.$ By a local biholomorphic change of coordinates,
we assume $(h_{1},...,h_{n-1})=(z_{1},...,z_{n-1}),$ and still write
the last component as $h_{n}.$ Then equation (\ref{eqn94}) is
reduced to
$$a(z_{1}^2+...+z_{n-1}^2+ h_{n}^2)=0.$$
To cancel the $z_{1}^2, z_{2}^2$ terms, it yields that $h_{n}$ has
linear $z_{1}, z_{2}$ terms. But then $h_{n}^2$ would produce  a
$z_{1}z_{2}$ term, which cannot be canceled out. This is a
contradiction. \ \ $\endpf$

This also establishes Proposition \ref{thm62}. \ \ $\endpf$

\medskip
\begin{remark}
By Proposition \ref{thm62}, there exist multiindices
$\tilde{\beta}^1,...,\tilde{\beta}^{n+1}$ with  $|\tilde{\beta}^j|
\leq 2$ and
$$z^0=(z_{1}^0,...,z_{n}^0)~\text{with}~~ \sum_{i=1}^{n-1}\mu_{i}z^0_{i}+z^0_{n} \neq 0$$
such that
\begin{equation}
\left| \begin{matrix}
 L^{\tilde{\beta}^1}(\psi_{1}(F)) & ...  & L^{\tilde{\beta}^1} (\psi_{n+1}(F))\\
... & ... & ... \\
L^{\tilde{\beta}^{n+1}} (\psi_{1}(F)) & ... &
L^{\tilde{\beta}^{n+1}} (\psi_{n+1}(F))
\end{matrix}\right|(z^0) \neq 0.
\end{equation}
\end{remark}
We then choose $\xi^0=(\xi_{1}^0,...,\xi_{n}^0)$ as in Lemma
\ref{lemma61}. That is
$$1+z_{1}^0\xi_{1}^0+...+z_{n}^0\xi_{n}^0=0;~~
\sum_{i=1}^{n} (\xi_{i}^0)^2 =0,~~\xi_{j}^0=\mu_{j}\xi_{n}^0, 1 \leq
j \leq n-1,~~\xi_{n}^0 \neq 0.$$ It is easy to see that $(z^0,
\xi^0) \in \mathcal{M}.$
We now define
\begin{equation}
\mathcal{L}_{i}=\frac{\partial}{\partial z_{i}}-
\frac{\frac{\partial \rho}{\partial z_{i}} (z,\xi)}{\frac{\partial
\rho}{\partial z_{n}}(z,\xi)}\frac{\partial}{\partial z_{n}}, 1 \leq
i \leq n-1
\end{equation}
for $(z,\xi) \in \mathcal{M}$ near $(z^0, \xi^0).$ They are
well-defined holomorphic tangent vector fields along $\mathcal{M}.$
Moreover, $\frac{\partial \rho}{\partial z_{n}}(z,\xi)$ is nonzero
near $(z^0, \xi^0).$

We define for any multiindex $\alpha=(\alpha_{1},..,\alpha_{n-1}),$
$\mathcal{L}^{\alpha}=
\mathcal{L}_{1}^{\alpha_{1}}...\mathcal{L}_{n-1}^{\alpha_{n-1}}.$
Then for any $(n+1)$ collection of $(n-1)-$multiindices, set
$\{\beta^1,...,\beta^N\},$
\begin{equation}
\Lambda(\beta^1,...,\beta^{n+1})(z,\xi):=\left|\begin{matrix}
 \mathcal{L}^{\beta^1}(\psi_{1}(F)) & ...  & \mathcal{L}^{\beta^1}(\psi_{n+1}(F)) \\
... & ...  & ...  \\
\mathcal{L}^{\beta^{n+1}}(\psi_{1}(F))  & ...  &
\mathcal{L}^{\beta^{n+1}}(\psi_{n+1}(F))
\end{matrix}\right|(z,\xi).
\end{equation}

\bigskip

By the fact that $\sum_{i=1}^n (\xi_{i}^0)^2=0,$ one can check that,
for any multiindex
 $\alpha=(\alpha_{1},..,\alpha_{n}),$ $\mathcal{L}^{\alpha}=L^{\alpha}$ when
 evaluated
at $(z^0, \xi^0).$ Then we get the following:

\begin{theorem}
There exist multiindices $\{{\beta}^1,...,{\beta}^N\}$ such that
$$\Lambda({\beta}^1,...,{\beta}^N)(z,\xi) \neq 0,$$
for $(z,\xi)$ in a small neighborhood of $(z^0, \xi^0)$, where
$\beta^1=(0,0,...,0).$
\end{theorem}
Proofs for the other types are similar and will be left to
Appendix II.



\section{Transversality and flattening of  Segre families: Proof of Proposition (II)}

In this section, we prove Proposition (II). We still use the notations we have set up so far. We equip the space $M$ with the canonical K\"ahler-Einstein metric $\omega$ as
described before.  We start with the following lemma:

\begin{lemma}\label{deg of isometry} Let $\wh{\sigma}: (M,\o)\rightarrow (M,\o)$
be a holomorphic isometry. In the affine space ${\mathcal A}$,
its components consist of  rational functions with its degree
bounded only by a constant depending on $(M,\o)$.
\end{lemma}
{\it Proof of Lemma \ref{deg of isometry}}: Notice that $M$ has been isometrically
 embedded into $\mathbb{CP}^N$ through the canonical map defined before.
Hence $\wh{\sigma}$  is the  restriction of a unitary
transformation.
 Hence $\wh\sigma$ can be identified with a map of the form:
$$(\tilde\psi_0,\tilde\psi_1,\tilde\psi_2,...,\tilde\psi_N)
=(\sum_{j=0}^Na_{0j}\psi_j,...,\sum_{j=0}^Na_{ij}\psi_j,...,\sum_{j=0}^Na_{Nj}\psi_j),$$
where $\psi_0=1$ and $(a_{ij})$ is a unitary matrix.  Write
$$\Psi(z): z(\in {\mathcal A})\mapsto
[1,\kappa_1z_1,\cdots,\kappa_iz_i,\cdots,\kappa_nz_n,o(z^2)]\in
\mathbb{CP}^N$$ for the embedding, where $\kappa_i=1$ or $\sqrt 2$.
$\wh{\sigma}$ induces a birational  self-action $\sigma$ of
$\mathcal A$ such that $\Psi(\sigma(z))=\wh{\sigma}(\Psi(z))$. Then,
from the special form of $\Psi$,
$\sigma(z)=\left(\frac{\tilde\psi_1}{\kappa_1\tilde\psi_0},\frac{\tilde\psi_2}{\kappa_2\tilde\psi_0},...,
\frac{\tilde\psi_n}{\kappa_n\tilde\psi_0}\right).$ Apparently $\tilde\psi_0 \not  \equiv
0$.
$\endpf$

\medskip

\begin{theorem}\label{flatten}Suppose $\xi^0\in\mathbb{C}^n\sm\{0\}$.
Then for a generic smooth point $z^0$ on the Segre variety
$Q_{\xi^0}$ and a small neighborhood $U\subset {\mathbb
C}^n$ of $z^0,$
 there is a  point $z^1 \in U\cap Q_{\xi^0}$, such that $Q_{z^0}$ and $Q_{z^1}$ are both smooth at $\xi^0$ and 
 intersect transversally there.
  Moreover,
  there is a   biholomorphic parametrization
   $\mathcal{G}(\tilde\xi_1,\tilde\xi_2,...,\tilde\xi_n)
   =(\xi_1,\xi_2,...,\xi_n),$ with
$   (\tilde\xi_1,\tilde\xi_2,...,
   \tilde\xi_n)\in U_1\times U_2\times...\times U_n\subset\mathbb C^n$. Here when $1 \leq j \leq 2, $  $U_j$ is a small neighborhood of $1\in {\mathbb C}$. When $3 \leq j \leq n,$ $U_j$ is a small neighborhood of $0 \in \mathbb{C}$  with
   ${\mathcal G}(1,1,0,\cdots,0)=\xi^0$,
   such that $\mathcal{G}(\{\tilde\xi_1=1\}\times U_2\times...\times U_n)\subset Q_{z^0},
   \mathcal{G}(U_1\times\{\tilde\xi_2=1\}\times U_3\times...\times U_n)\subset Q_{z^1},$
    and $\mathcal{G}(\{\tilde\xi_1=t\}\times U_2\times...\times U_n),
   \mathcal{G}(U_1\times\{\tilde\xi_2=s\}\times U_3\times...\times U_n),s\in U_1,t\in U_2$
    are open pieces of  Segre varieties. Also,
     $\mathcal{G}$  consists of algebraic functions with total degree
     bounded by a constant depending only on $(M,\o)$.
\end{theorem}
We first claim that, due to the invariance of the Segre family, we
need only to prove the theorem
for a special point $0\neq\xi^0\in\mathbb C^n\subset M$.
Indeed,  by the invariance property mentioned in $\S 2$, for an
isometry $\sigma$, $(\sigma,\Ol{\sigma})$ preserves the Segre family
${\mathcal M}\subset M\times M.$ Here for $p\in \mathbb{CP}^N$,
$\Ol{\sigma}(p):=\Ol{\sigma(\Ol p)}$ as  before. Here, we mention
that in the statement of the theorem, we assume that $z^0$ is a
generic smooth point because under this transformation, some smooth
points on $Q_{\xi^0}$ may be mapped into the hyperplance of $M$ at
infinity, which can not  be chosen as our $z^0$.

We now proceed to the proof of Theorem $\ref{flatten}$ by choosing a
good point $\xi^0$.  We only carry out the proof for the case of
hyperquadrics and Grassmannian spaces here. The proof for the
remaining cases is  similar and will be  included in Appendix III.

\medskip
{\it Proof of Theorem \ref{flatten}}: \textbf{Case 1.
Hyperquadrics}:
 Suppose $M$ is the hyperquadric. Then the defining equation for
the Segre family is
$$\rho(z,\xi)=1+\sum_{i=1}^nz_i\xi_i+\frac{1}{4}(\sum_{i=1}^nz_i^2)(\sum_{i=1}^n\xi_i^2)=0.$$
 We choose $\xi^0=(1,0,0,...,0).$  Hence
  $Q_{\xi^0}=\{z:\ \rho(z,\xi^0)=1+z_1+\frac{1}{4}(\sum_{i=1}^nz_i^2)=0\}.$
  We compute the gradient of $\rho(z,\xi^0)$ as follows:  $\nabla \rho(z,\xi^0)=(1+\frac{1}{2}z_1,\frac{1}{2}z_2,...,
  \frac{1}{2}z_n).$ Notice that  $Q_{\xi^0}$ is smooth  except at
  $(-2,0,...,0),$ namely, we have  $\nabla \rho(z,\xi^0)\not
   =0$ away from $(-2,0,\cdots,0)$.  For a smooth point $z^0(\not = (-2,0,\cdots,0))$  of $ Q_{\xi^0}$ , we  choose a neighborhood $ U$ of $z^0$ in ${\mathbb C}^n$ such that
   $ U\cap Q_{\xi^0}$ is a smooth piece of
   $Q_{\xi^0}$.
 Pick  also $z^1 (\not = z^0)\in U\cap Q_{\xi_0}$ and compute the
gradient of the defining function of $Q_{z^0}$ and $Q_{z^1}$ at
$\xi^0=(1,0,...,0),$ respectively.  Recall
$$Q_{z^s}=\{\xi|\rho(z^s,\xi)=1+\sum_{i=1}^nz_i^s \xi_i+\frac{1}{4}(\sum_{i=1}^n(z_i^s)^2)
(\sum_{i=1}^n\xi_i^2)=0\},~\text{for}~s=0,1.$$

\begin{displaymath}\left(\begin{matrix}\nabla \rho(z^0,\xi)|_{\xi^0=(1,0,...,0)}\\
\nabla \rho(z^1,\xi)|_{\xi^0=(1,0,...,0)}

\end{matrix}\right)
=\left(\begin{matrix}z_1^0+\frac{1}{2}\sum_{i=1}^n(z_i^0)^2&z_2^0&z_3^0&...&z_n^0\\
z_1^1+\frac{1}{2}\sum_{i=1}^n(z_i^1)^2&z_2^1&z_3^1&...&z_n^1
\end{matrix}\right)
=\left(\begin{matrix}-2-z_1^0&z_2^0&z_3^0&...&z_n^0\\
-2-z_1^1&z_2^1&z_3^1&...&z_n^1
\end{matrix}\right)
\end{displaymath}
The second equality is simplified by making use of the fact that
$z^0,z^1\in Q_{\xi^0=(1,0,...,0)},$ which implies that
$0=1+z_1^0+\frac{1}{4}\sum_{i=1}^n(z_i^0)^2=1+z_1^1+\frac{1}{4}\sum_{i=1}^n(z_i^1)^2.$
Hence,
\begin{displaymath}\rm{rank}\left(\begin{matrix}\nabla \rho(z^0,\xi)|_{\xi^0=(1,0,...,0)}\\
\nabla \rho(z^1,\xi)|_{\xi^0=(1,0,...,0)}

\end{matrix}\right)
=\rm{rank}\left(\begin{matrix}-2-z_1^0&z_2^0&...&z_n^0\\
-2-z_1^1&z_2^1&...&z_n^1

\end{matrix}\right)
=\rm{rank}\left(\begin{matrix}-2-z_1^0&z_2^0&...&z_n^0\\
-\Delta z_1^1&\Delta z_2^1&...&\Delta z_n^1
\end{matrix}\right)
\end{displaymath}

\begin{displaymath}
=\rm{rank}\left(\begin{matrix}2+z_1^0&z_2^0&...&z_n^0\\
\Delta z_1^1&\Delta z_2^1&...&\Delta z_n^1
\end{matrix}\right)
=\rm{rank}\left(\begin{matrix}&&\nabla \rho(z,\xi^0)|_{z^0}\\
\Delta z_1^1&\Delta z_2^1&...&\Delta z_n^1
\end{matrix}\right),
\end{displaymath}
where $ \Delta z_i^1:=z_i^1-z_i^0.$ Notice that $z^0$ is a smooth
point on $Q_{\xi_0}$. Hence $\nabla\rho(z,\xi^0)$ is transversal to
the  tangent space of $Q_{\xi^0}$ at $z^0.$  If we choose $z^1\in
Q_{\xi^0}$ close enough to $z^0$, which ensures $(\Delta
z_1^1,...,\Delta z_n^1)$ close enough to tangent space of
$Q_{\xi^0}$ at $z^0$, we  then get
\begin{displaymath}
\rm{rank}\left(\begin{matrix}\nabla \rho(z^0,\xi)|_{\xi^0=(1,0,...,0)}\\
\nabla \rho(z^1,\xi)|_{\xi^0=(1,0,...,0)}
\end{matrix}\right)
=\rm{rank}\left(\begin{matrix}&&\nabla \rho(z,\xi^0)|_{z^0}\\
\Delta z_1^1&\Delta z_2^1&...&\Delta z_n^1.
\end{matrix}\right)=2.
\end{displaymath} \\
We  assume, without loss of generality, that   $\frac{\partial
(\rho(z^0,\xi),\rho(z^1,\xi))}{\partial (\xi_1,\xi_2)}\neq 0$ at
$\xi^0$.
Now we introduce new variables $\tilde\xi_1,...,\tilde\xi_n$ and
consider the following system of equations:
\begin{displaymath}
\begin{cases}P_1:1+\sum_{i=1}^n(\tilde\xi_1z_i^0)\xi_i+\frac{1}{4}(\sum_{i=1}^n(\tilde\xi_1)^2(z_i^0)^2)
(\sum_{i=1}^n\xi_i^2)&=0\\
P_2:1+\sum_{i=1}^n(\tilde\xi_2z_i^1)\xi_i+\frac{1}{4}(\sum_{i=1}^n(\tilde\xi_2)^2(z_i^1)^2)(\sum_{i=1}^n\xi_i^2)&=0\\
P_3:\quad\quad\quad\quad\quad\quad\quad\quad\quad\quad\quad\quad\quad\quad\quad\quad\quad\,\,\tilde\xi_3-\xi_3&=0\\
...&...\\
P_n:\quad\quad\quad\quad\quad\quad\quad\quad\quad\quad\quad\quad\quad\quad\quad\quad\quad\,\,\tilde\xi_n-\xi_n&=0
\end{cases}
\end{displaymath}
Then we have
$\frac{\partial(P_1,...,P_n)}{\partial(\xi_1,...,\xi_n)}|_{A}\neq 0$
and
$\frac{\partial(P_1,...,P_n)}{\partial(\tilde\xi_1,...,\tilde\xi_n)}|_{A}\neq
0$ where
$$A=(\tilde\xi_1,...,\tilde\xi_n;\xi_1,...,\xi_n)=(1,1,0,...,0;1,0,...,0).$$
 By  Lemma \ref{james-002},   we get the needed algebraic
 flattening with total degree bounded only by $(M,\omega)$. This completes the proof of Theorem $\ref{flatten}$ in the
 hyperquadric case.

\medskip
 {\bf Case 2. Grassmannians}: Pick
$\xi^0=(\xi_{11}^0,\xi_{12}^0,...,\xi_{pq}^0)=(1,0,...,0).$ The
defining function for the Segre family associated with this point is
as follows:

$\rho(z,\xi)=1+z_{11}\xi_{11}+z_{12}\xi_{12}+...+z_{1q}\xi_{1q}+z_{21}\xi_{21}+...+z_{p1}\xi_{p1}
+\sum_{i,j\neq 1}z_{ij}\xi_{ij}+\sum_{i,j\geq
2}(z_{11}z_{ij}-z_{i1}z_{1j})(\xi_{11}\xi_{ij}-\xi_{i1}\xi_{1j})+\sum_{(i,j),
    (k,l)\neq (1,1)}(z_{ij}z_{kl}-z_{il}z_{jk})(\xi_{ij}\xi_{kl}-\xi_{il}\xi_{jk})+\hbox{ higher order terms}.$

Then
$Q_{\xi^0}=\{z|\rho(z,\xi^0)=1+z_{11}=0\},\nabla\rho(z,\xi^0)=(1,0,0,...,0).$
Hence $Q_{\xi_0}$
is smooth. 
For $z\in Q_{\xi^0}$, we have
$z=(-1,z_{12},...,z_{1q},z_{21},...,z_{p1},...,z_{ij},...,z_{pq}).$
Pick $z^0,z^1\in Q_{\xi^0}$. Then

$Q_{z^s}=\{\xi|0=\rho(z^s,\xi)=1+z_{11}^s\xi_{11}+z_{12}^s\xi_{12}+...+z_{1q}^s\xi_{1q}+z_{21}^s\xi_{21}+...
+z_{p1}^s\xi_{p1}+\sum_{i,j\neq 1}z_{ij}^s\xi_{ij}+\sum_{i,j\geq
    2}(z_{11}^sz_{ij}^s-z_{i1}^sz_{1j}^s)(\xi_{11}\xi_{ij}-\xi_{i1}\xi_{1j})+\sum_{(i,j),
    (k,l)\neq (1,1)}(z_{ij}^sz_{kl}^s-z_{il}^sz_{jk}^s)(\xi_{ij}\xi_{kl}-\xi_{il}\xi_{jk})+{\rm high\, order\, terms}\}$,\\ for $s=0,1.$
We then compute their gradients as follows:
\begin{displaymath}\left(\begin{matrix}\nabla \rho(z^0,\xi)|_{\xi^0}\\
\nabla \rho(z^1,\xi)|_{\xi^0}
\end{matrix}
\right) =\left(\begin{matrix}\frac{\partial
\rho(z^0,\xi)}{\partial\xi_{11}}&\frac{\partial
\rho(z^0,\xi)}{\partial\xi_{12}}&...&\frac{\partial \rho(z^0,\xi)}
{\partial\xi_{1q}}&\frac{\partial
\rho(z^0,\xi)}{\partial\xi_{21}}&...&\frac{\partial
\rho(z^0,\xi)}{\partial\xi_{p1}}&...
&\frac{\partial \rho(z^0,\xi)}{\partial\xi_{pq}}\\
\frac{\partial \rho(z^1,\xi)}{\partial\xi_{11}}&\frac{\partial
\rho(z^1,\xi)}{\partial\xi_{12}}&...&\frac{\partial
\rho(z^1,\xi)}{\partial\xi_{1q}}
&\frac{\partial \rho(z^1,\xi)}{\partial\xi_{21}}&...&\frac{\partial \rho(z^1,\xi)}{\partial\xi_{p1}}&...&\frac{\partial \rho(z^1,\xi)}{\partial\xi_{pq}}\\
\end{matrix}\right)\big{|}_{\xi^0}
\end{displaymath}
\begin{displaymath}
=\left(\begin{matrix}-1&z_{12}^0&...&z_{1q}^0&z_{21}^0&...&z_{p1}^0&-z_{i1}^0z_{1j}^0&...\\
-1&z_{12}^1&...&z_{1q}^1&z_{21}^1&...&z_{p1}^1&-z_{i1}^1z_{1j}^1&...\\
\end{matrix}\right).
\end{displaymath}
Thus, we have
\begin{displaymath}
\rm{rank}\left(\begin{matrix}\nabla \rho(z^0,\xi)\big{|}_{\xi^0}\\
\nabla \rho(z^1,\xi)\big{|}_{\xi^0}
\end{matrix}\right)
=
\rm{rank}\left(\begin{matrix}-1&z_{12}^0&...&z_{p1}^0&-z_{i1}^0z_{1j}^0&...\\
0&\Delta z_{12}^1&...&\Delta z_{p1}^1&(-z_{i1}^0\Delta z_{1j}^1-
z_{1j}^0\Delta z_{i1}^1-\Delta z_{i1}^1\Delta z_{1j}^1)&...\\
\end{matrix}\right),
\end{displaymath}
where $\Delta z_{ij}^1=z_{ij}^1-z_{ij}^0.$  Hence, if we choose
$z^1$ such that $z_{12}^1\neq z_{12}^0$, Then the  rank equals to
$2$. Hence $Q_{z^0}$ and $Q_{z^1}$ are smooth and intersect transversally at
$\xi^0$.

Without loss of generality, assume $\frac{\partial
(\rho(z^0,\xi),\rho(z^1,\xi))}{\partial (\xi_{11},\xi_{12})}\neq 0$ at $\xi^0.$ Now we introduce new variables $\tilde\xi_{11},...,\tilde\xi_{pq}$
and set up the system:
\begin{displaymath}
\begin{cases}P_{11}:\,\,&\rho(z^0,\tilde\xi_{11}\xi)=0\\
P_{12}:\,\,&\rho(z^1,\tilde\xi_{12}\xi)=0\\
P_{13}:
\,\,&\tilde\xi_{13}-\xi_{13}=0\\
\,...&...\\
P_{pq}:
\,\,& \tilde\xi_{pq}-\xi_{pq}=0
\end{cases}
\end{displaymath}
Then
$\frac{\partial(P_{11},...,P_{pq})}{\partial(\xi_{11},...,\xi_{pq})}|_{A},\
 \frac{\partial(P_{11},...,P_{pq})}{\partial(\tilde\xi_{11},...,\tilde\xi_{pq})}|_{A}\neq 0$, where
$A=(\tilde\xi_{11},...,\tilde\xi_{pq},\xi_{11},...,\xi_{pq})=(1,1,0,...,0,1,0,...,0)$.
By Lemma \ref{james-002}, we get the needed algebraic flattening.

The proof is similar in the  other cases. We include a detailed argument for the remaining cases in Appendix III. $\endpf$

\section{Irreducibility of Segre varieties: Proof of Proposition (III)}
In this section we will establish Proposition (III). We prove  results
on the irreducibility of the potential function $\rho$, Segre
varieties and the Segre family. We still adapt the previously used
notation and assume that $M$ is an irreducible Hermitian symmetric
space of compact type of dimension $n$, which has been
minimally embedded into a projective space as described before.

\begin{lemma}\label{singleirre} Each Segre variety  is an irreducible algebraic subvariety.
\end{lemma}

{\it Proof of lemma \ref{singleirre}}:  For a minimally embedded
Hermitian symmetric space, since all Segre varieties are unitarily
equivalent, it suffices to prove the lemma for a single Segre
variety. Without lost of generality, we take $z=(0,...,0)\in\mathcal
A\subset M$. Therefore, the corresponding Segre variety $Q^*_{z}$ is
the hyperplane section $M\sm {\mathcal A}$, which is of pure
dimension. From the classical algebraic geometry \cite{GH}, when $M$
is an irreducible Hermitian symmetric space of compact type, the
hyperplane section at infinity in the minimal canonical embedding
case is a union of Schubert cells. Moreover as shown in \cite{CMP},
the top dimensional piece is equivalent to $\mathbb C^{n-1}$ and the
others are of codimension at least two. Hence, the smooth points of
$Q_{z}$ are connected and thus $Q_z$ is irreducible. $\endpf$

\medskip

As a corollary of this lemma, we conclude that for each $z\in
{\mathbb C}^n$, the defining function $\rho(z,\cdot)$ of $Q_z$ has
to be a power of  one irreducible factor. However, as in the proof
of Theorem \ref{flatten}, for some $a(\not =0)\in {\mathbb C}^n$,
$d_{\xi}\rho(a,\xi)$ is not identically zero along $Q_a$. Next, we
use this property and the symmetric property of $M$ to prove the
following:

\begin{proposition}\label{irreducibility} For any $b\in  {\mathcal A}$ with $ b\not =(0,...,0)$, $\rho(b,\xi)$ ($\rho(z,b)$, respectively) is irreducible as a
polynomial of $\xi$ (as a polynomial in $z$, respectively).

\end{proposition}

{\it Proof of proposition \ref{irreducibility}}: Since
$\rho(z,\xi)=\rho(\xi,z)$, we need just to verify the first
statement.
Let $a$ be as above. For $b\in\mathcal A$, there is
$\widehat\sigma\in$ ${\rm Isom}(M,\omega)\cap SU(N+1,\mathbb C)$
such that $\widehat\sigma(a)=b$. (Notice that $\wh\sigma$ is
represented by a unitary action.) By Lemma \ref{deg of isometry},
let $\sigma=(\frac{l_1}{\kappa_1l_0},...,\frac{l_n}{\kappa_nl_0})$
be the representation of $\wh{\sigma}$  in $\mathcal A$ with $l_j's$
polynomials in $z$. Write $\Psi=[1,r_z]$ for the embedding of
$\mathcal A$ in $\mathbb P^N$. Then from the definition of
$\rho(z,\ov z),$ we have
$$\rho(z,\ov z)=||\Psi(z)||^2=\Psi\cdot\ov \Psi^t=(\widehat\sigma \Psi)\cdot\ov {(\wh\sigma \Psi)}^t.$$ 

\begin{lemma}\label{lemmasigl}
$(\widehat\sigma \Psi)\cdot\ov {(\wh\sigma \Psi)}^t=|l_0(\Psi)|^2\cdot||\Psi(\sigma(z))||^2
=|l_0(\Psi)|^2\cdot\rho(\sigma(z),\ov{\sigma(z)}).$ 
\end{lemma}

\begin{proof}
Writing $\Psi(z)=[1, r_z]=[1, \psi_1(z), \cdots, \psi_N(z)].$ Then the identity $\Psi(\sigma(z))=\wh{\sigma}(\Psi(z))$ obtained in the proof
of Lemma \ref{deg of isometry} yields that,
$$(\psi_1(\sigma(z)), \cdots, \psi_N(\sigma(z)))=\left(\frac{\tilde\psi_1(\Psi(z))}{\tilde\psi_0(\Psi(z))}, \cdots, \frac{\tilde\psi_N(\Psi(z))}{\tilde\psi_0(\Psi(z))}  \right).$$
Here $\tilde\psi_j=l_j$ for $0 \leq j \leq n$ and $\widehat\sigma(z)=[\tilde\phi_0, \cdots, \tilde\phi_N]$ as in the proof of Lemma \ref{deg of isometry}. Then
$$(\hat\sigma \Psi) \cdot \ov {(\wh\sigma \Psi)}^t=\sum_{j=0}^N |\tilde\psi_j(\Psi(z))|^2=\left(1+\sum_{j=1}^N |\psi_j (\sigma(z))|^2 \right) |\tilde\psi_0(\Psi(z))|^2=|l_0(\Psi)|^2\cdot||\Psi(\sigma(z))||^2.$$ This establishes the lemma.
\end{proof}

The Lemma \ref{lemmasigl} yields $\rho(z,\ov z)=|l_0(\Psi)|^2\cdot\rho(\sigma(z),\ov{\sigma(z)}).$
Complexifying
the identity  and substituting $z$ by $a$, we have:
\begin{equation}\label{trans}
l_0(\Psi)(a)\cdot\ov{l_0(\Psi)}(\xi)\cdot\rho(b,\overline{\sigma}(\xi))=\rho(a,\xi),
\end{equation}
where $l_0(\Psi)(a)\not =0,$ $l_0(\Psi)(\xi)$, $\rho(a,\xi)$ are
polynomials in $\xi$ and $\sigma(\xi)$ is a rational map in $\xi$.
Now supposing $\rho(b,\xi)=f^l(\xi),l\geq 2,$ we have
$\rho(b,{\ov\sigma}(\xi))=(f({\ov
\sigma}(\xi)))^l=(\frac{f_1(\xi)}{f_2(\xi)})^l$, where $f_1$ and
$f_2$ are coprime polynomials. Since $a,b\not=(0,...,0)$, $f_1$ is a
non-constant polynomial.  Therefore in  \eqref{trans}, even after
cancellation,
we still have a factor $f_1^l(\xi)$. However as shown in $\S 6$, the
right hand side of the identity \eqref{trans} must be an irreducible
polynomial, which is a contradiction.\ $\endpf$






\begin{proposition}\label{globalirre} $\rho(z,\xi)$ is an irreducible polynomial over ${\mathbb C}^n\times {\mathbb C}^n$. Thus, the Segre family $\mathcal M$ restricted to $\mathbb C^n\times\mathbb C^n=\mathcal A\times\mathcal A\subset M\times M$ is an irreducible subvariety of dimension $2n-1.$
\end{proposition}

We also have the following slightly strong version of the above
proposition, which was used for applying a monodromy argument:

\begin{proposition}\label{connectivity}Suppose $U$ is an connected open set in $\mathbb{C}^n\sm \{0\}$.
 Then the  Segre family $\mathcal{M}$ restricted to $U\times\mathbb{C}^n$ or restricted to
 $\mathbb{C}^n\times U$ is an irreducible analytic variety.
\end{proposition}

{\it Proof of Proposition \ref{connectivity}}: We need only to prove
the first statement.
Recall the notations we defined before:
 $\mathcal{M}_{\rm{SING}}=\{(z,\xi):
  \frac{\partial \rho}{\partial \xi_{j}}=0,\forall j \}\cup\{(z,\xi): \frac{\partial \rho}
  {\partial z_{j}}=0,\forall j \},$ and $\mathcal{M}_{\rm{REG}}=\mathcal{M}\backslash\mathcal{M}_{\rm{SING}}.$
Since $\rho(z,\xi)$ is an irreducible polynomial and $\frac{\partial
\rho}{\partial \xi_{j}}, \frac{\partial \rho}{\partial
z_{j}},j=1,...,n$ are polynomials with lower degrees,
$\frac{\partial \rho}{\partial \xi_{j}},\frac{\partial
\rho}{\partial z_{j}},j=1,...,n$ are not identically zero on
$\mathcal{M}=\{\rho(z,\xi)=0\}.$ Each of $\frac{\partial
\rho}{\partial \xi_{j}}, \frac{\partial \rho}{\partial z_{j}}$
defines a proper subvariety inside $\mathcal{M}.$
By Proposition \ref{irreducibility},
for each $\tilde z (\not = 0)\in\mathbb{C}^n,$ there is a certain
point $\tilde \xi$ on $Q_{\tilde z}$ such that a  partial derivative
of $\rho(\tilde{z},\xi)$ in $\xi$ at $(\tilde z,\tilde\xi)$ does not
vanish. Hence ${\mathcal M}_{\rm SING}$ does not contain any Segre
variety. Also the standard projection from $\mathcal{M}_{\rm REG}$
into the
$z$-space 
 is a submersion.
 Since $Q_z$ is irreducible for
$z\in \mathbb{C}^n\backslash(0,...,0)$,
$Q_z\cap\mathcal{M}_{\rm{REG}}$ is connected.
To prove the theorem,
we just need to  show that
$\mathcal{M}_{\rm{REG}}|_{U\times\mathbb{C}^n}$
is connected. Write the above projection map to the $z$-space  as
$\Phi:\mathcal{M}_{\rm{REG}}|_{U \times \mathbb{C}^n}\rightarrow U.$
Since it is a submersion, it is an open mapping. Suppose $z^0$ is a
point in $U.$   As mentioned above, we know that each fiber of
$\Phi$ is connected.  For any $(z^0,\xi^0)\in {\mathcal M}_{\rm
REG}$ in the fiber above $z^0$,
we can choose a connected neighborhood $V$ of $(z^0,\xi^0)$  on
${\mathcal M}_{\rm{REG}}|_{U \times \mathbb{C}^n}$ such that
$\Phi(V)$ is neighborhood of $z_0$. Hence, for any $z\in\Phi(V)$,
any point  in $Q_z\cap {\mathcal M}_{\rm {REG}}$ can be connected by
a smooth curve inside $\mathcal{M}_{\rm{REG}}|_{V \times
\mathbb{C}^n}$  to $(z^0,\xi^0)$. Since $U$ is connected, by a
standard open-closeness argument, we see that
$\mathcal{M}_{\rm{REG}}|_{U \times \mathbb{C}^n}$ is connected.
 $\endpf$

\medskip

\appendix


\section{Appendix I:\
Affine cell coordinate functions for  two  exceptional classes of
the Hermitian symmetric spaces of compact type}

Define the multiplication law of octonions with the standard basis
$\{e_0=1,e_1,\cdots,e_7\}$  by the following table:

\begin{table}[!hbp]
\begin{tabular}{|c|c|c|c|c|c|c|c|}
\hline
\hline
{} & $e_1$ & $e_2$ & $e_4$ & $e_7$ & $e_3$ & $e_6$ & $e_5$ \\
\hline
 $e_1$& $-1$&$e_4$ &$-e_{2}$ &$-e_{3}$&  $e_{7}$&  $-e_{5}$&  $e_{6}$\\
\hline
 $e_2$& $-e_{4}$&$-1$ &$e_{1}$ &$-e_{6}$&  $e_{5}$&  $e_{7}$&  $-e_{3}$\\
\hline
 $e_4$& $e_2$&$-e_1$ &$-1$ &$-e_{5}$&  $-e_{6}$&  $e_{3}$&  $e_{7}$\\
\hline
 $e_7$ & $e_3$&$e_6$ &$e_{5}$ &$-1$&  $-e_{1}$&  $-e_{2}$&  $-e_{4}$\\
\hline
 $e_3$& $-e_7$&$-e_5$ &$e_{6}$ &$e_{1}$&  $-1$&  $-e_{4}$&  $e_{2}$\\
\hline
 $e_6$& $e_5$&$-e_7$ &$-e_{3}$ &$e_{2}$&  $e_{4}$&  $-1$&  $-e_{1}$\\
\hline
 $e_5$& $-e_6$&$e_3$ &$-e_{7}$ &$e_{4}$&  $-e_{2}$&  $e_{1}$&  $-1$\\
\hline
\end{tabular}
\end{table}

{\bf $\bf \clubsuit 1.$ {Case $\bf M_{16}$}}: Define
\begin{displaymath}
\begin{matrix}
x&=&&(x_0,x_1,x_2,x_3,x_4,x_5,x_6,x_7),\\
y&=&&(y_0,y_1,y_2,y_3,y_4,y_5,y_6,y_7).\\
\end{matrix}
\end{displaymath}

Define $A_j(x,y),j=0,\dots,7,$ such that
\begin{displaymath}
x\bar y=\sum_{j=0}^{7}A_j(x,y)e_j, {\,\,\rm where\,\,} x=\sum_{j=0}^{7}x_je_j {\,\,\rm and\,\,} y=\sum_{j=0}^{7}y_je_j.
\end{displaymath}

Define $B_j(x,y),j=0,1$ such that
\begin{displaymath}
x\bar x=B_0(x,y)e_0 {\,\,\rm and\,\,} y\bar y=B_1(x,y)e_0.
\end{displaymath}

Then by computation, we have the following formulas:
\begin{displaymath}
\begin{matrix}
A_0=&A_0(x,y)=&&y_0x_0+y_1x_1+y_2x_2+y_3x_3+y_4x_4+y_5x_5+y_6x_6+y_7x_7,\\
A_1=&A_1(x,y)=&-&y_0x_1+y_1x_0-y_2x_4+y_4x_2-y_3x_7+y_7x_3-y_5x_6+y_6x_5,\\
A_2=&A_2(x,y)=&-&y_0x_2+y_2x_0-y_4x_1+y_1x_4-y_3x_5+y_5x_3-y_6x_7+y_7x_6,\\
A_3=&A_3(x,y)=&-&y_0x_3+y_3x_0+y_1x_7-y_7x_1+y_2x_5-y_5x_2-y_4x_6+y_6x_4,\\
A_4=&A_4(x,y)=&-&y_0x_4+y_4x_0-y_1x_2+y_2x_1+y_3x_6-y_6x_3-y_5x_7+y_7x_5,\\
A_5=&A_5(x,y)=&-&y_0x_5+y_5x_0+y_1x_6-y_6x_1-y_2x_3+y_3x_2+y_4x_7-y_7x_4,\\
A_6=&A_6(x,y)=&-&y_0x_6+y_6x_0-y_1x_5+y_5x_1+y_2x_7-y_7x_2-y_3x_4+y_4x_3,\\
A_7=&A_7(x,y)=&-&y_0x_7+y_7x_0-y_1x_3+y_3x_1-y_2x_6+y_6x_2-y_4x_5+y_5x_4,\\
B_0=&B_0(x,y)=&&x_0^2+x_1^2+x_2^2+x_3^2+x_4^2+x_5^2+x_6^2+x_7^2,\\
B_1=&B_1(x,y)=&&y_0^2+y_1^2+y_2^2+y_3^2+y_4^2+y_5^2+y_6^2+y_7^2.
\end{matrix}
\end{displaymath}

Then the embedding functions of a Zariski open subset  ${\mathcal A}$,
which is identified with ${\mathbb C}^{16}$ with coordinates
$z:=(x_0,\cdots,x_7, y_0,\cdots, y_7)$, of
$M_{16}:=\frac{E_6}{SO(10)\times SO(2)}$
 into $\mathbb{CP}^{26}$
are given by:
$$z\mapsto [1,x_0,x_1,x_2,x_3,x_4,x_5,x_6,x_7,y_0,y_1,y_2,y_3,y_4,y_5,y_6,y_7,A_0,A_1,A_2,A_3,A_4,A_5,A_6,A_7,B_0,B_1].$$
\medskip

$\bf \clubsuit 2.$ {\bf Case $\bf M_{27}$:} Similarly we define
\begin{displaymath}
\begin{matrix}
x&=&&(x_1,x_2,x_3),\\
y&=&&(y_0,y_1,y_2,y_3,y_4,y_5,y_6,y_7),\\
t&=&&(t_0,t_1,t_2,t_3,t_4,t_5,t_6,t_7),\\
\omega&=&&(\omega_0,\omega_1,\omega_2,\omega_3,\omega_4,\omega_5,\omega_6,\omega_7).\\
\end{matrix}
\end{displaymath}

Define functions $A,B,C,D_0,\dots,D_7,E_0\dots,E_7,F_0\dots,F_7$ and $G$ such that,
\begin{displaymath}{\rm Com}(X)=X\times X=\left(
\begin{matrix}A&D&\overline E\\
\overline D&B&F\\
 E&\overline F&C\\
\end{matrix}\right),\,\,\,G={\rm det}(X),
\end{displaymath}
where $D=\sum_{j=0}^{7}D_je_j$, $E=\sum_{j=0}^{7}E_je_j$, $F=\sum_{j=0}^{7}F_je_j$ and the matrix $X$ corresponding to the point  $(x,y,t,w)\in\mathbb C^{27}$ is given by
\begin{displaymath}X=\left(
\begin{matrix}x_1&y&\overline t\\
\overline y&x_2&w\\
t&\overline w&x_3\\
\end{matrix}\right)\in\mathcal J_3(\mathbb O). 
\end{displaymath}
Recall the formulas in \cite{O}, we have
\begin{displaymath}X\times X=\left(
\begin{matrix}x_2x_3-w\overline w&\overline w \overline t-x_3y&yw-x_2\overline t\\
wt-x_3\overline y&x_3x_1-t\overline t&\overline t\overline y-x_1w\\
\overline y\overline w-x_2t&ty-x_1\overline w&x_1x_2-y\overline y\\
\end{matrix}\right)\in\mathcal J_3(\mathbb O),
\end{displaymath}
\begin{displaymath}
{\rm det}(X)=x_1x_2x_3-x_1w\overline w-x_2 t\bar t-x_3y\overline y+2 \Re^c(wty),
\end{displaymath}
where $ \Re^c(x)=x_0$ for any $ x=\sum_{i=0}^{7}x_ie_i\in\mathbb O.$

By further computation, we have the explicit expressions as follows:
\begin{displaymath}
\begin{matrix}
A=A(x,y,t,\omega)&=&&x_2x_3-(\omega_0^2+\omega_1^2+\omega_2^2+\omega_3^2+\omega_4^2+\omega_5^2+\omega_6^2+\omega_7^2),\\
B=B(x,y,t,\omega)&=&&x_1x_3-(t_0^2+t_1^2+t_2^2+t_3^2+t_4^2+t_5^2+t_6^2+t_7^2),\\
C=C(x,y,t,\omega)&=&&x_1x_2-(y_0^2+y_1^2+y_2^2+y_3^2+y_4^2+y_5^2+y_6^2+y_7^2),\\
D_0=D_0(x,y,t,\omega)&=&&t_0\omega_0+t_1\omega_1+t_2\omega_2+t_3\omega_3+t_4\omega_4+t_5\omega_5+t_6\omega_6+t_7\omega_7-x_3y_0,\\
D_1=D_1(x,y,t,\omega)&=&-&t_0\omega_1+t_1\omega_0-t_2\omega_4+t_4\omega_2-t_3\omega_7+t_7\omega_3-t_5\omega_6+t_6\omega_5-x_3y_1,\\
D_2=D_2(x,y,t,\omega)&=&-&t_0\omega_2+t_2\omega_0-t_4\omega_1+t_1\omega_4-t_3\omega_5+t_5\omega_3-t_6\omega_7+t_7\omega_6-x_3y_2,\\
D_3=D_3(x,y,t,\omega)&=&-&t_0\omega_3+t_3\omega_0+t_1\omega_7-t_7\omega_1+t_2\omega_5-t_5\omega_2-t_4\omega_6+t_6\omega_4-x_3y_3,\\
D_4=D_4(x,y,t,\omega)&=&-&t_0\omega_4+t_4\omega_0-t_1\omega_2+t_2\omega_1+t_3\omega_6-t_6\omega_3-t_5\omega_7+t_7\omega_5-x_3y_4,\\
D_5=D_5(x,y,t,\omega)&=&-&t_0\omega_5+t_5\omega_0+t_1\omega_6-t_6\omega_1-t_2\omega_3+t_3\omega_2+t_4\omega_7-t_7\omega_4-x_3y_5,\\
D_6=D_6(x,y,t,\omega)&=&-&t_0\omega_6+t_6\omega_0-t_1\omega_5+t_5\omega_1+t_2\omega_7-t_7\omega_2-t_3\omega_4+t_4\omega_3-x_3y_6,\\
D_7=D_7(x,y,t,\omega)&=&-&t_0\omega_7+t_7\omega_0-t_1\omega_3+t_3\omega_1-t_2\omega_6+t_6\omega_2-t_4\omega_5+t_5\omega_4-x_3y_7,\\
E_0=E_0(x,y,t,\omega)&=&&y_0\omega_0-y_1\omega_1-y_2\omega_2-y_3\omega_3-y_4\omega_4-y_5\omega_5-y_6\omega_6-y_7\omega_7-x_2t_0,\\
E_1=E_1(x,y,t,\omega)&=&&y_0\omega_1+y_1\omega_0+y_2\omega_4-y_4\omega_2+y_3\omega_7-y_7\omega_3+y_5\omega_6-y_6\omega_5-x_2t_1,\\
E_2=E_2(x,y,t,\omega)&=&&y_0\omega_2+y_2\omega_0+y_4\omega_1-y_1\omega_4+y_3\omega_5-y_5\omega_3+y_6\omega_7-y_7\omega_6-x_2t_2,\\
E_3=E_3(x,y,t,\omega)&=&&y_0\omega_3+y_3\omega_0-y_1\omega_7+y_7\omega_1-y_2\omega_5+y_5\omega_2+y_4\omega_6-y_6\omega_4-x_2t_3,\\
E_4=E_4(x,y,t,\omega)&=&&y_0\omega_4+y_4\omega_0+y_1\omega_2-y_2\omega_1-y_3\omega_6+y_6\omega_3+y_5\omega_7-y_7\omega_5-x_2t_4,\\
E_5=E_5(x,y,t,\omega)&=&&y_0\omega_5+y_5\omega_0-y_1\omega_6+y_6\omega_1+y_2\omega_3-y_3\omega_2-y_4\omega_7+y_7\omega_4-x_2t_5,\\
E_6=E_6(x,y,t,\omega)&=&&y_0\omega_6+y_6\omega_0+y_1\omega_5-y_5\omega_1-y_2\omega_7+y_7\omega_2+y_3\omega_4-y_4\omega_3-x_2t_6,\\
E_7=E_7(x,y,t,\omega)&=&&y_0\omega_7+y_7\omega_0+y_1\omega_3-y_3\omega_1+y_2\omega_6-y_6\omega_2+y_4\omega_5-y_5\omega_4-x_2t_7,\\
F_0=F_0(x,y,t,\omega)&=&&y_0t_0+y_1t_1+y_2t_2+y_3t_3+y_4t_4+y_5t_5+y_6t_6+y_7t_7-x_1\omega_0,\\
F_1=F_1(x,y,t,\omega)&=&&y_0t_1-y_1t_0-y_2t_4+y_4t_2-y_3t_7+y_7t_3-y_5t_6+y_6t_5-x_1\omega_1,\\
F_2=F_2(x,y,t,\omega)&=&&y_0t_2-y_2t_0-y_4t_1+y_1t_4-y_3t_5+y_5t_3-y_6t_7+y_7t_6-x_1\omega_2,\\
F_3=F_3(x,y,t,\omega)&=&&y_0t_3-y_3t_0+y_1t_7-y_7t_1+y_2t_5-y_5t_2-y_4t_6+y_6t_4-x_1\omega_3,\\
F_4=F_4(x,y,t,\omega)&=&&y_0t_4-y_4t_0-y_1t_2+y_2t_1+y_3t_6-y_6t_3-y_5t_7+y_7t_5-x_1\omega_4,\\
F_5=F_5(x,y,t,\omega)&=&&y_0t_5-y_5t_0+y_1t_6-y_6t_1-y_2t_3+y_3t_2+y_4t_7-y_7t_4-x_1\omega_5,\\
F_6=F_6(x,y,t,\omega)&=&&y_0t_6-y_6t_0-y_1t_5+y_5t_1+y_2t_7-y_7t_2-y_3t_4+y_4t_3-x_1\omega_6,\\
F_7=F_7(x,y,t,\omega)&=&&y_0t_7-y_7t_0-y_1t_3+y_3t_1-y_2t_6+y_6t_2-y_4t_5+y_5t_4-x_1\omega_7.\\
\end{matrix}
\end{displaymath}

\begin{align*}
G=G(x,y,t,\omega)&=x_1x_2x_3-x_1(\omega_0^2+\omega_1^2+\omega_2^2+\omega_3^2+\omega_4^2+\omega_5^2+\omega_6^2+\omega_7^2)\\
\quad&-x_2(t_0^2+t_1^2+t_2^2+t_3^2+t_4^2+t_5^2+t_6^2+t_7^2)\\
\quad&-x_3(y_0^2+y_1^2+y_2^2+y_3^2+y_4^2+y_5^2+y_6^2+y_7^2)\\
\quad&+2\{(y_0\omega_0-y_1\omega_1-y_2\omega_2-y_3\omega_3-y_4\omega_4-y_5\omega_5-y_6\omega_6-y_7\omega_7)t_0\\
\quad&+(y_0\omega_1+y_1\omega_0+y_2\omega_4-y_4\omega_2+y_3\omega_7-y_7\omega_3+y_5\omega_6-y_6\omega_5)t_1\\
\quad&+(y_0\omega_2+y_2\omega_0+y_4\omega_1-y_1\omega_4+y_3\omega_5-y_5\omega_3+y_6\omega_7-y_7\omega_6)t_2\\
\quad&+(y_0\omega_3+y_3\omega_0-y_1\omega_7+y_7\omega_1-y_2\omega_5+y_5\omega_2+y_4\omega_6-y_6\omega_4)t_3\\
\quad&+(y_0\omega_4+y_4\omega_0+y_1\omega_2-y_2\omega_1-y_3\omega_6+y_6\omega_3+y_5\omega_7-y_7\omega_5)t_4\\
\quad&+(y_0\omega_5+y_5\omega_0-y_1\omega_6+y_6\omega_1+y_2\omega_3-y_3\omega_2-y_4\omega_7+y_7\omega_4)t_5\\
\quad&+(y_0\omega_6+y_6\omega_0+y_1\omega_5-y_5\omega_1-y_2\omega_7+y_7\omega_2+y_3\omega_4-y_4\omega_3)t_6\\
\quad&+(y_0\omega_7+y_7\omega_0+y_1\omega_3-y_3\omega_1+y_2\omega_6-y_6\omega_2+y_4\omega_5-y_5\omega_4)t_7\}.\\
\end{align*}
Hence the embedding functions of a Zariski open subset  ${\mathcal A}$,
which is identified with ${\mathbb C}^{27}$ with coordinates
$z:=(x,y,t,\omega)=(x_1,x_2,x_3,y_0\cdots,y_7, t_0,\cdots,
t_7,\omega_0,\cdots,\omega_7)$, of $M_{27}:=\frac{E_7}{E_6\times
SO(2)}$
 into $\mathbb{CP}^{55}$
are given by:
$z\mapsto [1, x, y, t, \omega,
A, B, C, D_0,D_1,D_2,D_3,
D_4,D_5,D_6,D_7,E_0,E_1,E_2,E_3,E_4,E_5,$
$E_6,E_7,F_0,F_1,F_2,F_3,F_4,F_5,F_6,F_7,G].$
The detailed discussions related to this Appendix can be found in
\cite{CMP}, \cite{Fr} and \cite{O}.

\bigskip\bigskip
\section {Appendix II:\ Proof of Proposition (I) for other
types}

In this Appendix, we complete the proof of Proposition (I) for
spaces of the other type.
\subsection{Spaces of type II}

In this subsection, we establish Proposition (I) for the orthogonal
Grassmannians $G_{II}(n,n)$. As shown in $\S 2$, we have a Zariski
open affine chart $\mathcal{A} \subset G_{II}(n,n)$  of elements of the form:
\begin{displaymath}
\left(
\begin{matrix}I_{n\times n}&Z
\end{matrix}
\right)= \left(
\begin{matrix}1&0&0&\cdot\cdot\cdot&0&0&z_{12}&\cdot\cdot\cdot&z_{1n}\\
0&1&0&\cdot\cdot\cdot&0&-z_{12}&0&\cdot\cdot\cdot&z_{2n}\\
&&&\cdot\cdot\cdot&&&\cdot\cdot\cdot&\\
0&0&0&\cdot\cdot\cdot&1&-z_{1n}&-z_{2n}&\cdot\cdot\cdot&0\\
\end{matrix}
\right)
\end{displaymath}
Here $z=(z_{12},z_{13},...,z_{(n-1)n})$ is the local coordinates for
$\mathcal A \cong \mathbb C^{\frac{n(n-1)}{2}}.$
Its conjugate $\mathcal A^*\subset
(G_{II}(n,n))^*$ is also a copy of  $\mathbb C^{\frac{n(n-1)}{2}}.$
We write the local coordinates  for $A^{*}$ as
$\xi=(\xi_{12},...,\xi_{(n-1)n}).$

The canonical embedding is given by
\begin{equation*}
(1,
,...,\mathrm{pf} (Z_{\sigma}),...).
\end{equation*}
The defining function for the  Segre family  (in the product of such
affine pieces) is given by
\begin{equation*}\rho(z,\xi)=1+\sum_
{\begin{subarray}{c}\sigma \in S_{k},\\
2 \leq k \leq n ,2| k
 \end{subarray}}
{\rm{Pf}}(Z_{\sigma}){\rm{Pf}}
   (\Xi_{\sigma})
\end{equation*}

Write
\begin{equation}
r_{Z}= \Big({\mathrm{Pf}} (Z_{\sigma})_{\sigma \in S_{k}}\Big)_{2
\leq k \leq n ,2| k}.
\end{equation}









The local biholomorphic  map $F$ defined near $0\in U$ with $F(0)=0$
can be represented as a matrix:

$$F=\left(
      \begin{array}{cccc}
        0 & f_{12} & ... & f_{1n} \\
        -f_{12} & 0 & ... & f_{2n} \\
        ... & ... & ... & ...\\
        -f_{1n} & ...& ... & 0 \\
      \end{array}
    \right).\ \
$$
Let $r_{F}$  be
\begin{equation}\label{eq71}
r_{F}=\Big(\mathrm{pf} ((F)_{\sigma})_{\sigma \in S_{k}}\Big)_{2
\leq k \leq n ,2| k}.
\end{equation}
Under the notation of $\S 2$, it is easy to see
$r_{Z}=(\psi_{1},...,\psi_{N}),$
$r_{F}=(\psi_{1}(F),...,\psi_{N}(F)).$

We write $\widetilde{z}$ for the $z$ with the last component
$z_{(n-1)n}$ dropped. More precisely,
\begin{equation}\label{eq72}
\widetilde{z}=(z_{12},...,z_{1n},z_{23},...,z_{2n},...,z_{(n-2)(n-1)},z_{(n-2)n}),
\end{equation}
Recall $z$ has $n'=n(n-1)/2$ independent variables. Thus
$\widetilde{z}$ has $(n'-1)$ components. We define the
$\widetilde{z}-$rank and  $\widetilde{z}-$ nondegeneracy as in
Definition \ref{df21}  using $\psi=r_{F}$ in (\ref{eq71}) and
$\widetilde{z}$ as in (\ref{eq72}) with $m=n',$ respectively.  We
now prove the following:
\begin{proposition}\label{thm45}
$r_{F}$ is $\widetilde{z}-$nondegenerate near $0.$ More precisely,
$\mathrm{rank}_{1+N-n'}(r_{F}, \widetilde{z})=N.$
\end{proposition}
{\it Proof of Proposition \ref{thm45}}: Suppose not. Without loss of
generality, we  assume that
$$\mathrm{rank}_{1+N-n'}(r_{F}, \widetilde{z}) < N.$$
As a consequence of Theorem \ref{thm211}, there exist $c_{\sigma,k}
\in \mathbb{C}, 4 \leq k \leq n, 2|k, \sigma \in S_{k}$, which are
not all zero, such that
\begin{equation}\label{eq73}
\sum_{4 \leq k \leq n, 2|n} \sum_{\sigma \in S_{k}}
c_{\sigma,k}~{\rm pf}((F)_{\sigma})(z_{12},...,z_{(n-2)n},0)) \equiv
0.
\end{equation}

However, (\ref{eq73})  cannot hold by the following lemma, which
gives a contradiction:

\begin{lemma}\label{lemma46}
Let
$$H=\left(
      \begin{array}{cccc}
        0 & h_{12} & ... & h_{1n} \\
        -h_{12} & 0 & ... & ... \\
        ... & ... & ... & ... \\
        -h_{1n} & ... & ... & 0 \\
      \end{array}
    \right)$$
be an anti-symmetric matrix-valued holomorphic function in a
neighborhood $U$ of $0$ in $\tilde{z}=(z_{12},...,z_{(n-2)n)}) \in
\mathbb{C}^{n'-1}$ with $H(0)=0.$ Assume that $H$ is of full rank at
$0.$ Set $r_{H}$ similar to the definition of $r_{F},$
\begin{equation}
r_{H}=\Big(\mathrm{pf} (H_{\sigma})_{\sigma \in S_{k}}\Big)_{2 \leq
k \leq n ,2| k}.
\end{equation}
Assume that $a_{\sigma,k}, \sigma \in S_{k}, 4 \leq k \leq n, $ are
complex numbers such that
\begin{equation}
\sum_{4 \leq k \leq n, 2|k} \sum_{\sigma \in S_{k}}
a_{\sigma,k}~{\rm pf}(H_{\sigma})(z_{12},...,z_{(n-2)n})) \equiv
0~\text{for all}~\widetilde{z} \in U.
\end{equation}
Then $$a_{\sigma, k}=0$$ for all $\sigma \in S_{k}, 4 \leq k \leq n,
2|k.$
\end{lemma}
{\it Proof of Lemma \ref{lemma46}}: Suppose not. We will prove the
lemma by seeking a contradiction. Note that $H$ has full rank at
$0.$ Hence there exist $(n'-1)$ components $\widehat{H}$ of $H$ that
forms a local biholomorphism from $\mathbb{C}^{n'-1}$ to
$\mathbb{C}^{n'-1}.$ We assume that these $(n'-1)$ components
$\widehat{H}$ are $H$ with $h_{i_{0}j_{0}}$ being dropped,  where
$i_{0} < j_{0}.$ Without loss of generality, we assume $i_{0}=n-1,
j_{0}=n.$  By a local biholomorphic change of coordinates, we assume
$\widehat{H}=\widetilde{z}=(z_{12},...,z_{(n-2)n}).$ We still write the missing
component as $h_{(n-1)n}.$ Now we assume $2(m+1), m \geq 1,$ is the
least number $k$ such that $\{a_{\sigma, k}\}_{\sigma \in S_{k}}$
are not all zero. We then consider $\{a_{\sigma, 2(m+1)}\}_{\sigma
\in S_{2(m+1)}}.$ We first claim that $a_{\sigma, 2(m+1)}=0$ for
those $\sigma \in S_{2(m+1)}$ such that $\mathrm{pf} (H_{\sigma})$
involves $h_{(n-1)n}.$ More precisely, if $\mathrm{pf} (H_{\sigma}),
\sigma \in S_{2(m+1)}$ involves $h_{(n-1)n},$ then $\sigma=\{
i_{1},...,i_{2m}, (n-1), n\}$ for some $1 \leq i_{1}<...<i_{2m} \leq
n-2.$ Suppose its coefficient is not zero. Then $\mathrm{pf}
(H_{\sigma})$ will produce the monomial
$z_{i_{1}i_{2}}z_{i_{3}i_{4}}...z_{i_{2m-3}i_{2m-2}}z_{i_{2m-1}(n-1)}z_{i_{2m}n}.$
This term can only be canceled by the terms of form:
$z_{i_{2m-1}(n-1)}h_{(n-1)n}Q$ or $z_{i_{2m}n}h_{(n-1)n}Q.$ But
neither of them can appear in any other Pfaffians. This is a
contradiction. Once we know there are no $h_{(n-1)n}$ involved, then
the remaining Pfaffians have only  terms consisting of   the product
of some of  $z_{12},...,z_{(n-2)n}.$ Their sum cannot be zero unless
their coefficients are all zero. This is a contradiction. We thus
establishes Lemma \ref{lemma46}.\ $\endpf$

We thus also get a contradiction to equation (\ref{eq73}). This
establishes Proposition \ref{thm45}.  $\endpf$

\medskip

\begin{remark}
By Proposition \ref{thm45}, there exist multiindices
$\tilde{\beta}^1,...,\tilde{\beta}^N$ with all $|\tilde{\beta}^j|
\leq 1+N-n',$ and there is a point
$$z^0=\left(
        \begin{array}{ccccc}
          0 & z^0_{12} & ... & z^0_{1(n-1)} & z^0_{1n} \\
          -z^0_{12} & 0 & ... & z^0_{2(n-1)} & z^0_{2n} \\
          ... & ... & ... & ... & ... \\
          -z^0_{1(n-1)} & -z^0_{2(n-1)} & ... & 0 & z^0_{(n-1)n} \\
          -z^0_{1n} & -z^0_{2n} & ... & -z^0_{(n-1)n} & 0 \\
        \end{array}
      \right), z^0_{(n-1)n} \neq 0; $$
near $0$ such that
\begin{equation}
\left| \begin{matrix}
\frac{\partial^{|\beta^1|} (\psi_{1}(F))}{\partial \tilde{z}^{\tilde{\beta}^1}} & ...  & \frac{\partial^{|\beta^1|} (\psi_{N}(F))}{\partial \tilde{z}^{\tilde{\beta}^1}} \\
... & ... & ... \\
\frac{\partial^{|\beta^N|} (\psi_{1}(F))}{\partial
\tilde{z}^{\tilde{\beta}^N}} & ... & \frac{\partial^{|\beta^N|}
(\psi_{N}(F))}{\partial \tilde{z}^{\tilde{\beta}^N}}
\end{matrix}\right|(z^0) \neq 0.
\end{equation}
\end{remark}

We set
$$
\xi^0=\left(
        \begin{array}{ccccc}
          0 & 0 & ... & 0 & 0 \\
          0 & 0 & ... & 0 & 0 \\
          ... & ... & ... & ... & ... \\
          0 & 0 & ... & 0  & \xi^0_{(n-1)n} \\
          0 & 0 & ... & -\xi^0_{(n-1)n} & 0 \\
        \end{array}
      \right)
 \in \mathbb{C}^{n^2}, \xi^0_{(n-1)n}=-\frac{1}{z_{(n-1)n}^0}.$$
Then it is easy to see that $(z^0, \xi^0) \in
\mathcal{M}=\{\rho(z,\xi)=0\}.$

Write for each $1 \leq i < j \leq n, (i,j) \neq (n-1,n),$
\begin{equation}
\mathcal{L}_{ij}=\frac{\partial}{\partial z_{ij}}-
\frac{\frac{\partial \rho}{\partial z_{ij}}(z,\xi)}{\frac{\partial
\rho}{\partial z_{(n-1)n}}(z,\xi)}\frac{\partial}{\partial
z_{(n-1)n}}
\end{equation}
which are holomorphic tangent vector fields  along $\mathcal{M}$
near $(z^0, \xi^0).$ Here we note that $\frac{\partial
\rho}{\partial z_{(n-1)n}}(z,\xi)$ is nonzero near $(z^0, \xi^0).$
For any $(n'-1)$-multiindex $\beta=(\beta_{12},...,\beta_{(n-2)n})$,
we write
$$\mathcal{L}^{\beta}=\mathcal{L}_{12}^{\beta_{12}}...\mathcal{L}_{(n-2)n}^{\beta_{(n-2)n}}.$$

Now we define for any $N$ collection of $(n'-1)-$multiindices
$\{\beta^1,...,\beta^N\},$
\begin{equation}
\Lambda(\beta^1,...,\beta^N)(z,\xi):=\left|\begin{matrix}
 \mathcal{L}^{\beta^1}(\psi_{1}(F)) & ...  & \mathcal{L}^{\beta^1}(\psi_{N}(F)) \\
... & ...  & ...  \\
\mathcal{L}^{\beta^N}(\psi_{1}(F))  & ...  &
\mathcal{L}^{\beta^N}(\psi_{N}(F))
\end{matrix}\right|(z,\xi).
\end{equation}

Note that for any multiindex $\beta, \mathcal{L}^{\beta}$ evaluating
at $(z^0,\xi^0)$ coincides with $\frac{\partial}{\partial
\tilde{z}^{\beta}}.$ We thus again have
\medskip
\begin{theorem}
There exists multiindices $\{{\beta}^1,...,{\beta}^N\},$ such that
$$\Lambda({\beta}^1,...,{\beta}^N)(z,\xi) \neq 0,$$
for $(z,\xi)$ in a small neighborhood of $(z^0, \xi^0)$ and
$\beta^1=(0,...,0).$
\end{theorem}

\subsection{Spaces of type III}

Let $F$ be a local biholomorphic map at $0$. In this case, both $Z$ and $F$ are $n \times
n$ symmetric matrices. We write
$$Z=\left(
      \begin{array}{cccc}
        z_{11} & z_{12} &... & z_{1n} \\
        z_{12} & z_{22} & ... & z_{2n} \\
        ... & ... & ... & ... \\
        z_{1n} & z_{2n} & ... & z_{nn} \\
      \end{array}
    \right),\ \ z=(z_{11},z_{12},z_{13},...,z_{nn}).
$$
Similar notations are used for  $F.$

Recall from $\eqref{3r}$ of $\clubsuit 3$ in $\S 2$:
\begin{equation}\label{eq77}
r_{z}=\left(\psi_{1}^1(z),...,\psi_{N_{1}}^1(z),\psi_{1}^2(z),...,\psi_{N_{2}}^2(z),...,\psi_{1}^{n-1}(z),...,\psi_{N_{n-1}}^{n-1}(z),
\psi^{n}(z) \right),
\end{equation}
where $\psi_{j}^k$ is a homogeneous polynomial of degree $k, 1 \leq
j \leq N_{k}. \ \psi^n$ is a homogeneous polynomial of degree $n.$
Moreover, the components of $r_{z}$ are linearly independent.

We write the number of components in $r_{z}$ to be
$N=N_{1}+...+N_{n},$ where we set $N_{n}=1.$ We will also sometimes
write $\psi_{N_{n}}^n=\psi^{n}.$

We emphasize that for each fixed $k, \psi_{1}^k,...,\psi_{N_{k}}^k$
are linearly independent. Moreover,  each $\psi_{j}^k$ is a certain
linear combination of the determinants of  $k \times k$ submatrices
of $Z.$ This will be crucial for our argument later.

\bigskip

We define $r_{F}$ as the composition of $r_z$ with the map $F$:
\begin{equation}\label{eq78}
r_{F}=\left(\psi_{1}^1(F),...,\psi_{N_{1}}^1(F),\psi_{1}^2(F),...,\psi_{N_{2}}^2(F),...,\psi_{1}^{n-1}(F),...,
\psi_{N_{n-1}}^{n-1}(F), \psi^{n}(F) \right).
\end{equation}
In what follows, we write also $z_{ij}=z_{ji}.$ We write $\mathrm{det}(A)$ as the determinant of $A$ when $A$ is a square matrix.

Let $P, \widetilde{P}$ be monomials in $z_{ij}'$s, and $h$ a polynomial in $z_{ij}'$s. Let $a, b$ be two complex numbers.  In the following lemmas, when we say $h$ always
has the terms $aP, b\widetilde{P}$, we mean $h$ has the term  $aP$ if and only if it has the term
$b\widetilde{P}.$

\begin{lemma}\label{lemma51}
Fixing $1 \leq i, j < n,$ let $P=z_{in}z_{nj}Q$ and $\widetilde{P}=z_{ij}z_{nn}Q$ with $Q$  a monomial in $z_{ij}'$s.  The following statements are true.

\begin{itemize}
\item Let $A$ be a square submatrix of $Z.$
If $z_{ij}\nmid Q,$ then $\mathrm{det}(A)$ always has monomials of the form $cP,-c\widetilde{P}$ for some $c
\in \mathbb{C}$ depending on the submatrix $A$. (If $\mathrm{det}(A)$ does not have any multiple of $P$, it does not have any multiple of $\widetilde{P},$ either; vice versa). If $z_{ij} | Q,$ then $\mathrm{det}(A)$ always
has monomials $cP, -(c/2)\widetilde{P}$ for some $c \in \mathbb{C}$
depending on $A$.

\item Let $k \geq 1.$ Let $\psi_l^k(z)$ be as defined in (\ref{eq77}), $1 \leq l \leq N_k.$
If $z_{ij}\nmid Q,$
then $\psi_{l}^k(z)$ always has monomials $\lambda
P,-\lambda\widetilde{P}$ for some $\lambda \in \mathbb{C},$ If
$z_{ij} | Q,$ then $\psi_{l}^k(z)$ always has monomials $\lambda P,
-(\lambda /2)\widetilde{P}$ for some $\lambda \in \mathbb{C}.$

\end{itemize}

\end{lemma}
{\it Proof of Lemma \ref{lemma51}}: The first part is a consequence of the
Laplace expansion of a determinant by complementary minors. The second part is due to the fact that
$\psi_{j}^k$ is a linear combination of the determinants of
submatrices of $Z$ of order $k.$
$\endpf$

\medskip

Similarly, one can prove  in a similar way Lemmas
\ref{lemma53}-\ref{lemma56}.

\begin{lemma}\label{lemma53}
Fixing $1 \leq j < n-1,$ let $P=z_{jn}z_{(n-1)(n-1)}Q$
and  $\widetilde{P}=z_{j(n-1)}z_{(n-1)n}Q$ with $Q$ a monomial in $z_{ij}'$s.
\begin{itemize}
\item Let $A$ be a square
submatrix of $Z.$   If $z_{jn}\nmid Q,$ then $\mathrm{det}(A)$ always
has  monomials $cP,-c\widetilde{P}$ for some $c \in \mathbb{C}.$ If
$z_{jn} | Q,$ then $\mathrm{det}(A)$ always has monomials $cP,
-2c\widetilde{P}$ for some $c \in \mathbb{C}.$

\item Let $k \geq 1.$ Let $\psi_l^k(z)$ be as defined in (\ref{eq77}), $1 \leq l \leq N_k.$ If $z_{jn}\nmid Q,$ then $\psi_l^k(z)$ always
has  monomials $\lambda P,-\lambda\widetilde{P}$ for some $\lambda \in \mathbb{C}.$ If
$z_{jn} | Q,$ then $\psi_l^k(z)$ always has monomials $\lambda P,
-2\lambda \widetilde{P}$ for some $\lambda \in \mathbb{C}.$

\end{itemize}
\end{lemma}


\begin{lemma}\label{lemma55}
Fixing  $1 \leq i < n-1,$ let $P=z_{i(n-1)}z_{ni}Q$ and
$\widetilde{P}=z_{ii}z_{(n-1)n}Q$ with $Q$ a monomial in $z_{ij}'$s.

\begin{itemize}

\item Let $A$ be a square submatrix of
$Z.$   If $z_{(n-1)n}\nmid Q,$ then $\mathrm{det}(A)$ always has
monomials $cP,-c\widetilde{P}$ for some $c \in \mathbb{C}.$ If
$z_{(n-1)n} | Q,$ then $\mathrm{det}(A)$ always has monomials $cP,
-(c/2)\widetilde{P}$ for some $c \in \mathbb{C}.$

\item Let $k \geq 1.$ Let $\psi_l^k(z)$ be as defined in (\ref{eq77}), $1 \leq l \leq N_k.$ If $z_{(n-1)n}\nmid Q,$ then $\psi_l^k(z)$ always has
monomials $\lambda P,-\lambda \widetilde{P}$ for some $\lambda \in \mathbb{C}.$ If
$z_{(n-1)n} | Q,$ then $\psi_l^k(z)$ always has monomials $\lambda P,
-(\lambda /2)\widetilde{P}$ for some $\lambda \in \mathbb{C}.$

\end{itemize}

\end{lemma}

\begin{lemma}\label{lemma56}
Fixing $1 \leq i < n-1, 1 \leq j < n-1, i \neq j,$ let
$P_{1}=z_{i(n-1)}z_{nj}Q, P_{2}=z_{in}z_{j(n-1)}Q$, and
$\widetilde{P}=z_{ij}z_{(n-1)n}Q$ with $Q$ a monomial in $z_{ij}'$s.

\begin{itemize}

\item Let $A$ be a square submatrix of
$Z.$   If $z_{ij}\nmid Q, z_{(n-1)n}\nmid Q,$ then $\mathrm{det}(A)$
always has  terms $c_{1} P_{1} + c_{2} P_{2},
-(c_{1}+c_{2})\widetilde{P}$ for some $c_{1}, c_{2} \in \mathbb{C}.$
If $ z_{ij}\nmid Q, z_{(n-1)n} | Q,$ or $z_{ij}| Q, z_{(n-1)n} \nmid
Q,$ then $\mathrm{det}(A)$ always has terms $c_{1} P_{1} + c_{2}
P_{2}, -\frac{c_{1} +c_{2}}{2}\widetilde{P}$ for some $c_{1}, c_{2}
\in \mathbb{C}.$ If $z_{ij} | Q, z_{(n-1)n} | Q,$ then $\mathrm{det} (A)$ always has terms $c_{1} P_{1} + c_{2} P_{2},
-\frac{c_{1} +c_{2}}{4} \widetilde{P}.$

\item Let $k \geq 1.$ Let $\psi_l^k(z)$ be as defined in (\ref{eq77}), $1 \leq l \leq N_k.$ If $z_{ij}\nmid Q$ and $ z_{(n-1)n}\nmid Q$, then $\psi_{l}^k(z)$ always has  terms $\lambda_{1} P_{1} + \lambda_{2}
P_{2}, -(\lambda_{1}+\lambda_{2})\widetilde{P}$ for some
$\lambda_{1}, \lambda_{2} \in \mathbb{C}.$ If $ z_{ij}\nmid Q,
z_{(n-1)n} | Q,$ or $z_{ij}| Q, z_{(n-1)n} \nmid Q,$ then $\psi_{l}^k(z)$ always
has terms $\lambda_{1} P_{1} + \lambda_{2} P_{2}, -\frac{\lambda_{1}
+\lambda_{2}}{2} \widetilde{P}$ for some $\lambda_{1}, \lambda_{2}
\in \mathbb{C}.$ If $z_{ij} | Q, z_{(n-1)n} | Q,$ then $\psi_l^k(z)$ always has terms $\lambda_{1} P_{1} + \lambda_{2} P_{2},
-\frac{\lambda_{1} +\lambda_{2}}{4} \widetilde{P}$ for some
$\lambda_{1}, \lambda_{2} \in \mathbb{C}.$

\end{itemize}

\end{lemma}

We write $\widetilde{z}$ for $z$ with the last components $z_{nn}$
being dropped. More precisely,
\begin{equation}\label{eq79}
\widetilde{z}=(z_{11},...,z_{1n},z_{22},...,z_{2n},...,z_{(n-1)(n-1)},z_{(n-1)n}).
\end{equation}
Recall $z$ has $n'=n(n+1)/2$ independent variables. Thus
$\widetilde{z}$ has $(n'-1)$ components. We define
$\widetilde{z}-$rank and  $\widetilde{z}-$nondegeneracy in the same
way as before,  using $r_{F}$ in (\ref{eq78}) and $\widetilde{z}$ in
(\ref{eq79}) with $m=n'$. We now need to prove the following:
\begin{proposition}\label{thm510}
$r_{F}$ is $\widetilde{z}-$nondegenerate at $0.$ More precisely,
$\mathrm{rank}_{1+N-n'}(r_{F}, \widetilde{z})=N.$
\end{proposition}
{\it Proof of Proposition \ref{thm510}}: Suppose not. Then one
easily verifies that the hypothesis of Theorem \ref{thm211} is
satisfied. As a consequence of Theorem \ref{thm211}, there exist
$c_{j}^k \in \mathbb{C}, 2 \leq k \leq n, 1 \leq j \leq N_{k},$
which are not all zero such that
\begin{equation}
\sum_{k=2}^{n} \sum_{j=1}^{N_{k}}
c_{j}^k\psi_{j}^k(F(z_{11},...,z_{(n-1)n},0))\equiv 0.
\end{equation}
Here as before,  we write $N_{n}=1, \psi_{N_{n}}^n=\psi^{n}.$

Then we just need to show it can not happen by the following lemma:

\begin{lemma}\label{t2}
Let
$$H=\left(
      \begin{array}{cccc}
        h_{11} & h_{12} & ... & h_{1n} \\
        h_{12} & h_{22} & ... & h_{2n} \\
        ... & ... & ... & ... \\
        h_{1n} & ... & ... & h_{nn} \\
      \end{array}
    \right)
$$
be a symmetric matrix-valued holomorphic function near $0$  in
$\widetilde{z}=(z_{11},...,z_{1n},z_{22},...,z_{2n},..., z_{(n-1)n})
\in \mathbb{C}^{n'-1}$ with $H(0)=0.$ Assume that $H$ is of  full
rank at $0.$ Set $r_{H}$ in a similar way as in $(36):$
$$r_{H}=\left(\psi_{1}^1(H),...,\psi_{N_{1}}^1(H),\psi_{1}^2(H),...,\psi_{N_{2}}^2(H),...,\psi_{1}^{n-1}(H),...,
\psi_{N_{n-1}}^{n-1}(H), \psi^{n}(H) \right)$$ Again we  write
$N_{n}=1, \psi^n=\psi_{N_{n}}^n.$ Assume that $a_{j}^k, 2 \leq k
\leq n, 1 \leq j \leq  n $ are complex numbers such that
\begin{equation}\label{eq81}
\sum_{k=2}^{n} \sum_{j=1}^{N_{k}} a_{j}^{k} \psi_{j}^k
(H(\widetilde{z})) \equiv 0\ \ ~\text{for }~\widetilde{z} \in U.
\end{equation}
Then
$$a_{j}^k=0$$
for each $2 \leq k \leq n, 1 \leq j \leq N_{k}.$
\end{lemma}
{\it Proof of Lemma \ref{t2}}: Suppose not. We will prove the lemma
by seeking a contradiction. Notice that $H$ has full rank at $0.$
Hence there exist $(n'-1)$ components $\widehat{H}$ of $H$ that
gives a local biholomorphism from $\mathbb{C}^{n'-1}$ to
$\mathbb{C}^{n'-1}.$ We assume these $(n'-1)$ components
$\widehat{H}$ are $H$ with $h_{i_{0}j_{0}}$ being dropped, where
$i_{0} \leq j_{0}.$ Then we split our argument into two parts in
terms of $i_{0}=j_{0}$ or $i_{0} < j_{0}.$

\medskip
{\bf Case I:} Assume that $i_{0}=j_{0}.$ Without loss of generality,
we assume $i_{0}=j_{0}=n.$ By a local biholomorphic change of
coordinates, we assume $\widehat{H}=\widetilde{z}=(z_{11},...,z_{n(n-1)}).$ We
still write the last component as $h_{nn}.$ Now we assume $m \geq 2$
is the least number $k$ such that $\{ a_{1}^k,...,a_{N_{k}}^k\}$ are
not all zero. For any holomorphic $g,$ we define $T_{l}(g)$ to be
the homogeneous part of degree $l$ in the Taylor expansion of $g$ at
$0$. Now the assumption in (\ref{eq81}) yields:
\begin{equation}\label{eq82}
T_{m}\left(\sum_{j=1}^{N_{m}} a_{j}^{m} \psi_{j}^{m}
(H(\widetilde{z}))\right) \equiv 0.
\end{equation}

We first compute

$$\sum_{j=1}^{N_{m}} a_{j}^{m} \psi_{j}^{m} (H)=\sum_{j=1}^{N_{m}} a_{j}^{m} \psi_{j}^{m} (z_{11},...,z_{(n-1)n}, h_{nn})$$
formally. Namely,  we regard $h_{nn}$ as a formal variable and only
conduct formal cancellations. We write formally

\begin{equation}\label{eq83}
\sum_{j=1}^{N_{m}} a_{j}^{m} \psi_{j}^{m} (z_{11},...,z_{(n-1)n},
h_{nn})=P_{1}+ h_{nn}P_{2}.
\end{equation}
Here $P_{1}=P_{1}(z_{11},...,z_{(n-1)n})$ is a homogeneous
polynomial of degree $m,$ and $P_{2}=P_{2}(z_{11},...,z_{(n-1)n})$
is a homogeneous polynomial of degree $m-1.$ We claim $P_{2} \neq
0.$ Otherwise,

$$\sum_{j=1}^{N_{m}} a_{j}^{m} \psi_{j}^{m} (z_{11},...,z_{(n-1)n}, h_{nn})=P_{1}.$$
 This implies that $\sum_{j=1}^{N_{m}} a_{j}^{m} \psi_{j}^{m}
(z_{11},...,z_{(n-1)n}, h_{nn})$ does not depend on $h_{nn}$
formally. Then we can replace $h_{nn}$ by $z_{nn}.$ That is,
\begin{equation}\label{eq84}
\sum_{j=1}^{N_{m}} a_{j}^{m} \psi_{j}^{m} (z_{11},...,z_{(n-1)n},
z_{nn})=\sum_{j=1}^{N_{m}} a_{j}^{m} \psi_{j}^{m}
(z_{11},...,z_{(n-1)n}, h_{nn}(\widetilde{z}))=P_{1}.
\end{equation}

By (\ref{eq82}), we see that (\ref{eq84}) is identically zero. This
is a contradiction to the fact that $\{
\psi_{1}^m,...,\psi_{N_{m}}^m\}$ is linearly independent.

Now since $P_{2} \neq 0,$ thus by (\ref{eq83}), $\sum_{j=1}^{N_{m}}
a_{j}^{m} \psi_{j}^{m} (z_{11},...,z_{(n-1)n}, h_{nn})$ has a
monomial of the form $\mu \widetilde{P}=\mu z_{ij}h_{nn} Q$ of
degree $m$ for some  $ 1 \leq i,j < n, \mu \neq 0$ and some monomial
$Q.$ By Lemma \ref{lemma51}, we get that $\sum_{j=1}^{N_{m}}
a_{j}^{m} \psi_{j}^{m} (z_{11},...,z_{(n-1)n}, h_{nn})$  has also
the term $-\mu P$ or $-2 \mu P,$ where $P=z_{in}z_{nj}Q.$ This is a
contradiction to (\ref{eq82}). Indeed, $P$ can be only canceled by
the terms  of the forms: $z_{in}h_{nn}\widetilde{Q}$ or
$z_{nj}h_{nn}\widetilde{Q},$ where $\widetilde{Q}$ is of degree
$m-2.$ But they cannot appear in determinant of any submatrix of $H$
as $z_{in}$(or $z_{nj}$) can not appear with $h_{nn}.$

\bigskip

{\bf Case II:} Assume that $i_{0} \neq j_{0}.$ Without loss of
generality, we assume $i_{0}=(n-1),j_{0}=n.$ Then
$\widehat{H}=(h_{11},...,h_{(n-1)(n-1)},h_{nn})$ is a local
biholomorphism. By a local biholomorphic change of coordinates, we
assume $\widehat{H}=\widetilde{z}=(z_{11},...,z_{(n-1)n}).$ We will still write
the remaining component as $h_{(n-1)n}=h_{n(n-1)}.$ Note that the
fact we are using only is that $ \{z_{11},...,z_{(n-1)n}\}$ are
independent variables. Hence, to make our notation easier, we will
write
$$\widehat{H}=(z_{11},...,z_{(n-1)n})=(w_{11},...,w_{1n},w_{22},...,w_{2n},...,
w_{(n-1)(n-1)},w_{nn})$$ such that they have the same indices as
$h$'s in $\widehat{H}.$ Now we assume $m$ is the least number $k$
such that $\{ a_{1}^k,...,a_{N_{k}}^k\}$ are not all zero. Then
again assumption (\ref{eq81}) yields that
\begin{equation}
T_{m}\left(\sum_{j=1}^{N_{m}} a_{j}^{m} \psi_{j}^{m}
(H(\widetilde{Z}))\right) \equiv 0.
\end{equation}

Again we formally compute that
\begin{equation}\label{eq86}
\sum_{j=1}^{N_{m}} a_{j}^{m} \psi_{j}^{m} (w_{11},...,h_{(n-1)n},
w_{nn})=Q_{1}+ h_{(n-1)n}Q_{2}.
\end{equation}
Here $Q_{1}=Q_{1}(w_{11},...,w_{(n-1)(n-1)},w_{nn})$ is  a
homogeneous polynomial of degree $m.$  Similarly, we can show that
$Q_{2} \neq 0.$ We claim that (\ref{eq86}) does not have a monomial
of the form $h_{(n-1)n}h_{(n-1)n}Q.$ Otherwise, by Lemma \ref{lemma51}, we get that (\ref{eq86})  has
also a monomial of degree $m$ of the form: $w_{(n-1)(n-1)}w_{nn}Q.$
But note that in (\ref{eq86}) it can be canceled only by $h_{(n-1)n}h_{(n-1)n}Q.$
Then $h_{(n-1)n}$ will have a linear term $w_{(n-1)(n-1)}.$ But then
$h_{(n-1)n}h_{(n-1)n}Q$ will produce the term
$w_{(n-1)(n-1)}w_{(n-1)(n-1)}Q.$ This cannot be canceled out by any
other terms.

Now since $Q_{2} \neq 0,$ (\ref{eq86}) has a monomial of the form
$w_{ij}h_{(n-1)n}Q,$ where $Q$ is another monomial in $w$'s. Here $1
\leq i,j \leq n.$ Moreover, $(i,j)\neq (n-1, n-1), (n-1,n), (n,n-1)$
or $(n,n).$ We first assume $1 \leq i ,j <n-1, i \neq j.$ Then by
Lemma \ref{lemma56} , (\ref{eq86})  has either $P_{1}$ or $P_{2},$
where $P_{1}=w_{i(n-1)}w_{nj}Q, P_{2}=w_{in}w_{j(n-1)}Q.$ Note
$P_{1},P_{2}$ can only be canceled by the terms
$w_{i(n-1)}h_{(n-1)n}Q,
w_{nj}h_{(n-1)n}Q,w_{in}h_{(n-1)n}Q,w_{j(n-1)}h_{(n-1)n}Q.$ So one
of them will appear in (\ref{eq86}). Whichever case it is, by Lemma \ref{lemma51}, \ref{lemma53}, (\ref{eq86}) will
have either $P=w_{ln}w_{(n-1)(n-1)}Q, $ or
$\widehat{P}=w_{l(n-1)}w_{nn}Q$ for some $1 \leq l < n.$ We assume,
for instance, (\ref{eq86}) has the monomial $P.$ Then it also has
the monomial $\widetilde{P}=w_{l(n-1)}h_{(n-1)n}Q$ by Lemma \ref{lemma53}. Note that the only term that can cancel $P$ and appear in some determinant is $w_{ln}h_{n(n-1)}Q.$ Hence $h_{n(n-1)}$ has a linear
$w_{(n-1)(n-1)}$ term. Then $\widetilde{P}$ will have the monomial
$w_{l(n-1)}w_{(n-1)(n-1)}Q,$ which can not be canceled by any other
terms. This is a contradiction. The other cases can be proved
similarly.\ \  $\endpf$

\medskip
This establishes Proposition \ref{thm510}. \ \ $\endpf$
\medskip
\begin{remark}
By Proposition \ref{thm510}, there exist multiindices
$\tilde{\beta}^1,...,\tilde{\beta}^N$ with  $|\tilde{\beta}^j| \leq
1+N-pq,$ and there exist
$$z^0=\left( \begin{matrix}
 z^0_{11} & ...  & z^0_{1n}  \\
 ... & ... & ...  \\
 z^0_{1n} & ...  & z^0_{nn}
\end{matrix}\right),\ z_{nn}^0 \neq 0,$$
near $0$ such that
\begin{equation}
\left| \begin{matrix}
\frac{\partial^{|\beta^1|} (\psi_{1}(F))}{\partial \tilde{Z}^{\tilde{\beta}^1}} & ...  & \frac{\partial^{|\beta^1|} (\psi_{N}(F))}{\partial \tilde{Z}^{\tilde{\beta}^1}} \\
... & ... & ... \\
\frac{\partial^{|\beta^N|} (\psi_{1}(F))}{\partial
\tilde{Z}^{\tilde{\beta}^N}} & ... & \frac{\partial^{|\beta^N|}
(\psi_{N}(F))}{\partial \tilde{Z}^{\tilde{\beta}^N}}
\end{matrix}\right|(z^0) \neq 0.
\end{equation}
Here we simply write $r_{F}=(\psi_{1}(F),...,\psi_{N}(F)).$

\end{remark}

We then set
$$\xi^0=\left( \begin{matrix}
 0 & ...  & 0 & 0  \\
 0 & ...  & 0 & 0  \\
 0 & ... & 0  & \xi_{nn}^0
         \end{matrix}
 \right) \in \mathbb{C}^{n^2},\xi^0_{nn}=-\frac{1}{z_{nn}^0}.$$
It is easy to verify that $(z^0, \xi^0) \in
\mathcal{M}=\{\rho(z,\xi)=0\}.$

Write for each $1 \leq i \leq j \leq n, (i,j) \neq (n,n),$
\begin{equation}
\mathcal{L}_{ij}=\frac{\partial}{\partial z_{ij}}-
\frac{\frac{\partial \rho}{\partial z_{ij}}(z,\xi)}{\frac{\partial
\rho}{\partial z_{nn}}(z,\xi)}\frac{\partial}{\partial z_{nn}},
\end{equation}
which are holomorphic tangent vector fields along $\mathcal{M}$ near
$(z^0, \xi^0).$ Here we note that $\frac{\partial \rho}{\partial
z_{nn}}(z,\xi)$ is nonzero near $(z^0, \xi^0).$ For any
$(n'-1)$-multiindex $\beta=(\beta_{11},...,\beta_{(n-1)n})$, we
write
$$\mathcal{L}^{\beta}=\mathcal{L}_{11}^{\beta_{11}}...\mathcal{L}_{(n-1)n}^{\beta_{(n-1)n}}.$$

Now we define for any $N$ collection of $(n'-1)-$multiindices
$\{\beta^1,...,\beta^N\},$
\begin{equation}
\Lambda(\beta^1,...,\beta^N)(z,\xi):=\left|\begin{matrix}
 \mathcal{L}^{\beta^1}(\psi_{1}(F)) & ...  & \mathcal{L}^{\beta^1}(\psi_{N}(F)) \\
... & ...  & ...  \\
\mathcal{L}^{\beta^N}(\psi_{1}(F))  & ...  &
\mathcal{L}^{\beta^N}(\psi_{N}(F))
\end{matrix}\right|(z,\xi).
\end{equation}

Note $\mathcal{L}^{\beta}$ evaluating at $(z^0, \xi^0)$ coincides
with $\frac{\partial}{\partial \tilde{Z}^{\beta}}.$ We have
\medskip
\begin{theorem}
There exists multiindices $\{{\beta}^1,...,{\beta}^N\}$ such that
$\Lambda({\beta}^1,...,{\beta}^N)(z,\xi) \neq 0$ for $(z,\xi)$ in a
small neighborhood of $(z^0, \xi^0)$ and $\beta^1=(0,0,...,0).$
\end{theorem}

\subsection{The exceptional class ${ M}_{27}$}

In this setting, we use the coordinates
$$z=(x_{1},x_{2},x_{3}, y_{0},...,y_{7},t_{0},...,t_{7},w_{0},...,w_{7})\in {\mathbb C}^{27}.$$
The defining function of the Segre family  described in
$\eqref{e27rho}$ is :
$$\rho(z,\xi)=1+r_{z}\cdot r_{\xi}=1+\sum_{i=1}^N\psi_i(z)\psi_i(\xi),\ \ \text {where } N=55 \text{ and }$$
\begin{equation}\label{eqn104}
r_{z}=(x_{1},x_{2},x_{3},y_{0},...,y_{7},t_{0},...,t_{7},w_{0},...,w_{7},A,B,C,D_{0},...D_{7},
E_{0},...,E_{7},F_{0},...,F_{7},G).
\end{equation}
Here $A,B,C, D_{i},E_{i},F_{i}$ are homogeneous quadratic
polynomials in $z$ and   $G$ is a homogeneous cubic polynomial in
$z$:
\begin{equation}
A=x_{2}x_{3}-\sum_{i=0}^{7}w_{i}^2, B=x_{1}x_{3}-\sum_{i=0}^7
t_{i}^2, C=x_{1}x_{2}-\sum_{i=0}^7 y_{i}^2.
\end{equation}
For the expressions for $D_{i},E_{i},F_{i},G,$ see Appendix I. Let
$F$ be a local biholomorphic map near $0$.  We write
$$F=(\phi_{1},\phi_{2},\phi_{3},f_{10},...,f_{17},f_{20},...,f_{27},f_{30},...,h_{37}).$$
Also define $r_{F}$  to be the composition of $r_{z}$ with $F$:
\begin{equation}\label{eqt95}
r_{F}=r_z\circ
F=(\phi_{1},\phi_{2},\phi_{3},f_{10},...,f_{17},f_{20},...,f_{27},f_{30},...,f_{37},A(F),
B(F),C(F),....,G(F)).
\end{equation}
Notice that $r_{F}$ has $55$ components. We will also  write
$$r_{F}=(\psi_{1}(F),...,\psi_{55}(F)).$$
We write $\widetilde{z}$ for $z$ with $x_{3}$ being dropped. Namely,
\begin{equation}\label{eq96}
\widetilde{z}=(x_{1},x_{2},y_{0},...,y_{7},t_{0},...,t_{7},w_{0},...,w_{7}).
\end{equation}

We define the $\widetilde{z}-$rank and  $\psi-$nondegeneracy as in
Definition \ref{df21}  using $r_{F}$ in (\ref{eqt95}) and
$\widetilde{z}$ in (\ref{eq96}) with $m=27$.
\begin{proposition}\label{thm71}
$F$ is $\widetilde{z}-$nondegenerate near $0.$ More precisely,
$\mathrm{rank}_{29}(F, \widetilde{z})=55.$
\end{proposition}
{\it Proof of Proposition \ref{thm71}}: Suppose not. As a
consequence of Theorem \ref{thm211}, there exist $c_{1},...,c_{28}
\in \mathbb{C}$ that are not all zero, such that
\begin{equation}\label{eq91}
c_{1}A(F(x_{1},x_{2},0,y_{0},...,w_{7}))+...+c_{28}G(F(x_{1},x_{2},0,y_{0},...,w_{7}))
\equiv 0.
\end{equation}

We will  show that (\ref{eq91}) cannot hold by the following lemma:
\begin{lemma}\label{t3}
Let
$H=(\psi_{1},\psi_{2},\psi_{3},h_{10},...,h_{17},h_{20},...,h_{27},h_{30},...,h_{37})$
be a vector-valued holomorphic function in a neighborhood $U$ of $0$
in $\tilde{z}= (x_{1},x_{2},
y_{0},...,y_{7},t_{0},...,t_{7},w_{0},...,w_{7})\in \mathbb{C}^{26}$
with $H(0)=0.$ Assume that $H$  has full rank at $0.$ Assume that
$a_{1},...,a_{28}$ are complex numbers such that
\begin{equation}
a_{1}A(H(\widetilde{z}))+...+a_{28}G(H(\widetilde{z}))=0~\text{for
all}~\widetilde{z} \in U.
\end{equation}
Then $a_{i}=0$ for all $1 \leq i \leq 28.$
\end{lemma}
{\it Proof of Lemma \ref{t3}}: Suppose not. Notice that $H$ has full
rank at $0.$ Hence there exist $26$ components $\widehat{H}$ of $H$
that give a local biholomorphism from $\mathbb{C}^{26}$ to
$\mathbb{C}^{26}.$ We assume these $26$ components $\widehat{H}$ are
the $H$  with $\eta$  dropped, where $\eta \in \{
\psi_{1},\psi_{2},\psi_{3},h_{10},...,h_{17},h_{20},...,h_{27},h_{30},...,h_{37}\}.$
By a local biholomorphic change of coordinates, we assume
$$\widehat{H}=(x_{1},x_{2},y_{0},...,y_{7},t_{0},...,t_{7},w_{0},...,w_{7}).$$
We still write the remaining components as $\eta.$

{\bf Case I:} If $\eta \in \{\psi_{1},\psi_{2},\psi_{3} \}$, without
loss of generality, we can assume $\eta= \psi_{3}.$  We first claim
that the coefficients of $A,B,$ i.e., $a_{1}, a_{2}$ are zero. This
is due to the fact that $t_{i}^2, w_{i}^2, 0 \leq i \leq 7$ can only
be canceled by $t_{i}\psi_{3}, w_{i}\psi_{3}$, which do not appear
in the expressions of $A(H),...,G(H).$ We then claim the
coefficients of $C$ are zero, for $x_{1}x_{2}$ can not be canceled.
Then the coefficients of all $D$'s have to be zero, for each
$t_{i}w_{j}$  is unique and  can not be canceled. Then it follows
trivially that all other coefficients are zero.

{\bf Case II:} If $\eta \in
\{h_{10},...,h_{17},h_{20},...,h_{27},h_{30},...,h_{37} \},$ without
loss of generality, we assume $\eta=h_{37}.$ Notice that the only
fact we are using about $\widehat{H}$ is that its components are
independent variables. For simplicity of notation, we will write
$$\widehat{H}=(x_{1},x_{2},x_{3},y_{0},...,y_{7},t_{0},...,t_{7},w_{0},...,w_{6}).$$
We  first claim that the coefficient of $A$ is zero. This is due to
the fact that $x_{2}x_{3}$ cannot be canceled. We also claim that
the coefficient of $B$ is zero. Suppose not. Notice that $t_{i}^2$
can  only be canceled by $t_{i}h_{37}.$ Then the coefficient of each
$D_{i}$ is non zero for  $0 \leq i \leq 7.$ Moreover, $x_{1}x_{3}$
can  only be canceled by $x_{1}h_{37}.$ This implies $h_{37}$ has a
linear $x_{3}$-term. Then, in particular, the $t_{7}h_{37}$ term in
$D_{0}$ will produce a $t_{7}x_{3}$ term. It cannot be canceled by
any other terms. This is a contradiction. Similarly, one can show
that the coefficient of $C$ is zero. Then we claim the coefficient
of $D_{0}$ is zero. Otherwise, to cancel the $x_{3}y_{0}$ term,
$h_{37}$ needs have a linear $x_{3}$ term. Then the term
$t_{7}h_{37}$ in $D_{0}$ will produce a $t_{7}x_{3}$ term, which
cannot be canceled by any other term. By the same argument, one can
show that the coefficients of all $D_{i}, 0 \leq i \leq 7$, are
zero. Similarly, we can obtain  the coefficients of all $E_{i}, 0
\leq i \leq 7$, are zero. Then we claim the coefficients of all
$F$'s have to be zero. This is because each $y_{i}t_{j}$  is unique.
It can not be canceled out. Finally we get the coefficient of $G$ to
be zero.
\
$\endpf$

This also establishes Proposition \ref{thm71}.
\
$\endpf$

\medskip

\begin{remark}
By Proposition \ref{thm71}, there exist multiindices
$\tilde{\beta}^1,...,\tilde{\beta}^{55}$ with $|\tilde{\beta}^j|
\leq 29,$ and there exist
$$z^0=(x^0_{1},x^0_{2},x^0_{3}, y^0_{0},...,y^0_{7},t^0_{0},..,t^0_{7},w^0_{0},...,w^0_{7}),~~x^0_{3} \neq 0,$$
such that
$$\left| \begin{matrix}
\frac{\partial^{|\beta^1|} (\psi_{1}(F))}{\partial \tilde{z}^{\tilde{\beta}^1}} & ...  & \frac{\partial^{|\beta^1|} (\psi_{55}(F))}{\partial \tilde{z}^{\tilde{\beta}^1}} \\
... & ... & ... \\
\frac{\partial^{|\beta^{55}|} (\psi_{1}(F))}{\partial
\tilde{z}^{\tilde{\beta}^{55}}} & ... &
\frac{\partial^{|\beta^{55}|} (\psi_{55}(F))}{\partial
\tilde{z}^{\tilde{\beta}^{55}}}
\end{matrix}\right|(z^0) \neq 0.$$
\end{remark}
Then we set
$\xi^0=(0,0,\xi^0_{3},0,...0,0,...,0,0,...,0),\xi^0_{3}=-\frac{1}{x^0_{3}}.$
It is easy to see that $(z^0, \xi^0) \in \mathcal{M}=\{\rho(z,
\xi)=0\}.$ Write
$$\mathcal{L}_{i}=\frac{\partial}{\partial x_{i}}- \frac{\frac{\partial \rho}{\partial x_{i}}(z,\xi)}{\frac{\partial \rho}{\partial x_{3}}(z,\xi)}\frac{\partial}{\partial x_{3}}, 1 \leq i \leq 2;$$
$$\mathcal{L}_{3+i}=\frac{\partial}{\partial y_{i}}- \frac{\frac{\partial \rho}{\partial y_{i}}(z,\xi)}{\frac{\partial \rho}{\partial x_{3}}(z,\xi)}\frac{\partial}{\partial x_{3}}, 0 \leq i \leq 7;$$
$$\mathcal{L}_{11+i}=\frac{\partial}{\partial t_{i}}- \frac{\frac{\partial \rho}{\partial t_{i}}(z,\xi)}{\frac{\partial \rho}{\partial x_{3}}(z,\xi)}\frac{\partial}{\partial x_{3}}, 0 \leq i \leq 7;$$
$$\mathcal{L}_{19+i}=\frac{\partial}{\partial w_{i}}- \frac{\frac{\partial \rho}{\partial w_{i}}(z,\xi)}{\frac{\partial \rho}{\partial x_{3}}(z,\xi)}\frac{\partial}{\partial x_{3}}, 0 \leq i \leq 7.$$

For any $26$-multiindex $\beta=(\beta_{1},...,\beta_{26}),$ we write
$\mathcal{L}^{\beta}
=\mathcal{L}_{1}^{\beta_{1}}...\mathcal{L}_{26}^{\beta_{26}}.$ Now
we define for any $55$ collection of $26-$multiindices
$\{\beta^1,...,\beta^{55}\},$
\begin{equation}
\Lambda(\beta^1,...,\beta^{55})(z,\xi):=\left|\begin{matrix}
 \mathcal{L}^{\beta^1}(\psi_{1}(F)) & ...  & \mathcal{L}^{\beta^1}(\psi_{55}(F)) \\
... & ...  & ...  \\
\mathcal{L}^{\beta^{55}}(\psi_{1}(F))  & ...  &
\mathcal{L}^{\beta^{55}}(\psi_{55}(F))
\end{matrix}\right|(z,\xi).
\end{equation}

Note that for any multiindex, $\mathcal{L}^{\beta}$ evaluating at
$(z^0, \xi^0)$ coincides with $\frac{\partial}{\partial
\tilde{Z}^{\beta}}.$ We have,
\begin{theorem}
There exists multiindices $\{{\beta}^1,...,{\beta}^{55}\},$ such
that
$$\Lambda({\beta}^1,...,{\beta}^{55})(z,\xi) \neq 0$$
for $(z,\xi)$ in a small neighborhood of $(z^0, \xi^0)$ and
$\beta^1=(0,...,0).$
\end{theorem}

\subsection{The exceptional class ${ M}_{16}$}

This case is very similar to  the hyperquadric setting. 
In this case, we write the coordinates of ${\mathbb C}^{16}$ as
$$z:=(x_{0},...,x_{7}, y_{0},...,y_{7}).$$

The defining function of the Segre family as described in
$\eqref{e16rho}$ is
$$\rho(z,\xi)=1+r_{z}\cdot r_{\xi}=1+\sum_{i=1}^{N}\psi_i(z)\psi_i(\xi), \ \ \ \text{ where } N=26 \text{ and }$$
\begin{equation}\label{eqt111}
r_{z}=(x_{0},...,x_{7},y_{0},...,y_{7},A_{0},...A_{7},B_{0},B_{1}).
\end{equation}
Here $A_{i}, 0 \leq i \leq 7,B_{0},B_{1}$ are homogeneous quadratic
polynomials in  $z.$ For instance,
$$B_{0}=\sum_{i=0}^7 x_{i}^2, B_{1}=\sum_{i=0}^7 y_{i}^2.$$
For the expressions for  $A_{i},$ see Appendix I.

Let $F$ be as before. We write
$$F=(f_{0},...,f_{7},\widetilde{f}_{0},...\widetilde{f}_{7}).$$
And define $r_{F}$ as the composition of $r_{z}$ with $F:$
\begin{equation}\label{eq99}
r_{F}=r_z\circ
F=(f_{0},...,f_{7},\widetilde{f}_{0},...\widetilde{f}_{7},A_{0}(F),...A_{7}(F),B_{0}(F),B_{1}(F)).
\end{equation}
Notice that $r_{F}$ has $26$ components. 

We will need the following lemma:
\begin{lemma}\label{lemma75}
For each fixed $\mu_{0},...,\mu_{6}$ with $(\sum_{i=0}^6
\mu_{i}^2)+1=0$ and fixed $(y_{0},...,y_{7})$ with
$(\sum_{i=0}^6\mu_{i}y_{i})+y_{7} \neq 0,$ we can always find
$(\xi_{0},...,\xi_{7})$ such that
$$1+y_{0}\xi_{0}+...+y_{7}\xi_{7}=0;~~
\sum_{i=0}^7 (\xi_{i})^2 =0,~~\xi_{j}=\mu_{j}\xi_{7}, 0 \leq j \leq
6,~~\xi_{7} \neq 0.$$
\end{lemma}
{\it Proof of Lemma \ref{lemma75}}: The proof is similar to that
 as in the hyperquadric case.
$\endpf$

Take the complex hyperplane $\mathbb{H}:
y_{7}+\sum_{j=0}^{6}\mu_{j}y_{j}=0$ in
$(x_{0},...,x_{7},y_{0},...,y_{7}) \in \mathbb{C}^{16}.$
Write
$L_{0}=\frac{\partial}{\partial x_{0}},...,L_{7}=\frac{\partial}{\partial x_{7}}; L_{8}=\frac{\partial}{\partial y_{0}}-\mu_{1}\frac{\partial}{\partial y_{7}},...,
L_{14}=\frac{\partial}{\partial
y_{6}}-\mu_{6}\frac{\partial}{\partial y_{7}}.$

Then $\{L_{i}\}_{i=0}^{14}$ forms a basis of the tangent vector
fields of $\mathbb{H}.$ For any multiindex
$\alpha=(\alpha_{0},..,\alpha_{14}),$ we write
 $L^{\alpha}=L_{0}^{\alpha_{0}}...L_{14}^{\alpha_{14}}.$ We define the notion of $L-$rank and  $L-$nondegeneracy as
 in Definition \ref{df21} using $r_{F}$ in (\ref{eq99}) and $L^{\alpha}$ instead of $\widetilde{z}^{\alpha}$.
 We write the $k$th $L$-rank defined in this setting as $\mathrm{rank}_{k}(r_{F},L)$.  We now need to prove
the following:
\begin{proposition}\label{thm76}
$F$ is $L-$nondegenerate near $0.$ More precisely,
$\mathrm{rank}_{11}(r_{F},L)=26.$
\end{proposition}
{\it Proof of Proposition \ref{thm76}}: Suppose not. As in the
hyperquadric case, by a modified version of Theorem \ref{thm211}, we
have that there exist $26$ holomorphic functions
$g_{0}(w),...,g_{25}(w)$ defined near $0$ on the $w-$plane
 with $\{ g_{0}(0),...,g_{25}(0)\}$ not all zero such that
the following holds for  $z \in U:$
\begin{equation}
\sum_{i=0}^{25}
g_{i}(y_{7}+\mu_{0}y_{0}+...+\mu_{6}y_{6})\psi_{i}(F(z)) \equiv 0.
\end{equation}

Then since $F$ has full rank at $0$, one can similarly prove   that
$g_{0}(0)=0,...,g_{15}(0)=0.$ Hence we obtain:
\begin{lemma}
There exist $c_{0},...,c_{9} \in \mathbb{C}$ that are not all zero
such that
\begin{equation}\label{eq95}
c_{0}A_{0}(F(Z))+...+c_{7}A_{7}(F(Z))+c_{8}B_{0}(F(Z))+c_{9}B_{1}(F(Z))
\equiv 0,
\end{equation}
for all $Z \in U$ when restricted on
$y_{7}+\sum_{i=0}^6\mu_{i}y_{i}=0.$
\end{lemma}

We then just need to show that (\ref{eq95}) can not hold by the
following lemma after
applying a linear change of coordinates.

\begin{lemma}\label{lemma68}
Let $H=(h_{0},...,h_{7},g_{0},...,g_{7})$ be a vector-valued
holomorphic function in a neighborhood $U$ of $0$ in $\tilde{z}=
 (x_{0},...,x_{7},y_{0},...,y_{6})\in \mathbb{C}^{15}$ with
$H(0)=0.$ Assume that $H$ has  full rank at $0.$ Assume that
$a_{0},...,a_{9}$ are complex numbers such that
\begin{equation}\label{eqna0a9}
a_{0}A_{1}(H(\widetilde{z}))+...+a_{7}A_{7}(H(\widetilde{z}))+a_{8}B_{0}(H(\widetilde{z}))
+ a_{9}B_{1}(H(\widetilde{z})) =0~\text{for all}~\widetilde{z} \in
U.
\end{equation}
Then $a_{i}=0$ for  $1 \leq i \leq 10.$
\end{lemma}
{\it Proof of Lemma \ref{lemma68}}: Suppose not.  Notice that $H$
has  full rank at $0.$ Hence there exist $15$ components
$\widehat{H}$ of $H$ that gives a local biholomorphism from
$\mathbb{C}^{15}$ to $\mathbb{C}^{15}.$ We assume these $15$
components $\widehat{H}$ are $H$ with $\eta$ being dropped, where
$\eta \in \{h_{0},...,h_{7},g_{0},...,g_{7} \}.$ By a local
biholomorphic change of coordinates, we assume
$\widehat{H}=(x_{0},...,x_{7},y_{0},...,y_{6}).$ We still write the
remaining component as $\eta.$ Without loss of generality, we assume
$\eta=g_{7}.$

First we claim the coefficient $a_{9}$ of $B_{1}$ is zero. Suppose
not. Note that $y_{1}^2, y_{2}^2$ can be only canceled by $g_{7}^2$.
Then $g_{7}$ will have linear $y_{1}, y_{2}$ terms. Hence $g_{7}^2$
will produce a $y_{1}y_{2}$ term. It cannot be canceled by any other
terms. This is a contradiction. Now we consider the coefficients of
$A_{0},...,A_{7}.$ We claim $a_{i}=0, 0 \leq i \leq 7.$ Suppose not.
We write
$$y_{7}(\widetilde{Z})=\lambda_{0}y_{0}+...+\lambda_{6}y_{6}+\mu_{0}x_{0}+...
+\mu_{7}x_{7}+O(2),$$ for some $\lambda_{i}, \mu_{j} \in \mathbb{C},
0 \leq i \leq 6, 0 \leq j \leq 7.$ By collecting the terms of the
form $x_{0}y_{i}$ in the Taylor expansion of (\ref{eqna0a9}) we get
\begin{equation}
a_{i}+a_{7}\lambda_{i}=0, 0 \leq i \leq 6.
\end{equation}
By collecting the terms of the form $x_{1}y_{i}, 0 \leq i \leq 6,$
we get,
$$
a_{1}+a_{3}\lambda_{0}=0, -a_{0}+a_{3}\lambda_{1}=0,
-a_{4}+a_{3}\lambda_{2}=0, -a_{7}+a_{3}\lambda_{3}=0,$$
$$a_{2}+a_{3}\lambda_{4}=0, -a_{6}+ a_{3}\lambda_{5}=0,
a_{5}+a_{3}\lambda_{6}=0.$$ By collecting the terms of the form
$x_{2}y_{i}, 0 \leq i \leq 6,$ we get,
$$a_{2}+a_{6}\lambda_{0}=0, a_{4}+a_{6}\lambda_{1}=0, -a_{0}+a_{6}\lambda_{2}=0,-a_{5}+a_{6}\lambda_{3}=0.$$
$$-a_{1}+a_{6}\lambda_{4}=0, a_{3}+a_{6}\lambda_{5}=0, -a_{7}+a_{6}\lambda_{6}=0.$$

One can further write down all the coefficients for $x_{i}y_{j}, 0
\leq i \leq 7, 0 \leq j \leq 6.$ Once this is done, one easily sees
that $a_{i} \neq 0$ for any $0 \leq i \leq 7.$ Otherwise, all
$a_{i}=0, \ 0 \leq i \leq 7.$

Then by the above equations, we see that the matrix
\begin{equation}
\left(
  \begin{array}{ccccccc}
    a_{0} & a_{1} & a_{2} & a_{3} & a_{4} & a_{5} & a_{6} \\
    a_{1} & -a_{0} & -a_{4} & -a_{7} & a_{2} & -a_{6} & a_{5} \\
    a_{2} & a_{4} & -a_{0} & -a_{5} & -a_{1} & a_{3} & -a_{7} \\
  \end{array}
\right)
\end{equation}
is of rank one. Then one can get a contradiction by, for instance,
carefully checking the determinants of its $2 \times 2$ submatrices.
Hence $a_{i}=0, 0 \leq i \leq 7.$ Finally we easily get the
coefficient $a_{8}$ of $B_{0}$ is zero.
\
$\endpf$

This then establishes Proposition \ref{thm76}.
\
$\endpf$

\begin{remark}
First fix $\mu_{0},...,\mu_{6}$ with $(\sum_{i=0}^6 \mu_{i}^2)+1=0.$
By Proposition \ref{thm76}, there exists multiindices
$\tilde{\beta}^1,...,\tilde{\beta}^{26}$ with  $|\tilde{\beta}^j|
\leq 11,$  and
$$Z^0=(x_{0}^0,...,x_{7}^0,y_{0}^0,...,y_{7}^0)~\text{with}~~ \sum_{i=0}^6\mu_{i}y_{i}+y_{7} \neq 0,$$
such that
$$\left| \begin{matrix}
 L^{\tilde{\beta}^1}(\psi_{1}(F)) & ...  & L^{\tilde{\beta}^1}(\psi_{26}(F))\\
... & ... & ... \\
L^{\tilde{\beta}^{26}}(\psi_{1}(F)) & ... &
L^{\tilde{\beta}^{26}}(\psi_{26}(F))
\end{matrix}\right|(Z^0) \neq 0.$$
\end{remark}

We then let $\xi^0=(0,...,0, \xi_{0}^0,...,\xi_{7}^0),$ where
$(\xi_{0}^0,...,\xi_{7}^0)$ is choosen as
 in Lemma \ref{lemma75} associated with $(y_{0}^0,...,y_{7}^0).$ That is
$$1+y_{0}^0\xi_{0}^0+...+y_{7}^0\xi_{7}^0=0;~~
\sum_{i=0}^7 (\xi_{i}^0)^2 =0,~~\xi_{j}^0=\mu_{j}\xi_{7}^0, 0 \leq j
\leq 6,~~\xi_{7}^0 \neq 0.$$ It is easy to see that $(z^0, \xi^0)
\in \mathcal{M}.$

 We now define
\begin{equation}
\mathcal{L}_{i}=\frac{\partial}{\partial x_{i}}-
\frac{\frac{\partial \rho}{\partial x_{i}}(z,\xi)}{\frac{\partial
\rho}{\partial y_{7}}(Z,\xi)}\frac{\partial}{\partial y_{7}}, 0 \leq
i \leq 7;
\end{equation}
\begin{equation}
\mathcal{L}_{8+i}=\frac{\partial}{\partial y_{i}}-
\frac{\frac{\partial \rho}{\partial y_{i}}(z,\xi)}{\frac{\partial
\rho}{\partial y_{7}}(Z,\xi)}\frac{\partial}{\partial y_{7}}, 0 \leq
i \leq 6;
\end{equation}
for $(z,\xi) \in \mathcal{M}$ near $(z^0, \xi^0).$ They are tangent
vector fields along $\mathcal{M}.$ Moreover, $\frac{\partial
\rho}{\partial y_{n}}(z,\xi)$ is nonzero near $(z^0, \xi^0).$

We define for any multiindex $\alpha=(\alpha_{0},..,\alpha_{14}),$
$\mathcal{L}^{\alpha}=\mathcal{L}_{0}^{\alpha_{0}}...\mathcal{L}_{14}^{\alpha_{14}}.$
Define for any $26$ collection of $15$-multiindices $\{\beta^1,...,
\beta^{26}\},$
\begin{equation}
\Lambda(\beta^1,...,\beta^{26})(z,\xi)=\left|
                                       \begin{array}{ccc}
                                         \mathcal{L}^{\beta^1}(\psi_{1}(F)) & ... & \mathcal{L}^{\beta^1}(\psi_{26}(F)) \\
                                         ... & ... & ... \\
                                         \mathcal{L}^{\beta^{26}}(\psi_{1}(F)) & ... & \mathcal{L}^{\beta^{26}}(\psi_{26}(F)) \\
                                       \end{array}
                                     \right|(z, \xi).
\end{equation}

By the fact that $\sum_{i=0}^7 (\xi_{i}^0)^2=0,$ one can check that,
for any multiindex $\alpha=(\alpha_{0},..,\alpha_{14}),$
$\mathcal{L}^{\alpha}=L^{\alpha}$ when evaluated at $(z^0, \xi^0).$
Then as before, we get the following:

\begin{theorem}
There exists multiindices $\{{\beta}^1,...,{\beta}^{26}\}$ such that
$$\Lambda({\beta}^1,...,{\beta}^{26})(z,\xi) \neq 0,$$
for $(z,\xi)$ in a small neighborhood of $(z^0, \xi^0)$ and
$\beta^1=(0,0,...,0).$
\end{theorem}

\bigskip\bigskip

\section{Appendix III:\ Transversality and flattening of  Segre families for the remaining cases}

In this appendix, we will complete the proof of Theorem
$\ref{flatten}$ for the remaining cases.

\medskip
{\it Continuation of the proof of Theorem \ref{flatten}}: By the same method used before, we first establish the second part of Theorem $\ref{flatten}$ by assuming the first part of Theorem $\ref{flatten}$ is true.  Namely, suppose $\xi^0\in\mathbb{C}^n\sm\{0\}$ and
$z^0$  and $z^1$ are smooth points on the Segre variety
$Q_{\xi^0}$ such that $Q_{z^0}$ and $Q_{z^1}$ are both smooth at $\xi^0$ and 
intersect transversally there.  We  shall prove that
there is a   biholomorphic parametrization
$\mathcal{G}(\tilde\xi_1,\tilde\xi_2,...,\tilde\xi_n)
=(\xi_1,\xi_2,...,\xi_n),$ with
$   (\tilde\xi_1,\tilde\xi_2,...,
\tilde\xi_n)\in U_1\times U_2\times...\times U_n\subset\mathbb C^n$. Here when $1 \leq j \leq 2, $  $U_j$ is a small neighborhood of $1\in {\mathbb C}$. When $3 \leq j \leq n,$ $U_j$ is a small neighborhood of $0 \in \mathbb{C}$  with
${\mathcal G}(1,1,0,\cdots,0)=\xi^0$,
such that $\mathcal{G}(\{\tilde\xi_1=1\}\times U_2\times...\times U_n)\subset Q_{z^0},
\mathcal{G}(U_1\times\{\tilde\xi_2=1\}\times U_3\times...\times U_n)\subset Q_{z^1},$
and $\mathcal{G}(\{\tilde\xi_1=t\}\times U_2\times...\times U_n),
\mathcal{G}(U_1\times\{\tilde\xi_2=s\}\times U_3\times...\times U_n),s\in U_1,t\in U_2$
are open pieces of  Segre varieties. Also,
$\mathcal{G}$  consists of algebraic functions with total degree
bounded by a constant depending only on $(M,\o)$.  By the first part of Theorem \ref{flatten}, we have 
\begin{displaymath}
{\rm rank\,}\left(\begin{matrix}\nabla \rho(z^0,\xi)|_{\xi^0}\\
\nabla \rho(z^1,\xi)|_{\xi^0}\end{matrix}\right)=2.
\end{displaymath}

Without loss of generality, we assume $\frac{\partial
	(\rho(z^0,\xi),\rho(z^1,\xi))}{\partial (\xi_{1},\xi_{2})}\neq 0$ at $\xi^0.$ Now we introduce new variables $\tilde\xi_{1},...,\tilde\xi_{n}$
and set up the system:
\begin{displaymath}
\begin{cases}P_{1}:\,\,&\rho(z^0,\tilde\xi_{1}\xi)=0\\

P_{2}:\,\,&\rho(z^1,\tilde\xi_{2}\xi)=0\\

P_{3}:
\,\,&\tilde\xi_{3}-\xi_{3}=0\\
\,...&...\\
P_{n}:
\,\,& \tilde\xi_{n}-\xi_{n}=0
\end{cases}
\end{displaymath}
Then
$\frac{\partial(P_{1},...,P_{n})}{\partial(\xi_{1},...,\xi_{n})}|_{A},\
\frac{\partial(P_{1},...,P_{n})}{\partial(\tilde\xi_{1},...,\tilde\xi_{n})}|_{A}\neq 0$, where
$A=(\tilde\xi_{1},...,\tilde\xi_{n},\xi_{1},...,\xi_{n})=(1,1,0,...,0,1,0,...,0)$.
By Lemma \ref{james-002}, we get the needed algebraic flattening with the bound total degree.  

Next, we proceed to prove the first part of Theorem \ref{flatten}.  It suffices to find a sufficiently close point $z^1$ to $z^0$ such that 
\begin{displaymath}
{\rm rank\,}\left(\begin{matrix}\nabla \rho(z^0,\xi)|_{\xi^0}\\
\nabla \rho(z^1,\xi)|_{\xi^0}\end{matrix}\right)=2.
\end{displaymath} We shall establish the above equation case by case as follows:

 {\bf Case
3. Symplectic Grassmannians}:  Pick $\xi_0=(1,0,0,...,0)$. The
defining equation of the Segre family is
$\rho=1+\sum_{i=1}^nz_{ii}\xi_{ii}+2\sum_{i<j}z_{ij}\xi_{ij}+2\sum_{2\leq
i<j}(z_{11}z_{ij}-z_{1j}z_{i1})
(\xi_{11}\xi_{ij}-\xi_{i1}\xi_{1j})+\sum_{i=2}^{n}(z_{11}z_{ii}-z_{1i}^2)
(\xi_{11}\xi_{ii}-\xi_{1i}^2)+\sum_{i<k,j<l,(i,j)\neq
(1,1)}(z_{ij}z_{kl}-z_{il}z_{kj})(\xi_{ij}\xi_{kl}-\xi_{il}\xi_{kj})
+{\rm high\,\,order\,\,terms},$  where $z_{ji}:=z_{ij}$ for $j>i.$

$Q_{\xi^0}=\{z|\rho(z,\xi^0)=1+z_{11}=0\},\nabla\rho(z,\xi^0)=(1,0,...,0).$
Hence  $Q_{\xi^0}$ is smooth, and for  $z\in Q_{\xi^0}$ we have
$z=(-1,z_{12},z_{22},z_{13},...,z_{(n-1)n}).$  Pick $z^0,z^1\in
Q_{\xi^0}$.  Then

$Q_{z^s}=\{\xi|0=\rho(z^s,\xi)=1+\sum_{i=1}^nz_{ii}^s\xi_{ii}+2\sum_{i<j}z_{ij}^s\xi_{ij}+2\sum_{2\leq
i<j}(z^s_{11}z_{ij}^s-z_{1j}^sz_{i1}^s)(\xi_{11}\xi_{ij}-\xi_{i1}\xi_{1j})
+\sum_{i=2}^{n}(z^s_{11}z^s_{ii}-(z^s_{1i})^2)
(\xi_{11}\xi_{ii}-\xi_{1i}^2)+
\sum_{i<k,j<l,(i,j)\neq
(1,1)}(z_{ij}^sz_{kl}^s-z_{il}^sz_{kj}^s)(\xi_{ij}\xi_{kl}-\xi_{il}\xi_{kj})+{\rm high\,\,order\,\,terms}\},$ for $s=0,1.$


\begin{displaymath}\left(\begin{matrix}\nabla \rho(z^0,\xi)|_{\xi^0}\\
\nabla \rho(z^1,\xi)|_{\xi^0}\end{matrix}\right)
=\left(\begin{matrix}\frac{\partial
\rho(z^0,\xi)}{\partial\xi_{11}}&\frac{\partial
\rho(z^0,\xi)}{\partial\xi_{12}}&...&\frac{\partial
\rho(z^0,\xi)}{\partial\xi_{1n}}&...
&\frac{\partial \rho(z^0,\xi)}{\partial\xi_{ij}}&...&\frac{\partial \rho(z^0,\xi)}{\partial\xi_{nn}}\\
\frac{\partial \rho(z^1,\xi)}{\partial\xi_{11}}&\frac{\partial \rho(z^1,\xi)}{\partial\xi_{12}}&...&\frac{\partial \rho(z^1,\xi)}{\partial\xi_{1n}}&...
&\frac{\partial \rho(z^1,\xi)}{\partial\xi_{ij}}&...&\frac{\partial \rho(z^1,\xi)}{\partial\xi_{nn}}\\
\end{matrix}\right)\big{|}_{\xi^0}
\end{displaymath}
\begin{displaymath}
=\left(\begin{matrix}-1&2z_{12}^0&2z_{13}^0&...&2z_{1n}^0&-(z_{12}^0)^2&-2z_{12}^0z_{13}^0&...&-(2-\delta_{ij})z_{1j}^0z_{1i}^0&...\\
-1&2z_{12}^1&2z_{13}^1&...&2z_{1n}^1&-(z_{12}^1)^2&-2z_{12}^1z_{13}^1&...&-(2-\delta_{ij})z_{1j}^1z_{1i}^1&...
\end{matrix}\right).
\end{displaymath}
Hence, we have
\begin{displaymath}
\rm{rank}\left(\begin{matrix}\nabla \rho(z^0,\xi)|_{\xi^0}\\
\nabla \rho(z^1,\xi)|_{\xi^0}\end{matrix}\right)
=\rm{rank}\left(\begin{matrix}-1&2z_{12}^0&2z_{13}^0&...&2z_{1n}^0&-(z_{12}^0)^2&-2z_{12}^0z_{13}^0&...&-(2-\delta_{ij})z_{1j}^0z_{1i}^0&...\\
-1&2z_{12}^1&2z_{13}^1&...&2z_{1n}^1&-(z_{12}^1)^2&-2z_{12}^1z_{13}^1&...&-(2-\delta_{ij})z_{1j}^1z_{1i}^1&...
\end{matrix}\right)
\end{displaymath}

\begin{displaymath}
=\rm{rank}
\left(\begin{matrix}-1&2z_{12}^0&2z_{13}^0&...&2z_{1n}^0&-(2-\delta_{ij})z_{1j}^0z_{1i}^0&...\\
0&2\Delta z_{12}^1&2\Delta z_{13}^1&...&2\Delta z_{1n}^1&(2-\delta_{ij})\{z_{1j}^1\Delta z_{1i}^1+\Delta z_{1j}^1z_{1i}^1-\Delta z_{1j}^1\Delta z_{1i}^1\}&...
\end{matrix}\right).\\
\end{displaymath}
where $\Delta z_{ij}^1=z_{ij}^1-z_{ij}^0.$  If we pick $z_{12}^1\neq
z_{12}^0$, then the above rank is $2.$

\medskip
{\bf Case 4. Orthogonal Grassmannians}:  Here we use the Pfaffian
embedding stated in $\S 2$.  Fixing
$\xi^0=(\xi_{12}^0,\xi_{13}^0,\xi_{23}^0,...,\xi_{(n-1)n}^0)=(1,0,...,0)$,
the defining function of the Segre family is given by
$\rho=1+\sum_{i<j}z_{ij}\xi_{ij}+\sum_{2<i<j}(z_{12}z_{ij}-z_{1i}z_{2j}+z_{1j}z_{2i})(\xi_{12}\xi_{ij}-\xi_{1i}\xi_{2j}+\xi_{1j}\xi_{2i})+\sum_{i<j<k<l,\{1,2\}\not\subset\{i,j,k,l\}}
(z_{ij}z_{kl}-z_{ik}z_{jl}+z_{il}z_{jk})(\xi_{ij}\xi_{kl}-\xi_{ik}\xi_{jl}+\xi_{il}\xi_{jk})+{\rm
\,high\,order\,terms}.$ Note here we use the notation $z_{ji}:=-z_{ij}$ for $j>i.$



Note $Q_{\xi^0}=\{z|0=\rho(z,\xi^0)=1+z_{12}\}.$ Hence it is smooth.
Since $z\in Q_{\xi^0}$, we have $z=(-1,z_{13},...,z_{(n-1)n}).\quad$
Pick $z^0, \ z^1\in Q_{\xi^0}$. Then

$
Q_{z^s}=\{\xi|0=\rho(z^s,\xi)=1+\sum_{i<j}z_{ij}^s\xi_{ij}+\sum_{2<i<j}(z_{12}^sz_{ij}^s-z_{1i}^sz_{2j}^s+z_{1j}^sz_{2i}^s)(\xi_{12}\xi_{ij}-\xi_{1i}\xi_{2j}+\xi_{1j}\xi_{2i})\\
+\sum_{i<j<k<l,\{1,2\}\not\subset\{i,j,k,l\}}(z_{ij}^sz_{kl}^s-z_{ik}^sz_{jl}^s+z_{il}^sz_{jk}^s)(\xi_{ij}\xi_{kl}-\xi_{ik}\xi_{jl}+\xi_{il}\xi_{jk})+\,h.\,o.\,t.s.\}
$, for $s=0,1.$


\begin{displaymath}\left(\begin{matrix}
\nabla \rho(z^0,\xi)|_{\xi^0}\\
\nabla \rho(z^1,\xi)|_{\xi^0}
\end{matrix}\right)
=\left(\begin{matrix}\frac{\partial \rho(z^0,\xi)}{\partial\xi_{12}}&\frac{\partial \rho(z^0,\xi)}{\partial\xi_{13}}&...&\frac{\partial \rho(z^0,\xi)}{\partial\xi_{1n}}&...
&\frac{\partial \rho(z^0,\xi)}{\partial\xi_{ij}}&...&\frac{\partial \rho(z^0,\xi)}{\partial\xi_{(n-1)n}}\\
\frac{\partial \rho(z^1,\xi)}{\partial\xi_{12}}&\frac{\partial
\rho(z^1,\xi)}{\partial\xi_{13}}&...&\frac{\partial
\rho(z^1,\xi)}{\partial\xi_{1n}}&...
&\frac{\partial \rho(z^1,\xi)}{\partial\xi_{ij}}&...&\frac{\partial \rho(z^1,\xi)}{\partial\xi_{(n-1)n}}\\
\end{matrix}\right)\big{|}_{\xi^0}
\end{displaymath}
\begin{displaymath}
=\left(\begin{matrix}-1&z_{13}^0&...&z_{1n}^0&...&z_{2n}^0&(-z_{13}^0z_{24}^0+z_{14}^0z_{23}^0)a&...&(-z_{1i}^0z_{2j}^0+z_{1j}^0z_{2i}^0)a&...\\

-1&z_{13}^1&...&z_{1n}^1&...&z_{2n}^1&(-z_{13}^1z_{24}^1+z_{14}^1z_{23}^1)a&...&(-z_{1i}^1z_{2j}^1+z_{1j}^1z_{2i}^1)a&...
\end{matrix}\right).
\end{displaymath}
Hence,
\begin{displaymath}\rm{rank}\left(\begin{matrix}
\nabla \rho(z^0,\xi)|_{\xi^0}\\
\nabla \rho(z^1,\xi)|_{\xi^0}
\end{matrix}\right)
=\rm{rank}\left(\begin{matrix}-1&z_{13}^0&...&z_{1n}^0&...&z_{2n}^0&...\\

0&\Delta z_{13}^1&...&\Delta z_{1n}^1&...&\Delta z_{2n}^1&...
\end{matrix}\right).
\end{displaymath}
Here $\Delta z_{ij}^1=z_{ij}^1-z_{ij}^0.$  If we choose
$z_{13}^1\neq z_{13}^0,$  then the rank is $2.$ 

\medskip
{\bf Case 5. $M_{16}$}:  Pick
$\xi^0=(\kappa_0^0,\kappa_1^0,...,\kappa_7^0,\eta_0^0,\eta_1^0,...,
\eta_7^0)=(1,0,...,0),\ z^0\in Q_{\xi^0}.$  The defining equation of
the Segre family is

$1+x_0\kappa_0+x_1\kappa_1+...+x_7\kappa_7+y_0\eta_0+y_1\eta_1+...+y_7\eta_7+(x_0y_0+x_1y_1+...)(\kappa_0\eta_0+\kappa_1\eta_1+...)+(-y_0x_1+y_1x_0+...)(-\eta_0\kappa_1+\eta_1\kappa_0+...)+...+(x_0^2+x_1^2+...+x_7^2)({\kappa_0}^2+{\kappa_1}^2...+{\kappa_7}^2)+
(y_0^2+y_1^2+...+y_7^2)(\eta_0^2+\eta_1^2+...+\eta_7^2)=0.$

$Q_{\xi^0}=\{z|\rho(z,\xi^0)=1+x_0+(x_0^2+x_1^2+...+x_7^2)=0\},$ and
$\nabla\rho(z,\xi^0)|_{z^0}=(1+2x_0,2x_1,...,2x_7^0,0,...,0).$ Hence
$Q_{\xi^0}$ is smooth. Pick $z^0,z^1\in Q_{\xi^0}$.  Then

$Q_{z^s}=\{\xi|0=\rho(z^s,\xi)=1+x_0^s\kappa_0+x_1^s\kappa_1+...+x_7^s\kappa_7+y_0^s\eta_0+y_1^s\eta_1+...+y_7^s\eta_7+(x_0^sy_0^s+x_1^sy_1^s+...)(\kappa_0\eta_0+\kappa_1\eta_1+...)+(-y_0^sx_1^s+y_1^sx_0^s+...)(-\eta_0\kappa_1+\eta_1\kappa_0+...)+...+((x_0^s)^2+(x_1^s)^2+...+(x_7^s)^2)({\kappa_0}^2+{\kappa_1}^2+...+{\kappa_7}^2)+
((y_0^s)^2+(y_1^s)^2+...+(y_7^s)^2)(\eta_0^2+\eta_1^2+...+\eta_7^2)\}$, for $s=0,1.$


$$\rm{rank}\left(\begin{matrix}\nabla \rho(z^0,\xi)|_{\xi^0}\\
\nabla \rho(z^1,\xi)|_{\xi^0}
\end{matrix}\right)
\geq\rm{rank}\left(\begin{matrix}\frac{\partial \rho(z^0,\xi)}{\partial \kappa_0}&\frac{\partial \rho(z^0,\xi)}{\partial \kappa_1}&...&\frac{\partial \rho(z^0,\xi)}{\partial \kappa_7}\\
\frac{\partial \rho(z^1,\xi)}{\partial \kappa_0}&\frac{\partial \rho(z^1,\xi)}{\partial \kappa_1}&...&\frac{\partial \rho(z^1,\xi)}{\partial y_7}\\
\end{matrix}
\right)\big{|}_{\xi^0}
$$

$$
 =\rm{rank}
\left(\begin{matrix}-2-x_0^0&x_1^0&x_2^0&\cdots&x_7^0\\
-2-x_0^1&x_1^1&x_2^1&\cdots&x_7^1
\end{matrix}\right).\eqno{(C)}
$$

Since $(-2-x_0^0, x_1^0, x_2^0, \cdots, x_7^0)\neq (0,...,0),$ we
can pick $z^1$ sufficiently close to $z^0$, such that the above rank is
2.  That is because  $Q_{\xi^0}$ is irreducible  and the
subvarieties, defined by $2\times 2$ minors of the last matrix in
$(C)$,
are thin subsets of $Q_{\xi^0}$. 

\medskip
{\bf Case 6. $ M_{27}$}:   Take
$\xi^0=(\xi_1^0,\xi_2^0,\xi_3^0,\eta_0^0,\eta_1^0,...,\eta_7^0,\kappa_0^0,\kappa_1^0,...,\kappa_7^0,\tau_0^0,\tau_1^0,...,\tau_7^0)=(1,0,...,0).$
The defining function of the  Segre family is
$1+r_{z}\cdot{r}_{\xi}$ where \begin{displaymath}
r_{z}=(x_{1},x_{2},x_{3},y_{0},...,y_{7},z_{0},...,z_{7},w_{0},...,w_{7},A,B,C,D_{0},...D_{7},
E_{0},...,E_{7},F_{0},...,F_{7},G)\end{displaymath}
\begin{displaymath}
r_{\xi}=(\xi_1,\xi_2,\xi_3,...,\eta_7,...,\kappa_7,...,\tau_7,{A}(\xi),{B}(\xi),{C}(\xi),...,{D}_{7}(\xi),
...,{E}_{7}(\xi),...,{G}(\xi)).
\end{displaymath}
Here $A,B,C, D_{i},E_{i},F_{i}$ are homogeneous quadratic polynomials;
\ $G$ is a homogeneous cubic polynomial  defined in Appendix I.

For our purpose here, we present terms only involving $\xi_1,\xi_2,$
and omit those involving $\xi_3,\eta_0,\eta_1,$
$...,\eta_7,\kappa_0,\kappa_1,...,\kappa_7,\tau_0,\tau_1,...,\tau_7$
as follows: $\rho(z,\xi)=
1+x_1\xi_1+x_2\xi_2+...+(x_1x_2-(\sum_{i=0}^7y_i^2))(\xi_1\xi_2-(\sum_{i=0}^7(\tau_i)^2))+\cdots.$

$Q_{\xi^0}=\{z|\rho(z,\xi^0)=1+x_1=0\},\nabla\rho(z,\xi^0)=(1,0,0,...,0).$
Hence $Q_{\xi^0}$ is smooth and for $z\in Q_{\xi^0}$, we have
$z=(-1,x_2,x_3,...,)$.  Pick $z^0,z^1\in Q_{\xi^0}.$  Then

$Q_{z^s}=\{\xi|0=\rho(z^s,\xi)=1+x_1^s\xi_1+x_2^s\xi_2+...+(x_1^sx_2^s-(\sum_{i=0}^7(y_i^s)^2))(\xi_1\xi_2-(\sum_{i=0}^7(\tau_i)^2))+...\},$ for $s=0,1.$

\begin{displaymath}\rm{rank}\left(\begin{matrix}\nabla \rho(z^0,\xi)|_{\xi^0}\\
\nabla \rho(z^1,\xi)|_{\xi^0}
\end{matrix}\right)
=\rm{rank}\left(\begin{matrix}\frac{\partial \rho(z^0,\xi)}{\partial\xi_{1}}&\frac{\partial \rho(z^0,\xi)}{\partial\xi_{2}}&\frac{\partial \rho(z^0,\xi)}{\partial\xi_{3}}&...
&\frac{\partial \rho(z^0,\xi)}{\partial\eta_{7}}&...&\frac{\partial \rho(z^0,\xi)}{\partial\kappa_{7}}&...&\frac{\partial \rho(z^0,\xi)}{\partial\tau_7}\\
\frac{\partial \rho(z^1,\xi)}{\partial\xi_{1}}&\frac{\partial \rho(z^1,\xi)}{\partial\xi_{2}}&\frac{\partial \rho(z^1,\xi)}{\partial\xi_{3}}&...
&\frac{\partial \rho(z^1,\xi)}{\partial\eta_{7}}&...&\frac{\partial \rho(z^1,\xi)}{\partial\kappa_{7}}&...&\frac{\partial \rho(z^1,\xi)}{\partial\tau_7}\\
\end{matrix}\right)\big{|}_{\xi^0}
\end{displaymath}
\begin{displaymath}{}
\geq\rm{rank}\left(\begin{matrix}\frac{\partial \rho(z^0,\xi)}{\partial\xi_{1}}&\frac{\partial \rho(z^0,\xi)}{\partial\xi_{2}}\\
\frac{\partial \rho(z^1,\xi)}{\partial\xi_{1}}&\frac{\partial
\rho(z^1,\xi)}{\partial\xi_{2}}
\end{matrix}\right)\big{|}_{\xi^0}
=\rm{rank}\left(\begin{matrix}
-1&-(\sum_{i=0}^7(y_i^0)^2)\\
-1&-(\sum_{i=0}^7(y_i^1)^2)
\end{matrix}\right)\big{|}_{\xi^0}\geq 2,
\end{displaymath}
for those $z^1$'s  such that $\sum_{i=0}^7(y_i^1)^2\neq
\sum_{i=0}^7(y_i^0)^2.$ This can be done in any small neighborhood
of $z^0;$ for  $\{z|\sum_{i=0}^7(y_i)^2=B\}$ is a thin set in
$\{z|0=1+x_1\}$ for each fixed $B\in\mathbb{C}$.

 This
completes the proof of the flattening theorem. $\endpf$


\bigskip
\noindent H. Fang, Department of Mathematics, University of Wisconsin-Madison, Madison, WI 53706,
USA. (hfang35@wisc.edu)

\noindent X. Huang, Department of Mathematics, Rutgers University,
New Brunswick, NJ 08903, USA. (huangx$@$math.rutgers.edu)

\noindent M. Xiao, Department of Mathematics,
       University of California, San Diego,
9500 Gilman Dr. La Jolla, CA 92093, USA. (m3xiao$@$ucsd.edu)
\end{document}